\documentclass{article}

\usepackage{arxiv}

\usepackage[utf8]{inputenc} % allow utf-8 input
\usepackage[T1]{fontenc}    % use 8-bit T1 fonts
\usepackage{hyperref}       % hyperlinks
\usepackage{url}            % simple URL typesetting
\usepackage{booktabs}       % professional-quality tables
\usepackage{amsfonts}       % blackboard math symbols
\usepackage{nicefrac}       % compact symbols for 1/2, etc.
\usepackage{microtype}      % microtypography
\usepackage{doi}
\usepackage{enumitem}
\usepackage{style_primale}
\usepackage{graphicx}
\usepackage{tikz}
\newsavebox{\imagebox}

\title{A 3D-1D Virtual Element Method for Modeling Root Water Uptake}

\author{ \href{https://orcid.org/0000-0001-8642-4258}{\includegraphics[scale=0.06]{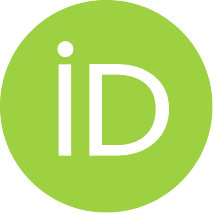}\hspace{1mm}Stefano~Berrone} \\
	Dipartimento di Scienze Matematiche\\
	``G. L. Lagrange''\\
	Politecnico di Torino, TO, 10129 \\
	\texttt{stefano.berrone@polito.it} \\
	\And
	\href{https://orcid.org/0000-0001-8544-6199}{\includegraphics[scale=0.06]{orcid.pdf}\hspace{1mm}Stefano~Ferraris} \\
	Dipartimento Interateneo di Scienze, \\
	Progetto e Politiche del Territorio\\
	Politecnico di Torino e Università di Torino, TO \\
	\texttt{stefano.ferraris@unito.it} \\
	\And
	\href{https://orcid.org/0000-0001-8544-6199}{\includegraphics[scale=0.06]{orcid.pdf}\hspace{1mm}Denise~Grappein} \\
	Dipartimento di Scienze Matematiche\\
	``G. L. Lagrange''\\
	Politecnico di Torino, TO, 10129 \\
	\texttt{denise.grappein@polito.it} \\
	\And
	\href{https://orcid.org/0000-0002-8540-3639}{\includegraphics[scale=0.06]{orcid.pdf}\hspace{1mm}Gioana~Teora} \\
	Dipartimento di Scienze Matematiche\\
	``G. L. Lagrange''\\
	Politecnico di Torino, TO, 10129 \\
	\texttt{gioana.teora@polito.it} \\
	\And
	\href{https://orcid.org/0000-0001-7123-9199}{\includegraphics[scale=0.06]{orcid.pdf}\hspace{1mm}Fabio~Vicini} \\
	Dipartimento di Scienze Matematiche\\
	``G. L. Lagrange''\\
	Politecnico di Torino, TO, 10129 \\
	\texttt{fabio.vicini@polito.it} \\
}

\renewcommand{\shorttitle}{3D-1D VEM for Modeling Root Water Uptake}

\newcommand{\Scoll}{\partial\Sigma_{\mathrm{collar}}}
\newcommand{\Stip}{\partial\Sigma_{\mathrm{tip}}^j}
\newcommand{\Lcoll}{\partial\Lambda_{\mathrm{collar}}}
\newcommand{\Ltip}{\partial\Lambda_{\mathrm{tip}}^j}
\hypersetup{
pdftitle={A 3D-1D Virtual Element Method for Modeling Root Water Uptake},
pdfsubject={q-bio.NC, q-bio.QM},
pdfauthor={Stefano~Berrone, Stefano~Ferraris, Denise~Grappein, Gioana~Teora, Fabio~Vicini},
pdfkeywords={Optimization, Soil-Roots, Coupling, Virtual Element Method},
}

\begin{document}
\maketitle

\begin{abstract}
An optimization-based strategy is proposed for coupling three-dimensional and one-dimensional problems (3D-1D coupling) in the context of soil-root interaction simulations. This strategy, originally designed to tackle generic 3D-1D coupled problems with discontinuous solutions, is here extended to the case of non-linear problems and applied, for the first time, along with a virtual element discretization of the 3D soil sample. This further enhances the capability of the method to handle geometrical complexities, allowing to easily mesh domains characterized, for instance, by the presence of stones and other impervious obstacles of arbitrary shape. A discrete-hybrid tip-tracking strategy is adopted to model both the root growth and the evolution in time of the water flux, the pressure head and the water content, both in the roots and in the surrounding soil sample. By choosing proper rules for the generation of branches, realistic root-network configurations are obtained. Several numerical examples are proposed, proving both the accuracy of the adopted method and its applicability in realistic and large scale simulations.
\end{abstract}

% keywords can be removed
\keywords{Optimization \and Soil-Roots \and Coupling \and Virtual Element Method}

\section{Introduction}\label{sec:introduction}

In recent years, the study of the interactions between root system architecture (RSA) and soil has gained significant interest, in particular in the framework of actions confronting climate change and drought \cite{Vereecken2022, Vetterlein2022, SchnepfCarminati2022, SchnepfLeitner2022, VanderborghtLeitner2024, VanderborghtCouvreur2024}. Investigating root activities and architecture during growth in situ can be very expensive and challenging, as the involved processes are often very slow and RSA usually exhibits considerable variability even within the same species and under identical environmental conditions. Therefore, modeling root water potential in various 3D soil geometries is crucial for guiding the design of experimental measurements \cite{Ganesan2024, Jin2020}.

In this paper, we aim to model water flow in an unsaturated soil and within the root xylem, which is the plant vascular tissue responsible for the distribution of water and nutrients taken up from the soil. The root water uptake occurs while the RSA expands over time, influenced by soil variables such as pressure head and the water content.

Numerical models describing root growth in the soil matrix include both continuum and discrete models \cite{Jin2020}. In continuum models, variables such as root length density and root apical meristems are defined using probability density functions that are updated through diffusion equations to model root growth. In contrast, discrete models explicitly describe the components of the RSA, representing the root system as a network of cylindrical inclusions immersed in the soil domain. 
To reduce the computational burden of simulations, such inclusions are often identified with their centerlines, hence avoiding the definition of a 3D mesh in such thin and elongated domains.
%In this case, physical quantities defined in the plant-root are either assumed to vary only in the longitudinal direction, or are averaged across the cross-section of each cylinder \cite{Koch2018, Schnepf2020}, leading to a 3D-1D coupled problem.}
Deriving a well-posed formulation for 3D-1D coupled problems is however non-trivial, as a bounded trace operator is not defined in standard Sobolev spaces between manifolds having codimension bigger than 1. Different approaches have been investigated in literature in order to tackle this issue, including the introduction of suitable weighted Sobolev spaces \cite{Dangelo2012_siam,Dangelo2012}, and the use of regularization \cite{Tornberg2004,Heltai2019,Koch2020} or splitting \cite{Gjerde2018a, Gjerde2019} techniques to treat the singular terms. In \cite{Laurino2019} proper averaging operators were introduced to perform a geometrical model reduction from a 3D-3D to a 3D-1D problem and in \cite{Kuchta2021} the coupling was tackled by means of Lagrange multipliers in a domain-decomposition setting.

In this work, we follow the approach proposed in \cite{Grappein2023}, where the 3D-3D problem is reduced to a 3D-1D problem by assuming that the quantities of interest have a negligible variation on the cross-sections of the inclusions, and proper subspaces of Sobolev spaces are introduced to reflect the main model assumptions. The resulting problem is then solved in a PDE-constrained optimization framework: auxiliary variables are introduced at the interface between soil and roots to fully decouple the two problems, and a quadratic cost functional, expressing the error committed in the fulfillment of interface conditions, is minimized, constrained by the set of 3D-1D equations. The resulting scheme has no mesh conformity requirements, making the method particularly flexible and suitable for the simulation on an evolving root network. Moreover, it provides a direct computation of interface variables, allowing us to obtain a straightforward quantification of water uptake, without any need of post-processing.

The 3D pressure head is discretized using the Virtual Element Method (VEM), whereas a mixed Finite Element Formulation with continuous pressure is adopted for the xylem variables. The use of VEM for the discretization of the soil sample, and hence of a polytopal grid, allows us to easily mesh complex geometries which can arise in the presence of obstacles, physical barriers, or layered and fractured soil \cite{Manzini2004, Morandage2021}, while facilitating the description of irregularly shaped geometries like plant pots or locally refined grids \cite{Koch2018, Fassino2024}. For what concerns the xylem, the mixed Finite Element formulation is derived by imposing strongly flux conservation at root branching points, sharply limiting the overall number of degrees of freedom. The root growth is described by means of a discrete-hybrid tip-tracking strategy, modeling the evolution of the root tip position. Rules for the creation of new branches are adapted from the ones available in literature (see among others \cite{Pages1989, Moraes2018, Koch2018, Jin2020}), to obtain a realistic representation of the RSA. In order to identify barriers we define an additional scalar field, which is constant in time but which takes its maximum value on the surface of obstacles. By introducing information on this field in the computation of the growth velocity we easily avoid roots from penetrating barriers.

The outline of the paper is as follows. In Section \ref{sec:3D3Dmodel} we present the 3D-3D coupled problem, while in Section \ref{sec:optimization}, after performing the model reduction, we introduce the PDE-constrained optimization problem in the case of a single root segment. In Section \ref{sec:RSA}, we present the Root System Architecture and describe the root growth process, considering multiple segments and branching. In Section \ref{sec:discretization}, we detail the time and space discretizations and describe our solving strategy to handle the arising non-linear problem through the conjugate gradient scheme. Finally, in Section \ref{sec:numericalresults}, we propose some numerical experiments to validate and show the potentiality of our procedure.

\subsection{Notation}\label{sec:not}

Let us consider a generic open subset $\genericset \subset \R^d$, $d=1,2,3$. Given two scalar functions $p,q \in \leb{2}{\genericset}$, two vector fields $\bm{a},\ \vv \in \vleb{2}{\genericset}{d}$ and two tensor fields $\bm{T},\ \bm{\sigma} \in \vleb{2}{\genericset}{d \times d}$, we denote by
\begin{align*}
\scal[\genericset]{p}{q} = \int_{\genericset} p q \ ~d\omega,\quad\scal[\genericset]{\bm{a}}{\vv} = \int_{\omega} \bm{a}\cdot \vv\ ~d\omega,\quad \scal[\genericset]{\bm{T}}{\bm{\sigma}} = \int_{\genericset} \bm{T} : \bm{\sigma} \ ~d\omega, 
\end{align*}
where $\bm{T} : \bm{\sigma} := \sum_{i,j = 1}^n \bm{T}_{ij} \bm{\sigma}_{ij}$. Furthermore, given a generic Sobolev space $\mathcal{H}$, we use the symbol $\norm[\mathcal{H}]{}$ to indicate the norm in $\mathcal{H}$, $\mathcal{H}'$ to denote the dual space of $\mathcal{H}$ and the symbol $\langle \cdot, \cdot \rangle_{\mathcal{H}, \mathcal{H}'}$ to denote the duality pairing between $\mathcal{H}$ and $\mathcal{H}'$.

Let us consider a finite time interval $(0,\tfin]$, i.e. $\tfin < + \infty$ and its discretization in $J$ sub-intervals $\Ij = (t_{j-1}, t_{j}]$, where the instants $t_j$ are such that
\begin{equation*}
0 = t_0 < t_1 < \dots < t_J = \tfin.
\end{equation*}

\begin{figure}[t]
\centering
\includegraphics[width=0.5\textwidth]{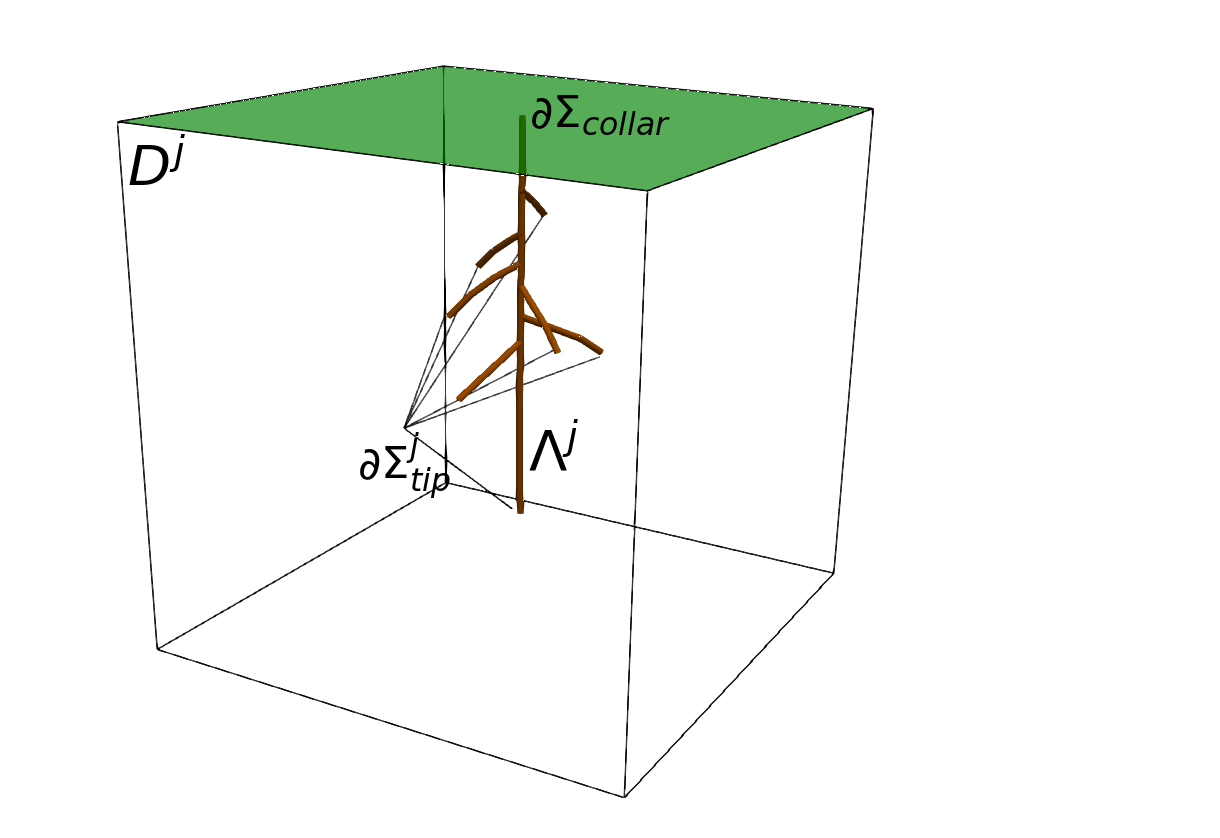}%
\includegraphics[width=0.4\textwidth]{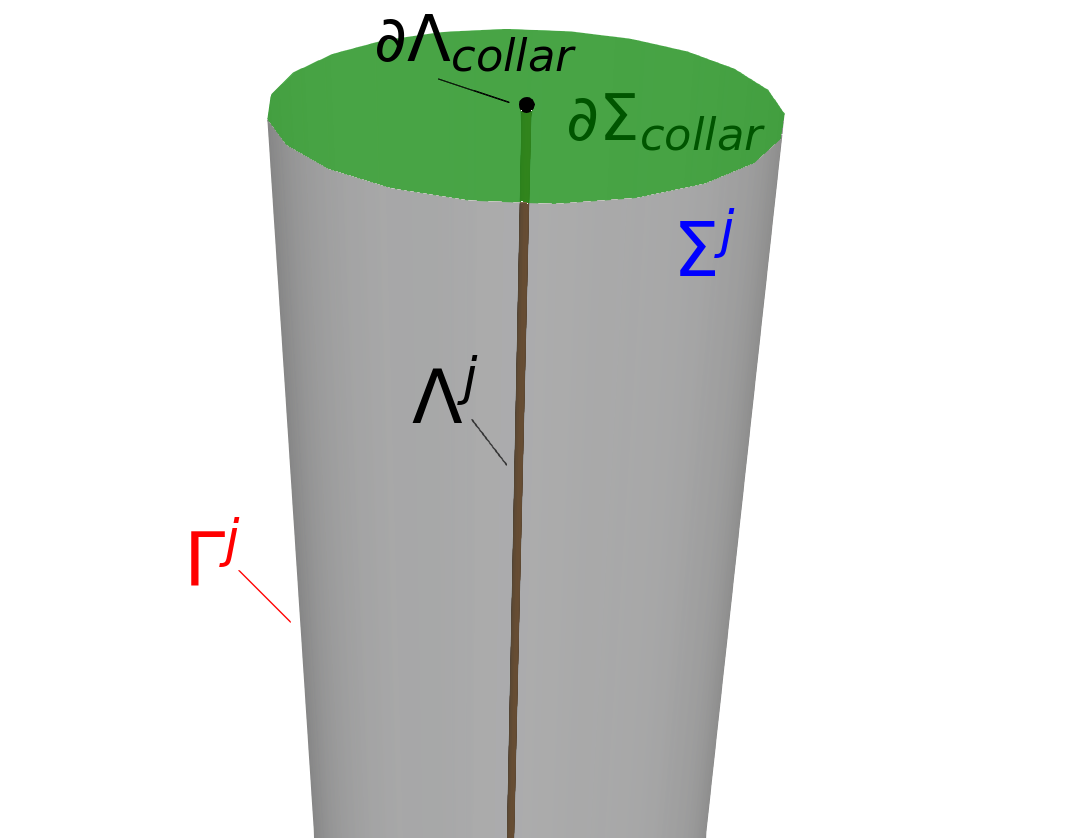}
\caption{Details on notation (zoom-in on the right).}
\label{fig:domain_scheme} 
\end{figure}

Let $\Omega \subset \R^3$ be a polyhedral convex domain representing the soil sample, and %and in which is embedded a generalized cylinder $\Sigmaj[ ](t) \in \R^3$ variable in time.
let $\Sigmaj[ ](t_j) = \Sigmaj \subset \Omega$ denote the root system architecture which is considered fixed while the quantities of interest vary during the time interval $\Ij$. In the following, we assume that $\Sigmaj$ is composed of thin tubular vessels of constant radius $R$, and that $R$ is much smaller than the root length and the characteristic dimension of the domain $\Omega$. Further we assume that
\begin{equation*}
\Sigmaj[0] \subseteq \Sigmaj[1] \subseteq \dots \subseteq \Sigmaj[J],
\end{equation*}
which means that we allow for the growth of the root network, but not for its regression or remodeling. The boundaries of $\Omega$ and $\Sigmaj$ are denoted respectively by $\partial \Omega$ and $\partial \Sigmaj$. In particular,
	$$\partial\Sigmaj=\Gammaj\cup \Stip\cup\Scoll,$$ where $\Gammaj$ denotes the lateral surface of the overall RSA, whereas $\Scoll$ and $\Stip$ refer to the collar cross-section and to the union of the apical cross-sections (see Figure \ref{fig:domain_scheme}). The root collar is defined as the interface between the stem and the roots and it is assumed to lie on the boundary of the soil-sample, i.e. $\Scoll\subset \partial \Omega$.  The portion of soil not including the cylinder is denoted by $\Dj = \Omega \setminus \overline{\Sigmaj}$. We define $\partial \mathcal{D}\subset \partial \Omega$ as
		$$\partial \mathcal{D}=\partial \Omega \setminus \Scoll.$$ Let us remark that $\partial \mathcal{D}$ and $\Scoll$ do not vary in time, hence we do not add the superscript $j$ to these symbols.
			Finally, we denote by $\Lambda^j$ the centerline of $\Sigmaj$, and we denote by $\Lcoll$ the center of the collar cross-section and by $\Ltip$ the union of the centers of the apical cross-sections, such that $\overline{\Lambdaj}=\Lambdaj\cup \Lcoll\cup\Ltip$.

\section{The 3D-3D model}\label{sec:3D3Dmodel}
The process of water uptake during the growth of the RSA is
here described by means of two partial differential equations, for the modeling of water flow in the soil sample and in the xylem, and an ordinary differential equation, modeling the evolution of the position of root tips. 

In the present section, we assume, for clarity of exposition, that $\Sigmaj$ is made by a single straight cylindrical root segment. Details on how multiple segments are handled and about root growth and branching will be provided in Section \ref{sec:RSA}. 

In what follows, the subscript $\varsigma$ is used to refer to soil variables, while the subscript $\chi$ denotes xylem variables.

\subsection{Water flow in soil}
The soil water flux $\fluxws{}$ is described by Darcy law. This equation, combined with mass conservation, leads to Richards equation, which models the flow of a wetting fluid (the water) in the unsaturated soil, under the effect of gravity and capillary action, and in the presence of a nonwetting fluid (the air) which is supposed to be at constant pressure equal to $0$. 
%Given an initial condition
%\begin{equation}
%\rph{} = \ri, \quad \text{ at } t_0, \text{ in } \Dj[0],
%\end{equation}
The Richards equation at each time $t \in \Ij$, for $j=1,\dots,J$, can be stated in the so-called head-based form as the non-linear, time-dependent, possibly degenerate, parabolic problem \cite{Celia1990}
\begin{equation}
\begin{cases}
\fluxws{}  = - \rK(\rph{})\nabla[\rph{} + z] & \text{ in } \Dj,\\
\rc({\rph{}})\frac{\partial \rph{}}{\partial t} + \div \fluxws{} = 0 & \text{ in } \Dj, \\
\rph{} = 0 & \text{ on } \partial \mathcal{D}.
\end{cases}
\label{eq:watersoil}
\end{equation}
The symbol $\rph{}$ denotes the soil pressure head, which is positive in the saturated zone, where only water is present, and negative in the unsaturated zone. The symbol
$z$ denotes the vertical coordinate with the $z$-axis direction $\ee_z$ oriented against the gravity direction, $\rK(\rph{}) = \rK(\wc(\rph{}))$ is the hydraulic conductivity, $\wc(\rph{})$ is the volumetric soil-water content and $\rc(\rph{}) = \rc( \wc(\rph{})) = \frac{d \wc}{d \rph{}}$ is the specific soil moisture capacity.
In $\Ij[1] = (t_0, t_1]$, we define an initial condition in $t_0$
\begin{equation}\label{eq:initial_cond}
\rph{} = \ri, \quad \text{ in } \Dj[0],
\end{equation}
while, for $\Ij = (t_{j-1}, t_j]$ with $j \geq 2$, the pressure head at time $t_{j-1}$ is directly provided by the value obtained at the end of the previous time interval.

Here, for the sake of clarity, we assume homogeneous Dirichlet boundary conditions for the soil pressure head, but other boundary conditions can be considered for Problem \eqref{eq:watersoil} to take into account natural phenomena like rainfall, snow melt, runoff and soil evaporation, as well as to prescribe fluxes along seepage faces or drainage boundaries.

The soil water content $\wc$ and, consequently, the hydraulic conductivity $\rK$ and the capacity term $\rc$ are modeled as a function of the soil pressure head $\rph{}$ using different characteristic curves which are usually determined experimentally for the different soils, as the van Genuchten's model \cite{VanGenuchten1980}.
\begin{comment}
An example is the van Genuchten's model \cite{VanGenuchten1980}:
\begin{equation}
\wc(\rph{}) = \begin{cases}
\wc_r + \frac{\wc_s - \wc_r}{[1 + (a \vert \rph{} \vert)^n]^m} & \text{if }  \rph{} < 0,\\
\wc_s & \text{if }  \rph{} \geq 0,
\end{cases}\quad \text{with } m = 1 - \frac{1}{n},
\end{equation}
\begin{equation}
\rK(\rph{}) = \begin{cases}
\rK_s  \frac{\left(1 - (a \vert \rph{}\vert)^{n-1}[1 + (a \vert \rph{} \vert)^n]^{-m}\right)^2}{[1 + (a \vert \rph{} \vert)^n]^{\frac{m}{2}}} & \text{if } \rph{} < 0,\\
\rK_s & \text{if }  \rph{} \geq 0,
\end{cases}
\end{equation}
\begin{equation}
\rc(\rph{}) = \begin{cases}
-nm(a\vert \rph{}\vert)^n \frac{\wc_s - \wc_r}{\rph{}[1 + (a \vert \rph{} \vert)^n]^{m+1}}  & \text{if } \rph{} < 0,\\
0 & \text{if }  \rph{} \geq 0,
\end{cases}
\end{equation}
where $a$ and $n$ are fitting parameters, $\wc_r$ and $\wc_s$ are the residual and the saturated volumetric water content, respectively, and $\rK_s $ is the saturated hydraulic conductivity. 
\end{comment}
Let us remark that the Richards equation degenerates into an elliptic equation if $\rc(\rph{}) =0$ at some points of the domain. On the other hand, if $ \rK(\rph{}) = 0$, then the equation becomes a first-order ordinary differential equation. 
Non-linearities and degeneracy pose significant challenges when dealing with this kind of problem  \cite{Mitra2021}. Here, we consider the non-degenerate case and we assume that data are sufficiently regular.

\subsection{Water flow in roots}

%Inside the plant root, the water flows through the root xylem. Here,
Water transport in the xylem is modeled via the following Stokes equation for incompressible Newtonian fluids \cite{Rand1996, Raimondi2017} with a reactive term due to friction forces:
\begin{equation}
\begin{cases}
- \div \xst + \ka \fluxwx{} = - \rho g \ee_z&  \text{ in } \Sigmaj,\\
\div \fluxwx{} = 0& \text{ in } \Sigmaj,\\
\xst \nnx = \bm{0} & \text{ on } \Scoll \cup\Stip.
\end{cases}
\label{eq:navierstokes}
\end{equation}
Here, $\fluxwx{}$ is the flow velocity in the xylem, $\ka$ is the resistance coefficient due to friction, $\rho$ is the fluid density, and $g$ is the gravity acceleration.
The tensor $\xst$ is the Cauchy pseudo-stress tensor, i.e.
\begin{equation}
\xst = - \xp \bm{I} +  \overline{\nu}  \nabla \fluxwx{},
\label{eq:pseudo_stress_tensor}
\end{equation}
where $\xp$ is the water pressure, $\overline{\nu}$ is the dynamic viscosity, and $\I \in \R^{3 \times 3}$ is the identity tensor.
We assume homogeneous Neumann boundary conditions, but non-homogeneous conditions can be considered to express, for instance, a non-zero transpiration rate at the root collar \cite{Mai2019}.

\subsection{The coupled problem}
The two problems \eqref{eq:watersoil} and \eqref{eq:navierstokes} are coupled by imposing the following interface conditions at $\Gammaj$ (\cite{Cattaneo2014, Carro2024}), namely
\begin{align}
&\label{eq:transmissivity} \fluxws{} \cdot \nns = -\fluxwx{} \cdot \nnx = \Lp (\rph{} - \xph{})\quad \text{on } \Gammaj,\\
&\fluxwx{} \cdot \ttau^{d} = \bm{0},\quad d=1,2, \quad \text{on } \Gammaj,
\end{align}
where $\xph{} = \frac{\xp}{\rho g}$ is the xylem pressure head, $\Lp$ expresses the conductivity of the root wall, $\nns$ is the unit outward normal vector to $\Gammaj$ directed from soil to root and $\nnx = -\nns$, while $\ttau^d$, $d=1,2$, are two unit tangential vectors to $\Gammaj$, such that
\begin{equation}
\nnx \nnx^T + \ttau^1 (\ttau^1)^T + \ttau^2 (\ttau^2)^T = \I.
\end{equation}
We observe that, during the daytime, the difference $\rph{} - \xph{}$ is usually positive, resulting in a water flow from soil to root \cite{Koch2018}.

Summarizing,
the 3D-3D problem for modeling the water exchange between soil and root reads as
\begin{align}
&\fluxws{}  = - \rK(\rph{})\nabla[\rph{} + z] & \text{in } \Dj,\label{eq:3Dproblem_soil_a}\\
&\rc({\rph{}})\frac{\partial \rph{}}{\partial t} + \div \fluxws{} = 0 & \text{in } \Dj ,\label{eq:3Dproblem_soil_b}\\
&\rph{} = 0 & \text{on }\partial \mathcal{D},\label{eq:3Dproblem_soil_c} \\
&\fluxws{} \cdot \nns = \Lp (\rph{} - \xph{}) & \text{on } \Gammaj, \label{eq:3Dproblem_soil_d} \\[1em]
&- \div \xst + \ka \fluxwx{} = - \rho g \ee_z& \text{in } \Sigmaj\label{eq:3Dproblem_root_a},\\
&\div \fluxwx{} = 0&  \text{in } \Sigmaj,\label{eq:3Dproblem_root_b}\\
&\xst \nnx = \bm{0} & \text{on }  \Scoll \cup\Stip, \label{eq:3Dproblem_root_c}\\
&\fluxwx{} \cdot \nnx  = \Lp (\xph{} - \rph{}) & \text{on } \Gammaj,\label{eq:3Dproblem_root_d} \\
&\fluxwx{} \cdot \ttau^d = 0,\ d=1,2 & \text{on } \Gammaj, \label{eq:3Dproblem_root_e}
\end{align} equipped with initial conditions for the soil pressure head $\rph{}$ on the initial geometry $\mathcal{D}^0$.

\section{The model reduction and the PDE-constrained optimization problem}\label{sec:optimization}
In order to perform the geometric model reduction leading to the 3D-1D problem, we now make some fundamental assumptions. 
Let us consider a cylindrical coordinate system around the axis of $\Sigmaj$, for all $j =0,\dots,J$. We denote by $\ee_r,\ \ee_{\beta}$, and $\ee_{s}$ the radial, the angular and the axial unit vectors, respectively, whereas $(r, \beta, s)$ are the related coordinates with $r \in [0,R),\ \beta \in [0,2\pi), s \in [0, \Sj]$.
Let us then assume that (see also \cite{Formaggia2003, Lawal2010}):
\begin{enumerate}[label=\textbf{A.\arabic*}]
\item The root segment is rigid and rectilinear. Thus, we neglect the effect of pressure on root volume.
\item Each cross-section is assumed to be circular with a constant radius $R$. In particular, $R$ does not depend on the axial coordinates $s$. 
\item \label{ass:pressure}We assume that all quantities are independent of the angular coordinate $\beta$ and, in particular, that the xylem pressure head $\xph{}$ is constant on each section, that is at each time $t \in \Ij$, $j=1,\dots,J$,
\begin{equation*}
\xph{}(t,r,\beta,s) = \xph{a}(t,s),\quad  \forall r \in [0,R),\ \beta \in [0,2 \pi),\ s \in [0,\Sj],
\end{equation*}
where $\xph{a}$ varies only with time and along the axial coordinate $s$.
\item \label{ass:velocity}The rotational velocity is negligible with respect to the axial velocity, i.e. in cylindrical coordinates, at each time $t \in \Ij$, $\forall j=1,\dots,J$, we can write 
\begin{equation}
\fluxwx{} = \begin{bmatrix}
(\fluxwx{})_r\\
(\fluxwx{})_{\beta}\\
(\fluxwx{})_{s}
\end{bmatrix} = \begin{bmatrix}
0\\
0\\
(\fluxwx{})_{s}
\end{bmatrix}, \quad \forall r \in [0, R),\ \beta \in [0, 2\pi),\ s \in [0, \Sj].
\label{eq:axial_vel}
\end{equation}
We further assume that the axial component of the xylem velocity can be defined as $(\fluxwx{})_{s}(t,r, \beta, s) = (\fluxwx{p})_{s}(t, s)f(rR^{-1})$, where  $(\fluxwx{p})_{s}(t,s)$ is the mean axial velocity on the cross-section $\Sigmaj(s)$ at the time $t$, i.e.
\begin{equation*}
(\fluxwx{p})_{s}(t, s) = \frac{1}{\vert \Sigmaj(s) \vert} \int_{\Sigmaj(s)} (\fluxwx{})_{s}\quad \forall s \in [0,\Sj],
\end{equation*}
and $f: [0,1] \to \R$, $f \in \con{\infty}{[0,1]}$ is a radial velocity profile which satisfies
\begin{equation*}
\int_0^1 f(y)y ~dy = \frac{1}{2}.
\end{equation*}
An example of such function $f$ which can be found in \cite{Formaggia2003} is given by
\begin{equation}
f(y) = \frac{\gamma + 2}{\gamma}(1-y^{\gamma})\quad \text{ for some $\gamma\in \R$.}
\label{eq:velocityprofile}
\end{equation}
We observe that, for $\gamma = 2$, the function \eqref{eq:velocityprofile} reduces to the well-known Hagen-Poiseuille solution in a pipe of uniform (circular) cross-section.
\end{enumerate}
Under Assumption \ref{ass:velocity}, at each time $t \in \Ij$, Equation \eqref{eq:3Dproblem_root_b}, i.e. the mass conservation equation for the xylem, reduces to
\begin{equation}
\frac{\partial (\fluxwx{})_s}{\partial s} = 0\quad  \text{ in } \Sigmaj.
\label{eq:null_divergence}
\end{equation}

Let us define, as in \cite{Grappein2022, Grappein2023}, the following two uniform extension operators
\begin{equation*}
\es : \sob{1}{\Lambdaj} \to \sob{1}{\Sigmaj} \text{ and } \eg : \sob{1}{\Lambdaj} \to \sob{1/2}{\Gammaj} 
\end{equation*}
defined such that, given $\qsx[a] \in \sob{1}{\Lambdaj}$, $\forall s \in [0,\Sj]$, $\es\qsx[a]$ is the uniform extension of the pointwise value $\qsx[a](s)$ to the cross-section $\Sigmaj(s)$ and $\eg\qsx[a]$ is the uniform extension of $\qsx[a](s)$ to the boundary $\Gammaj(s)$ of the same cross-section, i.e.
\begin{equation*}
\left(\es \qsx[a]\right)(\xx) = \qsx[a](s)\ \forall \xx \in \Sigmaj(s)\quad  \text{ and } \quad \left(\eg \qsx[a] \right)(\xx) = \qsx[a](s)\ \forall \xx \in \Gammaj(s).
\end{equation*}
Finally, we consider the trace operator
\begin{equation*}
\tg : \sob{1}{\Sigmaj} \to \sob{1/2}{\Gammaj} \text{ such that }  \tg \qsx = \qsx_{|\Gammaj} \ \forall \qsx \in \sob{1}{\Sigmaj}.
\end{equation*}
Now, we can introduce the following spaces
\begin{gather*}
\Vrxj = \sob{1}{\Lambdaj}, \quad \Vxj = \left\{\vvx{} = \begin{bmatrix}
0 \\ 
0 \\
(\vvx{})_s
\end{bmatrix} \in [\sob{1}{\Sigmaj}]^3:\ (\vvx{})_s = (\vvx{p})_s f(rR^{-1})  \text{ with } (\vvx{p})_s \in \Vrxj \right\},\\
 \Qrj = \sob[0]{1}{\Lambdaj},\quad \Qxj = \{\qx{} \in \sob{1}{\Sigmaj}:\ \qx{} = \es \qx{a}  ,\ \qx{a} \in \Qrj \},
\end{gather*}
and write the variational formulation of the root problem as: \textit{Find $\fluxwx{} \in\Vxj$, $\xph{} \in \Qxj$ such that}
\begin{align}
\addlinespace &\scal[\Sigmaj]{\viscr \nabla \fluxwx{}}{ \nabla \vvx{}} + \scal[\Sigmaj]{\kar \fluxwx{}}{\vvx{}} - \scal[\Sigmaj]{\xph{}}{ \div \vvx{}} = -  \scal[\Sigmaj]{ \ee_z}{\vvx{}}  & \forall \vvx{} \in \Vxj,\label{eq:3Dvarproblem_root_mom}\\
\addlinespace &\scal[\Sigmaj]{\div \fluxwx{}}{\qx{}} = 0 & \forall \qx{} \in \Qxj,\label{eq:3Dvarproblem_root_mass}
\end{align}
where we denote by $\viscr = \frac{\overline{\nu}}{\rho g}$ and $\kar = \frac{\ka}{\rho g}$.

Thanks to the above regularity assumptions, problem \eqref{eq:3Dvarproblem_root_mom}-\eqref{eq:3Dvarproblem_root_mass} can be easily reduced to a 1D problem following an approach similar to the one presented in \cite{Grappein2023}. Let us start by considering Equation \eqref{eq:3Dvarproblem_root_mass}: integration by parts and the interface condition \eqref{eq:3Dproblem_root_d} yield
\begin{align}
0 &= \scal[\Sigmaj]{\div \fluxwx{}}{\qx{}} = - \scal[\Sigmaj]{ \fluxwx{}}{\nabla \qx{}} + \scal[\partial \Sigmaj]{ \fluxwx{} \cdot \nn_{\xylem}}{\qx{}}\nonumber \\
&= - \scal[\Sigmaj]{ \fluxwx{}}{\nabla \qx{}} + \scal[\Stip]{\fluxwx{} \cdot \nn_{\xylem}}{ \qx{}}+ \scal[\Scoll]{\fluxwx{} \cdot \nn_{\xylem}}{ \qx{}}+ \scal[\Gammaj]{\Lp (\tg \xph{} - \tg \rph{})} {\tg \qx{}} \label{eq:byparts}
\end{align}
Let us first of all observe that, according to Assumption \ref{ass:velocity}, $\forall \qx{} \in \Qxj$ and $ \qx{a} \in \Qrj: \tg \qx{} = \es \qx{a}$
\begin{align}
	\scal[\Sigmaj]{ \fluxwx{}}{\nabla \qx{}}&=\int_0^{S^j}\int_{0}^{2 \pi} \int_{0}^R(\fluxwx{p})_{s}f(rR^{-1}) \frac{\partial \qx{a}}{\partial s}  r~dr~d\beta~ds\nonumber \\&= 2 \pi R^2 \Big(\int_{0}^1  f(y) y ~dy\Big)\Big(\int_0^{S^j} (\fluxwx{p})_{s}\frac{\partial \qx{a}}{\partial s}
	ds \Big) =\scal[\Lambdaj]{\pi R^2 (\fluxwx{p})_{s}}{\frac{\partial \qx{a}}{\partial s}}.
 \label{eq:u_nablaq}
\end{align}
Similarly, denoting by $\Sigmaj(\bar{s})$ a generic transversal section obtained by cutting $\Sigmaj$ at $\bar{s}\in[0,S_j]$ with a plane orthogonal to $\Lambdaj$ we have 
\begin{equation}
	\scal[\Sigmaj(\bar{s})]{\fluxwx{} \cdot \ee_{s}}{ \qx{}} =\pi R^2 (\fluxwx{p})_s(\bar{s})\ \qx{a}(\bar{s}),  \quad\forall \qx{} \in \Qxj,~ \qx{a} \in \Qrj: \tg \qx{} = \es \qx{a}.
 \label{eq:uq_sec}
\end{equation}
Furthermore, according to Assumption \ref{ass:pressure},
\begin{align}
\scal[\Gammaj]{\tg\xph{}}{\tg \qx{}}=\scal[\Lambdaj]{2\pi R\xph{a}}{\qx{a}}, \quad\forall \qx{} \in \Qxj,~ \qx{a} \in \Qrj: \tg \qx{} = \es \qx{a}. 
\label{eq:psixylem}
\end{align}
Finally, introducing the spaces
\begin{align*}
&\Hc = \{\phi \in \sob{1/2}{\Gammaj}:\ \phi = \eg \hat{\phi},\ \hat{\phi} \in \Qrj\}, \quad\Qsj = \{\qs{} \in \sob{1}{\Dj[ ]}:\ \tg {\qs{}}  \in \Hc \text{ and }  \qs{} = 0 \text{ on } \partial \Dj[ ]\}
\end{align*}
and assuming that $\rph{}\in \Qsj$, we similarly obtain 
\begin{align}
\scal[\Gammaj]{\tg\rph{}}{\tg \qx{}}=\scal[\Lambdaj]{2\pi R\rph{a}}{\qx{a}}, \quad\forall \qx{} \in \Qxj,~ \qx{a} \in \Qrj: \tg \qx{} = \es \qx{a} 
\label{eq:psisoil}
\end{align}
where $\rph{a}\in \Qrj:~\eg \rph{a}=\tg\rph{}$. Going now back to \eqref{eq:byparts} and exploiting \eqref{eq:u_nablaq}-\eqref{eq:psisoil}, along with the fundamental theorem of calculus, we obtain
\begin{align*}
0&=-\scal[\Lambdaj]{\pi R^2 (\fluxwx{p})_{s}}{\frac{\partial \qx{a}}{\partial s}} +\pi R^2 (\fluxwx{p})_s(S^j)\ \qx{a}(S^j) -\pi R^2 (\fluxwx{p})_s(0)\ \qx{a}(0)+\scal[\Lambdaj]{\Lp(\xph{a}-\rph{a})}{\qx{a}}\\
&= - \int_{0}^{\Sj} \pi R^2 (\fluxwx{p})_s \frac{\partial \qx{a}}{\partial s} ds + \int_{0}^{\Sj} \pi R^2\frac{\partial ((\fluxwx{p})_s\ \qx{a})}{\partial s} ds + \int_{0}^{S^j} 2 \pi R \Lp (\xph{a} - \rph{a}) \qx{a}ds\\
&= \scal[\Lambdaj] {\pi R^2 \frac{\partial ((\fluxwx{p})_s)}{\partial s} }{\qx{a}} +  \scal[\Lambdaj] { 2 \pi R \Lp (\xph{a} - \rph{a})}{\qx{a}}, \quad\forall \qx{} \in \Qxj,~ \qx{a} \in \Qrj: \tg \qx{} = \es \qx{a}.
\end{align*}

Concerning Equation \eqref{eq:3Dvarproblem_root_mom}, we have that
\begin{equation*}
\scal[\Sigmaj]{\xph{}}{ \div \vvx{}} =  \int_{\Lambdaj} \pi R^2 \frac{\partial (\vvx{p})_s}{\partial s} \xph{a} ,\quad
\scal[\Sigmaj]{ \ee_z}{\vvx{}} = \int_{\Lambdaj} \pi R^2 \ee_z \cdot  \ee_s (\vvx{p})_s. 
\end{equation*}
Moreover, due to \eqref{eq:null_divergence} and the fact that, for any $\vvx{}\in\Vxj$ 
\begin{equation*}
\frac{\partial (\vvx{})_s}{\partial r} = R^{-1} f'(rR^{-1}) (\vvx{p})_s,
\end{equation*}
we obtain
\begin{equation}
\nabla \fluxwx{}=\begin{bmatrix}
0 &0 & R^{-1} f'(rR^{-1}) (\fluxwx{p})_s \\
0 &0 &0\\
0 &0 &0
\end{bmatrix}, \quad \nabla \vvx{} = \begin{bmatrix}
0 &0 &R^{-1} f'(rR^{-1}) (\vvx{p})_s\\
0 &0 &0\\
0 &0 &\frac{\partial (\vvx{})_s}{\partial s} 
\end{bmatrix},
\label{eq:nabla_velocity}
\end{equation}
and hence
\begin{align*}
\scal[\Sigmaj]{\viscr \nabla \fluxwx{}}{ \nabla \vvx{}} + \scal[\Sigmaj]{\kar \fluxwx{}}{\vvx{}} &= \int_{\Lambdaj} \viscr  (\fluxwx{p})_s (\vvx{p})_s F' + \int_{\Lambdaj} \kar (\fluxwx{p})_s (\vvx{p})_s F\\
&= \scal[\Lambdaj]{\kx (\fluxwx{p})_s}{(\vvx{p})_s},
\end{align*}
with $F, F' \in \R$ defined as
\begin{align*}
	F := 2 \pi R^2 \int_0^1 [f(y)]^2 y ~dy, \quad F' := 2 \pi \int_{0}^1 [f'(y)]^2 y ~dy,
\end{align*}
and $\kx = F \kar + F' \viscr$.

Concerning the equations for the soil volume, we proceed as in \cite{Grappein2023} by rewriting the soil problem in the primal form. We hence look for $\rph{}\in \Qsj{}$ such that
\begin{align*}
\scal[\Dj]{\rc({\rph{}})\frac{\partial \rph{}}{\partial t}}{\qs{}}  + \scal[\Dj]{\rK(\rph{}) \nabla \rph{}}{\nabla \qs{}} + \scal[\Gamma]{\Lp(\tg \rph{} - \tg\xph{})}{\tg \qs{}} +  \scal[\Dj]{\rK(\rph{}) \ee_z}{\nabla \qs{}} =0\quad \forall \qs{} \in \Qsj.
\end{align*}
As for the mass balance equation, the integral over the lateral surface $\Gammaj$ can be reduced to an integral over the centerline $\Lambdaj$ as
\begin{align*}
\scal[\Gamma]{\Lp(\tg \rph{} - \tg\xph{})}{\tg \qs{}} =  \scal[\Lambdaj]{ 2 \pi R \Lp(\rph{a} - \xph{a}) }{\qs{a}}.
\end{align*}

Finally, we obtain the following 3D-1D coupled problem:  \textit{At each time $t \in \Ij$, for $j=1,\dots,J$, find $(\rph{},\ux, \xph{a}) \in~\Qsj \times \Vrxj \times \Qrj$ such that}
\begin{equation}
\resizebox{0.91\hsize}{!}{$
\begin{cases}
\scal[\Dj]{\rc({\rph{}})\frac{\partial \rph{}}{\partial t}}{\qs{}}  + \scal[\Dj]{\rK(\rph{}) \nabla \rph{}}{\nabla \qs{}} + \scal[\Lambdaj]{ 2 \pi R \Lp(\rph{a} - \xph{a}) }{\qs{a}}+ \scal[\Dj]{\rK(\rph{}) \ee_z}{\nabla \qs{}}=0&\forall \qs{} \in \Qsj,\\
\\
\addlinespace \scal[\Lambdaj]{ \kx \ux }{\vx} - \scal[\Lambdaj]{ \xph{a} }{\pi R^2 \frac{\partial \vx}{\partial s}} = - \scal[\Lambdaj]{ \pi R^2 \ee_z \cdot  \ee_s}{\vx} & \forall \vx \in \Vrxj,\\
\addlinespace \scal[\Lambdaj]{\pi R^2 \frac{\partial \ux}{\partial s}}{ \qx{a}} + \scal[\Lambdaj]{ 2\pi R \Lp (\xph{a} - \rph{a})}{\qx{a}} = 0& \forall \qx{a} \in \Qrj,\\
\end{cases}$}
\label{eq:coupled_3D1Dproblem}
\end{equation}
where, for the ease of the notation, we set $\ux := (\fluxwx{p})_s$, $\vx := (\vvx{p})_s$ .
We highlight that the reduction procedure results in a Darcy-like problem in mixed form for the xylem sap.

\subsection{PDE-constrained optimization problem}
Instead of solving the set of coupled equations summarized in Problem \eqref{eq:coupled_3D1Dproblem}, we here recast the problem into a PDE-constrained optimization problem, following the steps in \cite{Grappein2023}. First of all, adopting a domain decomposition approach, we introduce two auxiliary variables $\phi_{\soil},\phi_{\xylem} \in \Hc$ such that
\begin{align}
&\sdual[\Hc]{\tg\rph{} -\phi_\soil}{\eta} = 0\ &&\forall \eta \in {\Hc}'\label{eq:constr1},\\
&\sdual[\Hc]{\tg\xph{} - \phi_\xylem}{\eta} = 0 &&\forall \eta \in {\Hc}'\label{eq:constr2}.
\end{align}
Denoting by $\lambdas, \lambdax \in \Qrj$ the one dimensional variables such that $\eg\lambdas =\tg \phi_\soil$ and $\eg\lambdax =\tg \phi_\xylem$, Problem \eqref{eq:coupled_3D1Dproblem} can be rewritten as:
\textit{at each time $t \in \Ij$, for $j=1,\dots,J$, find $(\rph{},\ux, \xph{a}, \lambdas, \lambdax) \in \Qsj \times \Vrxj \times \Qrj  \times \Qrj \times \Qrj$ such that}
\begin{align}
&\scal[\Dj]{\rc({\rph{}})\frac{\partial \rph{}}{\partial t}}{\qs{}}+\scal[\Dj]{\rK(\rph{}) \nabla \rph{}}{\nabla \qs{}} + \scal[\Lambdaj]{2 \pi R \Lp( \rph{a} - \lambdax)}{\qs{a}}\nonumber  \\[-0.5em]&\hspace{8cm}+\scal[\Dj]{\rK(\rph{}) \ee_z}{\nabla \qs{}} = 0 && \forall \qs{} \in \Qsj\label{eq:3D1Dproblem1}, \\
&\scal[\Lambdaj]{ \kx \ux }{\vx} - \scal[\Lambdaj]{ \xph{a} }{  \pi R^2  \frac{\partial \vx}{\partial s}} = - \scal[\Lambdaj]{  \pi R^2 \ee_z \cdot  \ee_s}{\vx}  && \forall \vx \in \Vrxj\label{eq:3D1Dproblem2}, \\
&\scal[\Lambdaj]{ \pi R^2  \frac{\partial \ux}{\partial s}}{ \qx{a}} + \scal[\Lambdaj]{ 2\pi R \Lp (\xph{a} - \lambdas)}{\qx{a}} = 0 &&  \forall \qx{a} \in \Qrj\label{eq:3D1Dproblem3}, \\[0.7em]
&\sdual[\Qrj]{\rph{a} - \lambdas}{\etah} = 0\ && \forall \etah \in {\Qrj}'\label{eq:constr1D_1},\\
&\sdual[\Qrj]{\xph{a} - \lambdax}{\etah} = 0 && \forall \etah \in \Qrj'\label{eq:constr1D_2}.
\end{align}
Conditions \eqref{eq:constr1D_1}-\eqref{eq:constr1D_2}, i.e. \eqref{eq:constr1}-\eqref{eq:constr2}, are replaced by the minimization of a cost functional, measuring the error committed in the fulfillment of such constraints. Let
\begin{equation}
\mathcal{J}( \lambdas, \lambdax) = \frac{1}{2} \left( \norm[\Hc]{\tg \rph{}\!( \lambdax) - \eg\lambdas}^2 + \norm[\Hc]{\tg \xph{}\!( \lambdas) - \eg\lambdax}^2 \right).
\label{eq:costfunctional}
\end{equation}

The final optimal control problem reads as: 
\begin{equation}
\begin{gathered}
\textit{At each time } t \in \Ij, \textit{ for } j=1,\dots,J, \textit{ find } (\rph{},\ux, \xph{a}, \lambdas, \lambdax) \in \Qsj \times \Vrxj \times \Qrj  \times \Qrj \times \Qrj \textit{ such that} \\
 \textit{the functional \eqref{eq:costfunctional} is minimized subject to \eqref{eq:3D1Dproblem1}-\eqref{eq:3D1Dproblem3}.}
\end{gathered}
\label{eq:optimization_problem}
\end{equation}
There are two main advantages of solving the problem in a PDE-constrained optimization framework. First, there is great flexibility in the choice of the meshes: no mesh conformity is required between the 3D soil and the 1D roots and the different one-dimensional variables can be meshed independently \cite{Grappein2023}. Second, the introduction of the auxiliary variables allows the direct computation of interface variables without any need of post-processing. This can be interesting, for example, to easily monitor the flux exchanged at the interface as the root network grows.

\section{Modeling root growth}\label{sec:RSA}

\begin{figure}
\centering
\begin{subfigure}{0.45\textwidth}
\centering
    \includegraphics[width=1\textwidth]{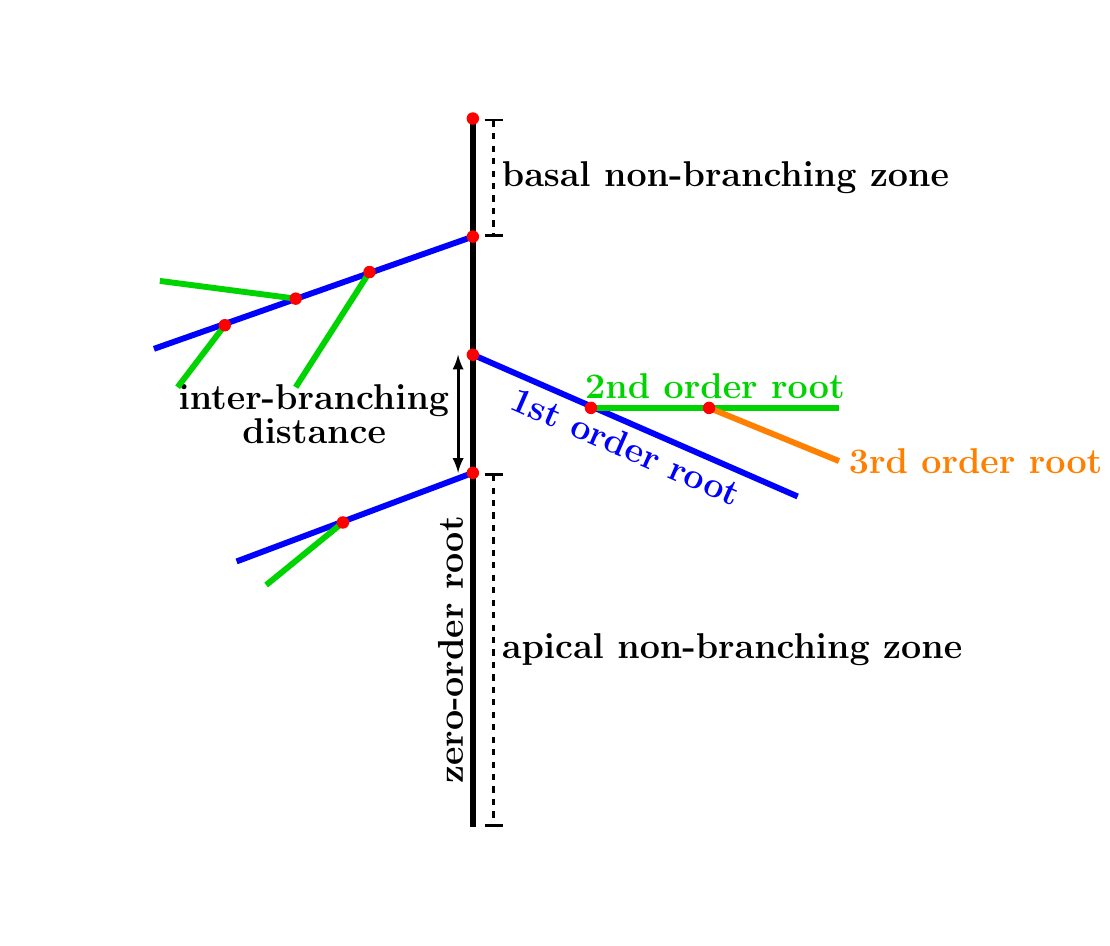}
    \caption{Nomenclature for the root system architecture.}
    \label{fig:RSA}
\end{subfigure}
\begin{subfigure}{0.45\textwidth}
\centering
\includegraphics[width=1\textwidth]{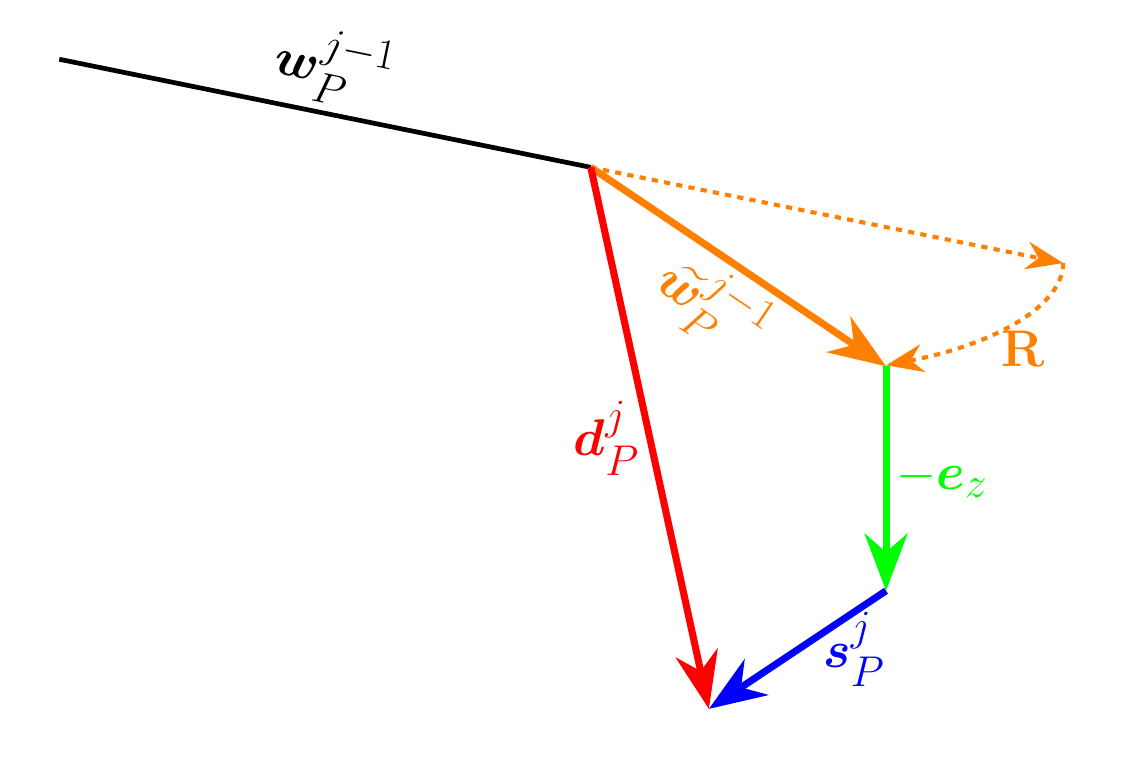}
\caption{Components of $\bm{d}_P^j$}
\label{fig:growth_direction}
\end{subfigure}
\caption{Nomenclature for root growth.}
\label{fig:root_growth}
\end{figure}

A root system may be seen as a set of axes, which are characterized by a ramification order: zero-order roots (or primary axes), which originate from the stem or the seed; first-order lateral roots, which are connected to the primary axes; second-order lateral roots, which are connected to the first-order lateral roots; and so forth \cite{Pages1989} (see Figure \ref{fig:RSA} for a graphical illustration).
The roots grow in an acropetal order, i.e. the youngest roots are the most lateral and closest to the apex. We always assume the seed to be in a known position and to be connected to the soil surface by a vertical mesocotyl \cite{Pages1989}.

To model root growth, we here adopt a \textit{tip-tracking} approach, monitoring the evolution in time of the position of root tips. Let us denote by $\mathcal{P}_\mathrm{tip}^j$ the set of root tips at time $t_j$. We assume that the position $\xx_{\tippoint}$ of a given root tip $\tippoint \in \mathcal{P}_\mathrm{tip}^j$ evolves according to
%(see also \cite{Fitter1991, Pages1989})
\begin{equation}
    \frac{d\xx_{\tippoint}}{dt} =\tipvelocityfun(t):= \bm{w}(\rph{}(t,\xx_{\tippoint}), \theta(t,\xx_{\tippoint}),\xx_{\tippoint} ; \rootparameterset),
    \label{eq:dxdt}
\end{equation}
where $\rootparameterset$ is a set of parameters linked to the genetics of the plant species under study \cite{Jin2020} and to the type of soil in the sample. More precisely, setting $\tipvelocity[](t_j)=\tipvelocity[j]$, we define
\begin{equation}
    \tipvelocity[j] := \growthrate[j]\frac{(1-\repulsionfunction^2(\xx_{\tippoint}))\primarydirection[j]+\repulsionfunction^2(\xx_{\tippoint})\repulsiondirection}{||(1-\repulsionfunction^2(\xx_{\tippoint}))\primarydirection[j]+\repulsionfunction^2(\xx_{\tippoint})\repulsiondirection||}.
    \label{eq:Wp}
\end{equation}
The growth rate $\growthrate[j]$ ($[\length \timem^{-1}]$) is defined as (see \cite{Moraes2018, Seidel2022})
\begin{equation}
\growthrate[j] = {V}_{a} (t_j) \ \mathrm{Imp}(\soilstrength(t_j, \xx_{\tippoint}))\ \mathrm{Imp}(\rph{}(t_j,\xx_{\tippoint})).
    \label{eq:growth_rate}
\end{equation}
The maximum root elongation rate ${V}_{a} \ [\length \timem^{-1}]$ can be a function of time through the root age, for example, whereas $\mathrm{Imp}(\soilstrength)$ and $\mathrm{Imp}(\rph{}) $ are dimensionless impedance factors modeling the impact of soil properties on the actual elongation rate. Let us introduce the soil-strength function $\soilstrength(t,\xx)\ ([\mathrm{MPa}])$, which can be defined empirically as \cite{Clausnitzer1994}:
\begin{equation}
    \soilstrength = \soilstrength_{\max}(1 - \Theta)^3,\label{eq:soil_strength}
\end{equation}	
where $\soilstrength[\max]\ [\mathrm{MPa}]$ is a parameter that incorporates the effects of soil texture and bulk density and 
\begin{equation*}
    \Theta=\frac{\wc-\wc_r}{\wc_s-\wc_r}
\end{equation*}
is the effective saturation, with $\wc_r$ and $\wc_s$ denoting the residual and the saturated volumetric water content, respectively.
The soil-strength impedance factor $\mathrm{Imp}(\soilstrength)$ is then defined as
\begin{equation}
\mathrm{Imp}(\soilstrength) = \begin{cases}
		0 & \text{if } \soilstrength \geq \soilstrength[\max],\\
		1 - \frac{\soilstrength}{\soilstrength_{\max}} &\text{if } \soilstrength < \soilstrength[\max],
		\end{cases}\label{eq:soil_strength_impedance}
\end{equation}
where a reduced mechanical impedance results in faster elongation. Under non-optimal conditions related to a water deficit (drought) or poor aeration (hypoxia), the root elongation is reduced by means of the stress factor, depending on $\rph{}(t,\xx)$ \cite{Moraes2018, Feddes1978}:
\begin{equation}
\mathrm{Imp}(\rph{}) = \begin{cases}
		0 & \text{if } \vert \rph{} \vert \leq \vert (\rph{})_1 \vert, \\
		\frac{\vert (\rph{})_1 \vert - \vert \rph{} \vert}{\vert (\rph{})_1 \vert - \vert (\rph{})_2 \vert} &\text{if } \vert (\rph{})_1 \vert < \vert \rph{} \vert \leq \vert (\rph{})_2 \vert,\\
		1 & \text{if } \vert (\rph{})_2 \vert < \vert \rph{} \vert \leq \vert (\rph{})_3 \vert,\\
		\frac{\vert (\rph{})_4 \vert - \vert \rph{} \vert}{\vert (\rph{})_4 \vert - \vert (\rph{})_3 \vert} &\text{if } \vert (\rph{})_3 \vert < \vert \rph{} \vert \leq \vert (\rph{})_4 \vert,\\
		0 & \text{if } \vert \rph{} \vert >\vert (\rph{})_4 \vert .\\
\end{cases}\label{eq:pressure_head_impedance}
\end{equation}
The values $\vert (\rph{})_1\vert$ and $\vert (\rph{})_4 \vert$ correspond to the hypoxia and drought thresholds, respectively, whereas the pressure head is assumed to be optimal for root growth between $\vert (\rph{})_2 \vert$ and $\vert (\rph{})_3 \vert$. 

Other factors, such as temperature or nutrient concentration, can influence the root elongation rate \cite{Clausnitzer1994, Wu2005}. However, for simplicity, we assume these factors to be constant and equal to their optimal values, thereby posing no impedance to root growth. 

The growth direction in \eqref{eq:Wp} is computed as the convex combination of two contributions: $\primarydirection[j]$, which considers plant genetics and environmental conditions, and $\repulsiondirection[j]$, which accounts for the presence of physical obstacles. In particular, we define
\begin{equation}
    \primarydirection[j] := \frac {k_s \soilstrengthdirection[j]-k_g\ee_z  +k_w\tipvelocitypert[j-1]}{||k_s \soilstrengthdirection[j]-k_g\ee_z  +k_w\tipvelocitypert[j-1]||}.
    \label{eq:standard_growth_direction}
\end{equation}

The vector field $\soilstrengthdirection[j]$ accounts for mechanical constraints and it is defined as
\begin{equation*}
    \soilstrengthdirection[j]=-\frac{\nabla \soilstrength(t_j,\xx_{\tippoint})}{\norm{\nabla\soilstrength(t_j,\xx_{\tippoint})}}.
\end{equation*} 
Since a higher soil water content results in lower soil strength, roots will tend to grow towards the wetter and hence less resistant soil regions. This tendency is called \textit{hydrotropism}. The positive dimensionless parameter $k_s$ in \eqref{eq:standard_growth_direction} allows to tune the weight of hydrotropism with respect to other tropisms: it can depend on the species, root order, and age.
		
The term $-k_g\ee_z$ in \eqref{eq:standard_growth_direction} allows us to account for \textit{geotropism}, which is the tendency of roots to grow downwards. The weight $k_g$ is dimensionless and it may depend again on species, root age, and order.

Finally, the term $k_w\tipvelocitypert[j-1]$ summarizes the plant roots tendency to preserve an already established growth direction (\textit{exotropism}) while exploring the surrounding environment \cite{Clausnitzer1994, Koch2018}. We set
\begin{equation}
    \tipvelocitypert[j-1] := \frac{\mathbf{R}\tipvelocity[j-1]}{\norm{\mathbf{R}\tipvelocity[j-1]}},
\end{equation} 
where matrix $\mathbf{R}$ is a random perturbation of the identity matrix $\mathbf{I}$ such that
\begin{equation*}
    \mathbf{R}=\mathbf{I}+m_a\begin{bmatrix}-m_2^2-m_3^2& m_1m_2 &m_1m_3\\m_1m_2& -m_1^2-m_3^2 &m_2m_3\\m_1m_3& m_2m_3 &-m_1^2-m_2^2\end{bmatrix}.
\end{equation*}
The coefficients $m_i$, $i=1,\dots,3$, normalized such that $\sum_{i=1}^3m_i^2=1$, are randomly chosen in $[0,1]$ according to a uniform distribution each time a growth direction is computed. The parameter $m_a$ weighs how much the identity is actually perturbed and it allows us to introduce some stochasticity in the root tip trajectory, taking into account the space-exploring nature of roots. As for the other types of tropism, the weight $k_w$ is dimensionless and may depend on species, root age and order. The same matrix $\mathbf{R}$ was used in the tip-tracking strategy adopted in \cite{GrappeinGiverso2023} to model the random orientation of extracellular matrix fibers in angiogenesis simulations.

Going back to \eqref{eq:Wp}, the vector field $\repulsiondirection[j]$ depends on the position of pot walls or other physical obstacles, such as stones, and it allows to avoid roots growing through impenetrable boundaries, modeling \textit{thigmotropism}. Generally, roots are negatively thigmotropic, meaning that they have a tendency to grow away from obstacles after contact. Let us denote by $\mathcal{B}$ the boundary of an obstacle and let $d_{\mathcal{B}}(\bm{x})$ be a given distance function measuring the distance from $\mathcal{B}$. Let us then define the \textit{repulsion function}
\begin{equation}
    \repulsionfunction(\xx):=
    \begin{cases}
    \cfrac{d_{\mathrm{max}}-d_{\mathcal{B}}(\bm{x})}{d_{\mathrm{max}}}&\text{if }d_{\mathcal{B}}(\bm{x})\leq d_\mathrm{max},
    \\0 &\text{otherwise},
    \end{cases}
    \label{eq:repulsionfunction}
\end{equation}
and 
\begin{equation*}
    \repulsiondirection := \begin{cases}
    -\cfrac{\nabla \repulsionfunction(\xx_{\tippoint})}{\norm{\nabla \repulsionfunction(\xx_{\tippoint})}} &\text{if }d_{\mathcal{B}}(\bm{x}_{\tippoint}) < d_\mathrm{max},\\
    \bm{0} \in \R^3 &\text{otherwise}.
    \end{cases}
\end{equation*}
According to \eqref{eq:Wp}, the combined contribution of exotropism, geotropism and hydrotropism dominates when $\repulsionfunction(\xx_{\tippoint})\ll1$, while when $\repulsionfunction(\xx_{\tippoint})\sim 1$, i.e. when the root tip is very close to an obstacle, the growth direction is almost completely deviated in order to avoid the obstacle. Furthermore, marking the soil surface as an impenetrable boundary, the contribution of $\repulsiondirection[j]$ allows us also to avoid the roots exiting through the air interface. More details about the repulsion function $\repulsionfunction$ will be given in Remark \ref{rem:repulsionfunction}.

\subsection{Branching}\label{sec:branching}

\begin{figure}
\centering
{\includegraphics[width=0.4\textwidth]{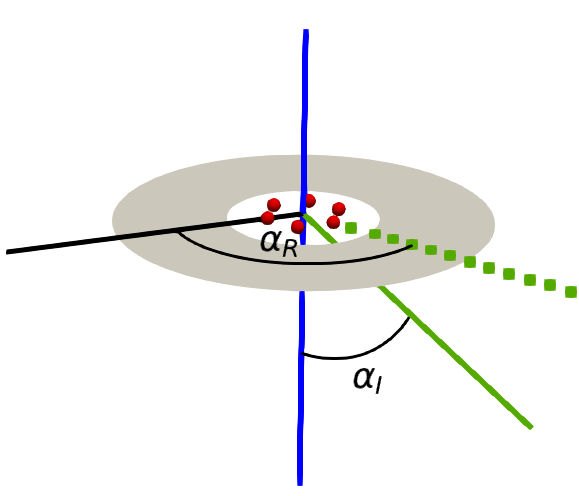}}
\caption{Branching angle and direction.}
\label{fig:branching_angle}
\end{figure}

As the RSA grows, new zero-order roots can emerge and existing roots can branch. The branching process can be split into two main steps: determining the branching position and the branching direction.
Under normal conditions, no branching either occurs at some distance from the axis basis (the length of the basal non-branching zone $L_B$) or from the root tip (length of the apical non-branching zone $L_A$). For primary axes, the axis basis corresponds with the position of the seed; for secondary axis, not originating directly from the seed, the axis basis corresponds instead with the originating branching point on the parent axis. Potential branching nodes can arise between these two zones at regular intervals defined by an inter-branching distance $I$ \cite{Pages1989}. Let us denote by $\ell_{\tippoint}^j$ the length of a given root axis ending in a root tip $P\in \mathcal{P}_\mathrm{tip}^j$.
The maximum number of potential nodes at time $t_j$ for the root axis ending in $\xx_P$ is given by
\begin{equation*}
N_b = \begin{cases}
0 & \text{if } \ell_{\tippoint}^j  < L_A + L_B\\
\left\lfloor \cfrac{\ell_{\tippoint}^j - (L_A + L_B)}{I} \right\rfloor + 1 & \text{if } \ell^j_{\tippoint}  \geq L_A + L_B.\\
\end{cases}
\end{equation*}
At each time step, each of these potential branching sites may either develop into an emerged lateral root or remain unbranched, according to a Bernoulli distribution of probability $p_{br}$ \cite{Pages1989,Pages2019}. As in \cite{Fitter1987}, we model the probability $p_{br}$ of generating a new branch from each potential node as a decreasing function of the order $\omega$ of the root on which the potential node is located. In particular, we set
\begin{equation}\label{eq:branching prob}
p_{br}(\omega)=\begin{cases}
    \frac{e^{-b_c(\omega+1)}}{\sum_{i=0}^{\omega_{\mathrm{max}}}e^{-b_c(i+1)}} & \text{if } \omega = 0,\dots,\omega_{\max}-1,\\
    0 & \text{if } \omega = \omega_{\max}.
\end{cases}
\end{equation}
For example, the root systems of agricultural plants like maize, wheat, etc. typically have $\omega_{\max}=2$ or $\omega_{\max}=3$ as maximum branching order \cite{Roose2004}.
If we choose $b_c \rightarrow 0$ branching becomes equiprobable at all potential nodes of all orders, while branching tends to be restricted to the nodes of the primary axes for bigger values of $b_c$. We remark that there are no potential nodes on roots of order $\omega_{\max}$. Of course, different choices for this branching probability function are possible. Also, different strategies to define potential nodes and to account for the emergence of branches from potential nodes can be adopted (see, for example, \cite{Ziegler2019}).

Denoting by $\mathcal{P}_\mathrm{branch}^j$ the set of nodes at which branching actually occurs at time $t_j$, the direction $\branchingdirection[j]$ of a branch originating from a point $P \in \mathcal{P}_\mathrm{branch}^j$ at time $t_j$, is determined by two angles: an insertion angle $\alpha_I \in [0, \pi)$ and a radial angle $\alpha_R \in [0, 2\pi)$ \cite{Pages1989, Fitter1991}.
The insertion angle $\alpha_I$ is defined as the angle between the parent root and the branch in the plane defined by these two axes. %) and it can be initially set by the user, but it then must decline progressively to a predetermined final value.
The radial angle $\alpha_R$ is defined as the angle between a given direction and the branch direction in the plane perpendicular to the parent root. It can be defined as
\begin{equation*}
\alpha_R := \frac{2 \pi N_R}{X},
\end{equation*}
where $X$ is the number of xylem poles where new lateral roots emerge, which are assumed to be uniformly distributed, and $N_R$ is a random integer number between $1$ and $X$
(see Figure \ref{fig:branching_angle} for a graphic illustration).

\subsection{Managing the intersections between root segments}
The growth of the RSA and the creation of root branches raises the question on how a more complex $\Sigmaj$, characterized by the presence of multiple intersecting root segments, can be tackled.
	
For each time interval $\mathcal{I}_j$, we assume that $\Sigmaj$ can be covered by a set of straight cylindrical root segments $\{\Sigmaj_i\}_{i\in \Yj}$, with $\Yj$ being the set of indexes of root segments at time $t_j$. We denote by $\Gammaj_i$ the lateral surface of $\Sigmaj_i$, and by $\Lambdaj_i = \{ \bm{\lambda}(s),\ s \in (0_i, S_i)\}$ its centerline. We assume that there is no overlap among the centerlines of different root segments and that their closure can intersect at most in one point $\bm {x}_b=\bm{\lambda}(s_b)$, $b \in \Bj$. Furthermore, let 
\begin{equation*}
    \Lambda^j = (\displaystyle\bigcup_{i \in \Yj} \Lambda_i^j) \cup \{\bm{x}_b\}_{b \in \Bj}.
\end{equation*}
We denote by $\Yj_b$ the set of indexes of the segments that are adjacent to an intersection point $\xx_b$ and by $\#\Yj_b$ its cardinality.
We use $\Lambda_i^{j,+} = \bm{\lambda}(s_b^+)$ and $\Lambda_i^{j,-} = \bm{\lambda}(s_b^-)$ to denote the inflow and outflow end-point of $\Lambda_i^j$, respectively. Clearly $s_b^+$ and $\ s_b^-$ can be either $0_i$ or $S_i$ for each $i \in \Yj$. Similarly, we split $\Yj_b$ as $\Yj_b=\Yj[j,+]_b \cup \Yj[j,-]_b$ to distinguish between the segments through which the flux enters in $\bm{x}_b$, indexed by $i \in \Yj[j,+]_b$  and the ones through which it exits, indexed by $i \in \Yj[j,-]_b$.

%We denote by $\Yj_b$ the set of indexes of the segments which are adjacent to an intersection point $\bm {x}_b$, by $\#\Yj_b$ its cardinality, and by $\Bj_i$ the set of the indexes of the intersection points which are the endpoints of a centerline $\Lambda^j_i$. We use $\Lambda_i^{j,+} = \bm{\lambda}(s_b^+)$ and $\Lambda_i^{j,-} = \bm{\lambda}(s_b^-)$ to denote respectively the inflow and outflow end-point of $\Lambda_i^j$. Clearly $s_b^+$ and $\ s_b^-$ can be either $0_i^j$ or $S_i$ for each $i \in \Yj$. Similarly, we split $\Yj_b$ as $\Yj_b=\Yj[j,+]_b \cup \Yj[j,-]_b$ to distinguish between the segments through which the flux enters in $\bm{x}_b$, indexed by $i \in \Yj[j,+]_b$  and the ones through which it exits, indexed by $i \in \Yj[j,-]_b$.
	 
\begin{figure}[t]
    \centering
    \begin{subfigure}{0.45\textwidth}
        \centering
        \includegraphics[width=0.7\textwidth]{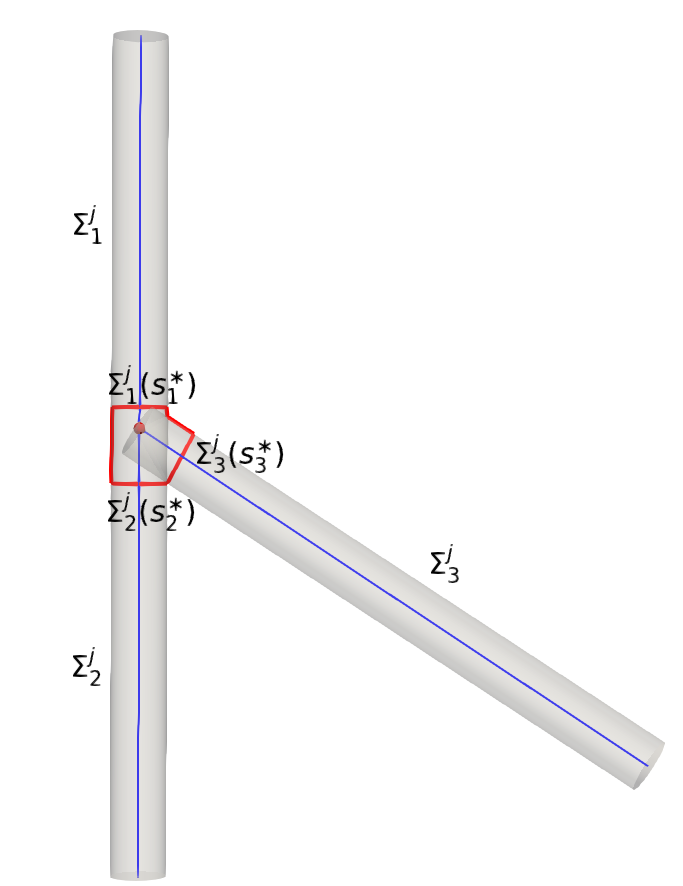}
        \caption{Branch originating from a parent root segment.}
        \label{fig:branching}
    \end{subfigure}
    \begin{subfigure}{0.45\textwidth}
        \centering
        \includegraphics[width=0.7\textwidth]{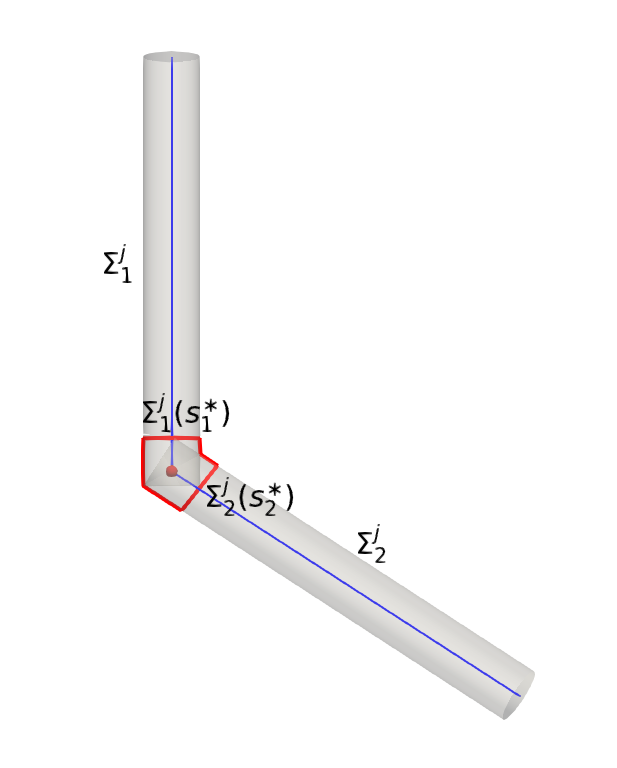}
        \caption{Adjacent root segments with different orientation.}
        \label{fig:linea_spezzata}
    \end{subfigure}
    \caption{Nomenclature at branching points. The boundaries of the extended intersection volume are highlighted in red.}
    \label{fig:branching_scheme}
\end{figure}

The model reduction proposed in Section \ref{sec:optimization} can be performed separately for each root segment, provided that proper conditions are imposed. As in \cite{Grappein2022}, we introduce, for each $\bm{x}_b$, an \textit{extended intersection volume} $\mathcal{V}_b$ (see Figure \ref{fig:branching_scheme}), whose diameter is assumed to be much smaller than the minimum length of the intersecting root segments. In particular, we assume that
\begin{enumerate}[label=\textbf{A.\arabic*}]
    \setcounter{enumi}{4}
    \item in the intersection volume, the xylem pressure $\xp$ has a unique constant value, coinciding with the value in $\xx_b$ and denoted by $\xp^b$.
\end{enumerate}
The two configurations reported in Figure \ref{fig:branching_scheme} refer to the most common cases of root segment intersection which can occur in the proposed growth model: in Figure \ref{fig:branching} a branch originates from a straight root segment, which is in turn split at the branching point into two root segments; in Figure \ref{fig:linea_spezzata} a change in the growth direction produces two adjacent segments with different orientation, which are treated as two intersecting straight segments. In both cases, the faces ${\Sigmaj_i}(s_i^*)$, $i \in \Yj_b$ (see Figure \ref{fig:branching_scheme}), which belong to the boundary of the intersection volume, can be seen as the extreme sections of disjoint sub-cylinders on which the trace operator $\tg$ and extension operators $\es$, $\eg$ are well defined. By an abuse of notation, we will still denote such sub-cylinders as $\Sigmaj_i$.  	
Denoting by $\uux[i]{}$ and by $\xph[i]{}$ the restrictions of functions in $\Vxj$ and in $\Qxj$ to $\Sigmaj_i$ for each $i \in \Yj$, respectively, Equations \eqref{eq:3Dvarproblem_root_mom}-\eqref{eq:3Dvarproblem_root_mass} can be written as
\begin{align*}
    &\sum_{i \in \Yj} \left[\scal[\Sigmaj_i]{\viscr \nabla \bm{u}_{\chi,i}}{ \nabla \bm{v}_{\chi,i}} + \scal[\Sigmaj_i]{\kar \bm{u}_{\chi,i}}{\bm{v}_{\chi,i}} - \scal[\Sigmaj_i]{\psi_{\chi,i}}{ \div \bm{v}_{\chi,i}} + \scal[\Sigmaj_i]{ \ee_z}{\bm{v}_{\chi,i}}\right]\\&\hspace{4.8cm}- \displaystyle\sum_{ b \in \Bj} \sum_{i \in \Yj_b} \scal [\Sigmaj_i(s_i^*)]{\frac{\nnx[i]^T({\xst[i]} \nnx[i])}{\rho g}}{\vvx[i]{} \cdot \nnx[i]}= 0 & \forall \vvx{} \in \Vxj,\\
    &\sum_{i \in \Yj} \scal[\Sigmaj_i]{\div \bm{u}_{\chi,i}}{q_{\chi,i}} = 0 & \forall \qx{} \in \Qxj.
\end{align*} 
These equations must be coupled with proper matching conditions at the interfaces with the extended intersection volume. To impose flux conservation at the intersection points, we impose flux conservation inside the intersection volumes, i.e.
\begin{gather}
    \sum_{i \in \Yj_b} (\fluxwx[i]{}  \cdot \nnx[i])_{|\Sigma_i^j(s_i^*)} = 0\quad \forall b \in \Bj.
    \label{eq:flux_conservation}
\end{gather}
Let us remark that $\bm{n}_{\chi,i}=\bm{e}_s$ on $\Sigmaj_i(s_i^*)$ and that, according to Assumption \ref{ass:velocity}, 
\begin{equation*}
    \bm{u}_{\chi,i}\cdot \bm{e}_s=f(rR^{-1})(\overline{\bm{u}}_{\chi,i})_s =f(rR^{-1})\widehat{u}_{\chi,i}.
\end{equation*}
Condition \eqref{eq:flux_conservation} can hence be rewritten as 
\begin{equation*}
    \sum_{i \in \Yj[j,+]_b} (\ux[i])_{|\Lambda_i^{j,-}}  -\sum_{i \in \Yj[j,-]_b}   (\ux[i])_{|\Lambda_i^{j,+}}  = 0,\quad~ \forall b \in \Bj, 
\end{equation*}
and, since it is essential for the flux, we can impose it directly into the function space, redefining $\Vrxj$ as
\begin{equation} 
    \Vrxj = \Big\{\vx \in \bigcup_{i \in \Yj} \sob{1}{\Lambdaj_i}: \left[\sum_{i \in \Yj[j,+]_b} (\vx[i])_{|\Lambda_i^{j,-}}  -\sum_{i \in \Yj[j,-]_b}   (\vx[i])_{|\Lambda_i^{j,+}}\right]  = 0,\quad~ \forall b \in \Bj \Big\}.
    \label{eq:newVj}
\end{equation}

Recalling that $\bm{\sigma}_\chi=-p_\chi \bm{I}+\overline{\nu}\nabla \bm{u}_\chi$ and exploiting \eqref{eq:nabla_velocity}, the balance of normal forces at the interfaces with the intersection volume can hence be written as
\begin{equation}
    (\nnx[i]^T({\xst[i]} \nnx[i]))_{|\Sigma_i^j(s_i^*)} = -p_\chi^b\quad \forall i \in \Yj_b\quad \forall b \in \Bj.
    \label{eq:normal_forces}
\end{equation}
Let us redefine $\Qrj$ as 
\begin{gather}
    \Qrj = \Big\{\qx{a} \in \con{0}{\Lambdaj}: \qx[i]{a} \in \sob{1}{\Lambdaj_i} \ \forall i \in \Yj \text{ and } \qx{a} = 0 \text{ at } \Lcoll\cup \Ltip\Big\},
    \label{eq:newQj}
\end{gather}
and let $\widehat{\psi}_\chi\in \Qrj$ be such that ${\xph{}}_{|{\Sigma_i^j}}=(\es\widehat{\psi}_{\chi})_{|{\Sigma_i^j}}$ and such that $\xph{}$ is the uniform extension to the whole intersection volume of $\widehat{\psi}_\chi^b = \frac{\xp^b}{\rho g}$. Then according to assumptions \ref{ass:pressure}--\ref{ass:velocity} and exploiting \eqref{eq:normal_forces} and the new definition of $\Vrxj$ we have that
\begin{equation}
    -\displaystyle \sum_{i \in \Yj_b} \scal [\Sigma_i^j(s_b)]{\frac{\nnx[i]^T({\xst[i]} \nnx[i])}{\rho g}}{\vvx[i]{} \cdot \nnx[i]} \!= \pi R^2\widehat{\psi}_{\chi}^b\left[\sum_{i \in \Yj[j,+]_b} (\vx[i])_{|\Lambda_i^{j,-}}  -\sum_{i \in \Yj[j,-]_b}   (\vx[i])_{|\Lambda_i^{j,+}}\right]=0 \quad \forall b \in \Bj.
    \label{eq:balanceofnormalstress}
\end{equation}
Hence, we can conclude that the optimal control problem in the case of multiple cylinders reads as problem \eqref{eq:optimization_problem} where the norms on $\Lambdaj$ are now defined as broken norms drawn according to the partition $\{\Lambdaj_i\}_{i \in \Yj}$ and the 1D functional spaces are defined as in \eqref{eq:newVj}-\eqref{eq:newQj}.

\section{The problem discretization and the solving strategy}\label{sec:discretization}
In this section, we provide the details on the discretization of problem \eqref{eq:optimization_problem}. As previously mentioned, root growth is tackled using a discrete tip-tracking approach: at time $t_{j-1}$, Equation \eqref{eq:dxdt} is solved using a forward Euler scheme, leading to the RSA configuration at time $t_j$; we then consider the evolution of the quantities of interest in the interval $\mathcal{I}_j=(t_{j-1},t_j]$ on the new fixed geometry, using a backward Euler scheme on a sub-partition of $\mathcal{I}_j$. The space discretization of the soil sample is instead tackled by the Virtual Element Method (VEM). Using a polytopal discretization enhances the method with great flexibility in handling geometrically complex samples, characterized by different layers or impervious obstacles, such as stones. Finally, the space discretization of the 1D variables is carried out by a mixed formulation with continuous pressure which imposes flux conservation at junctions in a strong way.

\subsection{Discrete root growth}
Let us assume that the soil pressure head and the water content have been computed at time $t_{j-1}$. Setting $\Delta \Ij=t_j-t_{j-1}$, the position of each root tip $\tippoint \in \mathcal{P}_{\mathrm{tip}}^{j-1}$ is updated as
\begin{equation}
    \xx_{\tippoint}^j=\xx_{\tippoint}^{j-1}+\Delta \Ij \tipvelocity[j],
    \label{eq:update_pos_tip}
\end{equation}
where $\tipvelocity[j]$ is defined as in \eqref{eq:Wp}.

If branching occurs during time interval $\Ij$, i.e. if $\mathcal{P}_\mathrm{branch}^{j-1} \neq \emptyset$, we need to compute also the position of the tips of the root segments originating from the points in $\mathcal{P}_\mathrm{branch}^{j-1}$. In particular, a new root tip $\tippoint^\star$ is created for each $\branchingpoint \in \mathcal{P}_\mathrm{branch}^{j-1}$. The position of the new root tip $\tippoint^\star$ at the time $t_j$ is computed as
\begin{equation}
    \xx_{\tippoint^\star}^j =\xx_{\branchingpoint}^{j-1} + \Delta \Ij V_{\tippoint^{\star}}^j \branchingdirection[j-1], \quad \forall P \in\mathcal{P}_\mathrm{branch}^{j-1},
    \label{eq:update_pos_branch}
\end{equation}
with $V_{\tippoint^{\star}}^j$ computed as in \eqref{eq:growth_rate} and $\branchingdirection[j-1]$ determined from the insertion and radial angles, as detailed in Section \ref{sec:branching}. The final number of root tips in $\mathcal{P}_\mathrm{tip}^j$ is hence given by the sum of the number of elements in $\mathcal{P}_\mathrm{tip}^{j-1}$ and in $\mathcal{P}_\mathrm{branch}^{j-1}$.
When the position of all the root tips in $\mathcal{P}_{\mathrm{tip}}^j$ has been computed, either according to \eqref{eq:update_pos_tip} or to \eqref{eq:update_pos_branch}, the new tips are connected to the originating tip or branching point by a straight line $\Lambda_{\tippoint}$, and the 1D root network is updated as
\begin{equation*}
    \Lambdaj=\Lambda^{j-1}\cup \bigcup_{\tippoint \in \mathcal{P}_\mathrm{tip}^j}\Lambda_{\tippoint},
\end{equation*}
which represents the fixed RSA on which the quantities of interest will evolve for $t \in \mathcal{I}_j$. 

\subsection{Time discretization of the constraint equations}\label{sec:time_discretization}
A backward Euler scheme is adopted for the time discretization of the constraint equations \eqref{eq:3D1Dproblem1}-\eqref{eq:3D1Dproblem3} for $t \in \Ij$. We consider a uniform partition of the interval $\Ij$ with a step $\Delta t \leq \Delta\Ij$ and $t_{j,n}=t_{j-1}+n\Delta t$, $n\geq0$. We use the superscript $n$ to refer to the value of a quantity of interest at time $t_{j,n}$.  The set of the time-discrete constraints in the time interval $\Ij$ hence reads as:
\textit{For each $n = 1,\dots,N$, find} $(\rph{}^n,\ux^n, \xph{a}^n, \lambdas^n, \lambdax^n) \in \Qsj \times \Vrxj \times \Qrj  \times \Qrj \times \Qrj$ \textit{such that}
\begin{equation}
\resizebox{0.91\hsize}{!}{$
\begin{cases}
\scal[\Omega]{\rc({\rph{}^n})\frac{\rph{}^n - \rph{}^{n-1}}{\Delta t}}{\qs{}} + \scal[\Omega]{\rK(\rph{}^n) \nabla \rph{}^n}{\nabla \qs{}} + \scal[\Lambdaj]{2 \pi R \Lp( \rph{a}^n - \lambdax^n)}{\qs{a}} + \scal[\Omega]{\rK(\rph{}^n) \ee_z}{\nabla \qs{}} =0& \forall \qs{} \in \Qsj,\\
\\
\scal[\Lambdaj]{ \kx \ux^n }{\vx} - \scal[\Lambda^j]{ \xph{a}^n }{ \pi R^2 \frac{\partial \vx}{\partial s}}  = - \scal[\Lambda^j]{  \pi R^2  \ee_z \cdot  \ee_s}{\vx} & \forall \vx \in \Vrxj,\\
\scal[\Lambdaj]{\pi R^2  \frac{\partial \ux^n}{\partial s}}{ \qx{a}} + \scal[\Lambdaj]{ 2\pi R \Lp (\xph{a}^n - \lambdas^n)}{\qx{a}} = 0 & \forall  \qx{a} \in \Qrj,
\end{cases}$}
\label{eq:timevardiscrete}
\end{equation}
where, for $j \geq 2$ and $n=1$, we take $\rph{}^{n-1}=\rph{}(t_{j-1})$, while for $j=1$ and $n=1$, $\rph{}^{n-1}=\rph{}(t_{0})=\ri$, according to \eqref{eq:initial_cond}.
%where, for example,
%\begin{align*}
%\rph{a}^n = \rph{a}(t_{n},\cdot)&\quad \forall n =1,\dots,N.
%\end{align*}

\subsection{Space discretization of the constraint equations}
For the discretization of the 3D-1D coupled problem, we extend the 3D soil domain $\Dj$ to the whole $\Omega$, and we consider a decomposition of $\Th$ of $\Omega$ into non overlapping polyhedrons $E$ which satisfy standard mesh assumptions \cite{BEIRAODAVEIGA20171110}. We denote by $\Fh[E]$ and $\Eh[E]$  the set of faces and of edges of a generic polytope $E$. Finally, we fix $\Fh = \bigcup_{E\in\Th} \Fh[E]$. 

Given an integer $k \geq 1$ and a polytope $E$, we denote by $\Poly[d]{k}{E}$ the space of polynomials in $\mathbb{R}^d$ of degree at most $k$ on $E$, with $n^d_k = \dim \Poly[d]{k}{E} = \frac{(k+1)\dots (k+d)}{d!}$. We further adopt the standard convention $ \Poly[d]{-1}{E} = \{0\}$ and $n_{-1}^d = 0$.
Moreover, let us introduce the local polynomial projectors $\proj{\nabla,E}{k} : \sob{1}{E} \to \Poly[d]{k}{E}$, $\proj{0}{k} : \leb{2}{E} \to \Poly[d]{k}{E}$, such that for any $q \in \sob{1}{E}$
\begin{equation*}
 \scal[E]{\nabla q - \nabla \proj{\nabla, E}{k} q }{\nabla p} = 0 \quad \text{and} \quad \int_E \proj{\nabla, E}{k} q = \begin{cases}
\int_{\partial E} q & \text{if } k =1,\\
\int_{E}  q & \text{if } k \geq 2,
\end{cases}\quad \forall p \in \Poly[d]{k}{E},
\end{equation*}
and for any $q \in \leb{2}{E}$
\begin{equation*}
\scal[E]{q - \proj{0, E}{k} q}{p} = 0,\quad \forall p \in \Poly[d]{k}{E}.
\end{equation*} 
For each face $F\in \Fh$,  we first introduce the space
\begin{equation*}
\Bk{}{\partial F} = \left\{q \in \con{0}{\partial F}: q_{|e}\!\in\!\Poly[1]{k}{e} \ \forall e \in \Eh[F] \right\},
\end{equation*}
which allows us to define the two-dimensional Virtual Element space as
\begin{equation*}
    \mathbb{W}(F) = \left\{q \in \sob{1}{F} \cap \con{0}{F}: q_{|\partial F}\!\in\!\Bk{}{\partial F},\ \Delta q \in \Poly[2]{k}{F},\  \scal[F]{q - \proj{\nabla,F}{k} q}{p} = 0 \ \forall p \in \Poly[2]{k}{F} \setminus\Poly[2]{k-2}{F}  \right\},
\end{equation*}
where $\Poly[2]{k}{F} \setminus\Poly[2]{k-2}{F}$ denotes the set of two-dimensional polynomials defined on $F$ of order $k-1$ or $k$.
For each polyhedron $E \in \Th$ we consider the boundary space
\begin{equation*}
    \mathbb{U}(\partial E) = \left\{q \in \con{0}{\partial E}: q_{|F}\!\in\!\mathbb{W}(F) \ \forall F \in \Fh[E]  \right\},
\end{equation*} and the three-dimensional virtual element space
\begin{equation*}
\Qsjh[E]{} = \{q \in \sob{1}{E}: q_{| \partial E} \!\in\! \mathbb{U}(\partial E),\ \Delta q \in \Poly[3]{k}{E},\  \scal[E]{q - \proj{\nabla,E}{k} q}{p} = 0 \ \forall p \in \Poly[3]{k}{E} \setminus\Poly[3]{k-2}{E}\}.
\end{equation*}
%where $\Poly[2]{k}{E} \setminus\Poly[2]{k-2}{E}$ is the set of polynomials of degree exactly equals to $k-1$ or $k$.
Gluing together the local spaces, we define the global Virtual Element space for the soil pressure variable as \cite{BEIRAODAVEIGA20171110}
\begin{equation*}
\Qsjh{} = \{q \in \sob[0]{1}{\Omega} \cap \con{0}{\Omega}: q \in \Qsjh[E]{} \ \forall E \in \Th~ \}.
\end{equation*}

Denoting by $\mathcal{P}_{k}^2(F)$ and $\mathcal{P}_{k}^3(E)$, suitable polynomial bases for $\Poly[2]{k}{F}$ and $\Poly[3]{k}{E}$, respectively, for each $\qs{h} \in \Qsjh{}$, we choose the following degrees of freedom (for more details see \cite{LBe14, Teora2024}): 
\begin{itemize}
\item the values of $\qs{h}$ at the vertices of $\Th$;
\item the values of $\qs{h}$ at the $k-1$ internal Gauss-Lobatto quadrature points on each internal edge of $\Th$;
\item the internal moments on each internal face $F\in\Fh$ defined as
\begin{equation}
{\frac{1}{\vert F \vert}}\int_F \qs{h} \ p \qquad \forall p \in \mathcal{P}_{k-2}^2(F);
\label{eq:soil_DOFF}
\end{equation}
\item the internal moments on each element $E \in \Th$ defined as:
\begin{equation}
 {\frac{1}{\vert E \vert}}\int_E \qs{h} \ p \qquad \forall p \in \mathcal{P}_{k-2}^3(E) \quad \forall E \in \Th.
\label{eq:soil_DOFE}
\end{equation}
\end{itemize}
In the following, we denote by $\lbrace\Bps[\alpha]\rbrace_{\alpha=1}^{\Nps}$ the set of Lagrangian basis functions associated with these degrees of freedom and spanning $\Qsjh{}$.

Concerning the 1D variables, we can build on each segment $\Lambdaj_i$, for all $i\in \Yj$, three possibly independent meshes $\xTh[i],\ \xxTh[i]$ and $\ssTh[i]$ for the discretization of $(\ux, \xph{a})$, $\lambdax$ and $\lambdas$, respectively. 
Discrete global spaces for the xylem velocity and pressure head are defined as
\begin{equation}
\Vxjh{} = \{\hat{v} \in \Vrxj:  \hat{v}_{|e} \in \Poly{k + 1}{e}\ \forall e \in \xTh[i], \forall i \in \Yj \},
\label{eq:xylem_velocity_space}
\end{equation}
\begin{equation}
\Qxjh{} = \{\hat{q} \in \Qrj: \hat{q}_{|e} \in \Poly{k}{e}\ \forall e \in \xTh[i],\ \forall i \in \Yj\}.
\label{eq:xylem_pressure_space}
\end{equation}

For the elements $\vvx{h} \in \Vxjh{}$, we choose the following degrees of freedom:
\begin{itemize}
\item the value of $\vvx[i]{h}$ at each internal vertex of $\xTh[i]$, $\forall i \in \Yj$, at the collar and at the root tips;
\item $\#\Yj_b-1$ values of $\vvx[i]{h}$, $i\in \Yj_b$, at each intersection point $\xx_b=\llambda(s_b)$, $b \in \Bj$. Let us indeed remark that the value that is left out can be retrieved from the others, according to \eqref{eq:flux_conservation};
\item if $k > 1$, for each $e \in \xTh[i]$, $\forall i \in \Yj$, we consider the scaled moments
\begin{equation}
\frac{1}{\vert e \vert} \int_e \vvx[i]{h}\ \widehat{g}^{\nabla}, \quad \forall \widehat{g}^{\nabla} \in \nabla  \mathcal{P}^1_{k}(e),
\label{eq:xylem_DOFNabla}
\end{equation}
\end{itemize}
where $\mathcal{P}_{k}^1(e)$ is a suitable polynomial basis for $\Poly[1]{k}{e}$, and $\nabla  \mathcal{P}^1_{k}(e)$ denotes the set containing the gradients of the basis functions in $\Poly[1]{k}{e}$. We denote by $\lbrace \Bux[\alpha]\rbrace_{\alpha=1}^{\Nux}$ the set of Lagrangian basis functions associated to the chosen degrees of freedom and spanning $\Vxjh{}$. 
%at each junction $s_b$, $\forall b \in \Bj$ for each $i \in \Yj$ with the exception of an element $i_0 \in \Yj$. Indeed, the value of $\vx[i_0]{}(s_b)$ can be retrieved as $\vx[i_0]{}(s_b) = \sum_{i \in \Yj \setminus \{i_0\}} \vx[i]{}(s_b)$;}

Let us now turn our attention to the xylem pressure-head space. For each element $\qx{h}\in\Qxjh{}$ we choose as degrees of freedom:
\begin{itemize}
\item the value of $\qx[i]{h}$ at each internal vertex of $ \xTh[i]$, for all $i \in \Yj$, at the collar and at the root tips;
\item one single value of $\qx[i]{h}$ at each intersection point $\xx_b$, since we are assuming the continuity of the pressure head;
\item if $k > 1$, the values at the $k-1$ internal Gauss-Lobatto quadrature points on each element of $\xTh[i]$, for all $i\in\Yj$.
\end{itemize}
The set of Lagrangian basis functions associated with the chosen degrees of freedom and spanning $\Qxjh{}$ is denoted by $\lbrace\Bpx[\alpha]\rbrace_{\alpha=1}^{\Npx}$.

For example, we assign $2$ degrees of freedom for the velocity variable in the configuration shown in Figure \ref{fig:branching}, while we define only one degree of freedom for the velocity variable in the configuration shown in Figure \ref{fig:linea_spezzata}. In both cases, we assign just one degree of freedom for the pressure head variable. 

For the control variables $\lambdas$ and $\lambdax$ we define the spaces $\Qspjh{}$ and $\Qxpjh{}$ in a similar fashion to what was done for the reduced xylem pressure variable, and we denote by $\lbrace\Bcs[\alpha]\rbrace_{\alpha=1}^{\Ncs}$ and $\lbrace\Bcx[\alpha]\rbrace_{\alpha=1}^{\Ncx}$ the respective set of associated Lagrangian basis functions.

\begin{remark}
For ease of notation, we assume the same order of accuracy for the 3D and the 1D variables, as well as for the control variables. However, different orders of accuracy could also be employed.
\end{remark}

From now on we will use uppercase symbols to refer to the space-discrete version of the quantities of interest. For instance, $$\rph{h}=\sum_{\alpha=1}^{\Nps} \rph[\alpha]{h}\Bps[\alpha], \quad \xph{h}=\sum_{\alpha=1}^{\Npx} \xph[\alpha]{h}\Bpx[\alpha], \quad \uux{h}=\sum_{\alpha=1}^{\Nux} \uux[\alpha]{h}\Bux[\alpha], \quad \phis{h}=\sum_{\alpha=1}^{\Ncs} \phis[\alpha]{h}\Bps[\alpha], \quad \phix{h}=\sum_{\alpha=1}^{\Ncx} \phix[\alpha]{h}\Bcx[\alpha] $$ will denote the discrete version of $\rph{},\ \xph{}, \ \ux, \ \phis{}$, and $\phix{}$, respectively. With an abuse of notation, we will use the same uppercase symbol to refer to the vector containing the values assigned to the degrees of freedom of the corresponding discrete variable, the meaning being clear from the context.
\begin{remark}
    The number of degrees of freedom for all the 1D variables actually depends on $j=1,...,J$. However, we decide to drop the index $j$, since we are focusing on the time and space variation of the quantities as $t\in \Ij$, assuming $j$ to be fixed.
\end{remark}

On each element $E \in \Th$ and for each $Z_{\soil} \in \Qsjh{}$, as in \cite{Cangiani2019, Vacca2015}, we define the following local discrete bilinear forms
\begin{align*}
\dbilinh[E]{\rph{h}}{\qs{h}}{Z_{\soil}} &= \scal[E]{\rK(\proj{0,E}{k} Z_{\soil})\proj{0,E}{k-1}\nabla \rph[]{h}}{\proj{0,E}{k-1}\nabla \qs[]{h}} + \stab[\tilde{a},E]{\rph[]{h} - \proj{\nabla,E}{k} \rph[]{h}}{\qs[]{h} - \proj{\nabla,E}{k} \qs[]{h}}{ Z_{\soil}},
\end{align*}
\begin{align*}
\reacbilinh[E]{\rph[]{h}}{\qs[]{h}}{Z_{\soil}} &= \scal[E]{\rc( \proj{0,E}{k} Z_{\soil}) \proj{0,E}{k} \rph[]{h}}{\proj{0,E}{k} \qs[]{h}} + \stab[\tilde{c}, E]{\rph[]{h} - \proj{0,E}{k} \rph[]{h}}{\qs[]{h} - \proj{0,E}{k} \qs[]{h}}{Z_{\soil}},
\end{align*}
where $\stab[\tilde{a},E]{}{}{Z_{\soil}}$ and $\stab[\tilde{c},E]{}{}{Z_{\soil}}$ are any computable, symmetric and positive definite bilinear forms such that $\dbilinh[E]{}{}{Z_{\soil}} $ and $\reacbilinh[E]{}{}{Z_{\soil}}$ scale like their continuous counterparts \cite{LBe13}
\begin{equation*}
\dbilin[E]{\rph[]{h}}{\qs[]{h}}{ Z_{\soil}} = \scal[E]{\rK(Z_{\soil}) \nabla \rph[]{h}}{ \nabla \qs[]{h}},\quad \reacbilin[E]{\rph[]{h}}{\qs[]{h}}{ Z_{\soil}} = \scal[E]{\rc(Z_{\soil})  \rph[]{h}}{  \qs[]{h}},
\end{equation*}
respectively. We then define the matrices
\begin{align*}
&\AA(\rph[]{h}^n) \in \R^{\Nps \times \Nps}: \quad [\AA(\rph[]{h}^n)]_{\alpha \beta} = \sum_{E \in \Th} \dbilinh[E]{\Bps[\beta]}{\Bps[\alpha]}{\rph[]{h}^n},\\
&\CC(\rph[]{h}^n) \in \R^{\Nps \times \Nps}: \quad  [\CC(\rph[]{h}^n)]_{\alpha\beta} = \sum_{E \in \Th} \reacbilinh[E]{\Bps[\beta]}{\Bps[\alpha]}{\rph[]{h}^n},
\end{align*}
where we again use the superscript $n$ to refer to the value of the related discrete variable at time $t_{j,n}$.
To work out the complete discrete matrix formulation of the constraint equations let us also define the matrices
\begin{align}
&\wAA \in \R^{\Nux \times \Nux}: & &\hspace{-2cm}[\wAA]_{\alpha \beta } = \scal[\Lambdaj]{ \kx \Bux[\beta] }{\Bux[\alpha] } ,\label{eq:matrix1D_start}\\
&\wBB \in \R^{\Npx \times \Nux}: & &\hspace{-2cm} [\wBB]_{\alpha \beta} = -\pi R^2\scal[\Lambdaj]{ \Bpx[\alpha] }{ \frac{\partial \Bux[\beta]}{\partial s}} ,\\
&\wMM[\Lp] \in \R^{\Npx \times \Npx}: & &\hspace{-2cm}[\wMM[\Lp]]_{\alpha \beta} = 2\pi R\scal[\Lambdaj]{\Lp \Bpx[\beta]}{\Bpx[\alpha]}, \\
&\MM[\Lp] \in \R^{\Nps \times \Nps}: & &\hspace{-2cm}[\MM[\Lp]]_{\alpha \beta} = 2\pi R\scal[\Lambdaj]{\ \Lp\ \Bps[\beta]}{\Bps[\alpha]} \\
&\sxDD[\Lp] \in \R^{\Nps \times \Ncx  }: & &\hspace{-2cm}[\sxDD[\Lp]]_{\alpha \beta } = 2\pi R\scal[\Lambdaj]{\Lp \Bcx[\beta]}{\Bps[\alpha]}, \\
&\xsDD[\Lp] \in \R^{\Npx \times \Ncs}: & &\hspace{-2cm}[\xsDD[\Lp]]_{\alpha \beta} = 2\pi R\scal[\Lambdaj]{\Lp\ \Bcs[\beta]}{\Bpx[\alpha]} \label{eq:matrix1D_end},
\end{align}
and the vectors
\begin{align*}
&\hspace{1.5cm}\bm{f}_\soil(\rph[]{h}^n) \in \R^{\Nps }: & &\hspace{-0cm} [\bm{f}_\soil(\rph[]{h}^n)]_{\alpha} = -\sum_{E \in \Th} \scal[E]{\rK(\proj{0,E}{k}\rph[]{h}^n) \ee_z}{\nabla \Bps[\alpha]},\\
&\hspace{1.5cm}\widehat{\bm{f}}_\xylem \in \R^{\Nux}: & & \hspace{-0cm}[\widehat{\bm{f}}_\xylem]_{\alpha} = - \pi R^2\scal[\Lambdaj]{ \ee_z \cdot  \ee_s}{\Bux[\alpha]} .
\end{align*}
The space-time discrete matrix formulation of the constraint equations \eqref{eq:3D1Dproblem1}-\eqref{eq:3D1Dproblem3} can finally be written as
\begin{equation}
    \begin{bmatrix} 
    \AA(\rph[]{h}^n) + \frac{1}{\Delta t} \CC(\rph[]{h}^n) + \MM[\Lp] 
    \end{bmatrix}  \rph[]{h}^n
    -\sxDD[\Lp] \phix[]{h}^n -\frac{1}{\Delta t}\CC(\rph[]{h}^n)\rph[]{h}^{n-1} -\bm{f}_\soil(\rph[]{h}^n) = \bm{0},
    \label{eq:discrete_system_3D}
\end{equation}
\begin{equation}
    \begin{bmatrix} 
    \addlinespace \wAA & \wBB^{T} \\
    \addlinespace \wBB & -\wMM[\Lp] 
    \end{bmatrix} \begin{bmatrix} 
    \addlinespace \uux[]{h}^n \\ 
    \addlinespace \xph[]{h}^n 
    \end{bmatrix} + 
    \begin{bmatrix}
    \addlinespace\mathbf{0} \\
    \addlinespace \xsDD[\Lp] 
    \end{bmatrix} \phis[]{h}^n 
     = \begin{bmatrix}
    \addlinespace \widehat{\bm{f}}_\xylem \\
    \addlinespace\mathbf{0} 
    \end{bmatrix}.
    \label{eq:discrete_system_1D}
\end{equation}

\begin{remark}\label{rem:matrix_concat}
    The size of the matrices \eqref{eq:matrix1D_start}-\eqref{eq:matrix1D_end} actually changes as the root network evolves. However, since we are allowing for the growth of the RSA but not for its remodeling or regression, we do not need to recompute the matrices at each time step. Indeed, by numbering correctly the degrees of freedom, we can update the matrices by concatenating the contributions of the new basis functions to the matrices computed at the previous time step.
\end{remark}

\subsection{The discrete cost functional}\label{sec:discr_cost_fun}
The discrete cost functional, to be minimized subject to the time-space discrete constraints \eqref{eq:discrete_system_3D}-\eqref{eq:discrete_system_1D}, is derived from \eqref{eq:costfunctional} by replacing the $\Hc$-norms with norms in $L^2(\Lambdaj)$. Namely, we introduce
\begin{equation}
    \widetilde{\mathcal{J}}( \phis[]{h}^{n}, \phix[]{h}^{n}) = \frac{1}{2} \left( \norm[\leb{2}{\Lambdaj}]{\rph[]{h}^{n}(\phix[]{h}^{n}) - \phis[]{h}^{n}}^2 + \norm[\leb{2}{\Lambdaj}]{\xph[]{h}^{n}(\phis[]{h}^{n}) - \phix[]{h}^{n}}^2\right),
\end{equation}
and a new set of matrices
\begin{align}
&\wMM \in \R^{\Npx \times \Npx  }:& & \hspace{-2cm}[\wMM]_{\alpha \beta} = \scal[\Lambdaj]{\Bpx[\beta]}{\Bpx[\alpha]}, \label{matrix1D_J_start}\\
&\MM \in \R^{\Nps \times \Nps  }:& &\hspace{-2cm}[\MM]_{\alpha \beta} = \scal[\Lambdaj]{\Bps[\beta]}{\Bps[\alpha]},\\
&\wGG \in \R^{\Ncs \times \Ncs  }:& &\hspace{-2cm}[\wGG]_{\alpha \beta} = \scal[\Lambdaj]{ \Bcx[\beta]}{\Bcx[\alpha]}, \\
&\GG \in \R^{\Ncx \times \Ncx  }:& &\hspace{-2cm}[\GG]_{\alpha \beta} = \scal[\Lambdaj]{\Bcs[\beta]}{ \Bcs[\alpha]}, \\
&\ssDD \in \R^{\Nps \times \Ncs  }:& &\hspace{-2cm} [\ssDD]_{\alpha \beta} = \scal[\Lambdaj]{\Bps[\alpha]}{\Bcs[\beta] }, \\
&\xxDD \in \R^{\Npx \times \Ncx  }:& &\hspace{-2cm} [\xxDD]_{\alpha \beta } = \scal[\Lambdaj]{\Bpx[\alpha]}{\Bcx[\beta]} \label{matrix1D_J_end},
\end{align}
such that $\widetilde{\mathcal{J}}$ can be written in matrix form as
\begin{align}
\widetilde{\mathcal{J}}( \phis[]{h}^{n}, \phix[]{h}^{n})=\frac{1}{2} \left( \begin{bmatrix} 
\addlinespace\rph[]{h}^{n}(\phix[]{h}^{n}) \\ 
\addlinespace\xph[]{h}^{n}(\phis[]{h}^{n}) \\
 \addlinespace\phis[]{h}^{n} \\  
 \addlinespace\phix[]{h}^{n,}
\end{bmatrix}^T\begin{bmatrix} 
\addlinespace \MM        & \mathbf{0}  & -\ssDD     & \mathbf{0} \\
\addlinespace \mathbf{0} & \wMM        & \mathbf{0} & -\xxDD \\
\addlinespace -(\ssDD)^T & \mathbf{0}  & \GG        & \mathbf{0} \\
\addlinespace\mathbf{0}  & -(\xxDD)^T  & \mathbf{0} & \wGG 
\end{bmatrix} 
\begin{bmatrix} 
\addlinespace\rph[]{h}^{n}(\phix[]{h}^{n}) \\ 
\addlinespace\xph[]{h}^{n}(\phis[]{h}^{n}) \\
 \addlinespace\phis[]{h}^{n} \\  
 \addlinespace\phix[]{h}^{n}
\end{bmatrix} 
\right).
\label{eq:discrete_cost_functional}
\end{align}
Considerations similar to those reported in Remark \ref{rem:matrix_concat} hold also for the matrices \eqref{matrix1D_J_start}-\eqref{matrix1D_J_end}.

\subsection{The solving strategy}\label{sec:solving_strategy}
Given a fixed $\Ij$, in each sub-interval $(t_{j,n-1},t_{j,n}]$, $n=1,...,N$, the non-linearities characterizing Equation \eqref{eq:discrete_system_3D} are tackled using the Picard non-linear iterative method, which is a globally convergent method that has been extensively analyzed in the case of Richards equation \cite{List2016}. At each non-linear iteration, indexed by $\ell \geq 1$, we introduce the matrix
\begin{gather*}
\AAstar = \begin{bmatrix} 
\AA(\rph[]{h}^{n,\ell-1}) +\frac{1}{\Delta t}\CC(\rph[]{h}^{n,\ell-1}) + \MM[\Lp] 
\end{bmatrix},
\end{gather*}
and the vector
\begin{equation*}
    \fstar = \frac{1}{\Delta t}\CC(\rph[]{h}^{n,\ell-1})\rph[]{h}^{n-1} + \bm{f}_\soil(\rph[]{h}^{n,\ell-1}),
\end{equation*}
where $\rph[]{h}^{n,\ell-1}$ is the soil pressure head obtained at the non-linear step $\ell-1$, 
%For the sake of compactness we also define
%\begin{gather*}
%\wAAstar  = \begin{bmatrix} 
%\addlinespace\wAA & \wBB^{T} \\
%\addlinespace\wBB & -\wMM[\Lp] \\
%\end{bmatrix}  ,\quad \SSstar = \begin{bmatrix}
%\addlinespace \mathbf{0}\\
%\xsDD[\Lp]\end{bmatrix}, \quad \wfstar = \begin{bmatrix}
%\addlinespace \widehat{\bm{f}}_\xylem \\
%\addlinespace \mathbf{0} \\
%\end{bmatrix}
%\end{gather*}
and we rewrite the discrete constraint equations \eqref{eq:discrete_system_3D}-\eqref{eq:discrete_system_1D} as
\begin{gather}
\AAstar \rph[]{h}^{n,\ell} 
= \sxDD[\Lp] \phix[]{h}^{n,\ell} +\fstar,\label{eq:picard_constraints_1}\\
\begin{bmatrix} 
\addlinespace \wAA & \wBB^{T} \\
\addlinespace \wBB & -\wMM[\Lp] 
\end{bmatrix} \begin{bmatrix} 
\addlinespace\uux[]{h}^{n,\ell} \\ 
\addlinespace\xph[]{h}^{n,\ell}
\end{bmatrix} =  
\begin{bmatrix}
\addlinespace\mathbf{0} \\
\addlinespace -\xsDD[\Lp] 
\end{bmatrix} \phis[]{h}^{n,\ell} 
 + \begin{bmatrix}
\addlinespace \widehat{\bm{f}}_\xylem \\
\addlinespace\mathbf{0}
\end{bmatrix}.
\label{eq:picard_constraints_2}
\end{gather}
We then solve the following constrained optimization problem:
\begin{equation}
\begin{gathered}
\textit{Find } (\rph[]{h}^{n,\ell}, \uux[]{h}^{n,\ell}, \xph[]{h}^{n,\ell}, \phis[]{h}^{n,\ell}, \phix[]{h}^{n,\ell}) \in \Qsjh{} \times \Vxjh{} \times  \Qxjh{}  \times \Qspjh{} \times \Qxpjh{} \\
 \textit{such that } \widetilde{\mathcal{J}}( \phis[]{h}^{n,\ell}, \phix[]{h}^{n,\ell}) \textit{ is minimized subject to \eqref{eq:picard_constraints_1}-\eqref{eq:picard_constraints_2}.}
\end{gathered}
\label{eq:discrete_opt_problem}
\end{equation}
The constrained optimization problem \eqref{eq:discrete_opt_problem} can actually be recast into an unconstrained optimization problem, by dropping the constraints into the discrete cost functional. For the sake of compactness, we define
\begin{gather*}
\wAAstar  = \begin{bmatrix} 
\addlinespace\wAA & \wBB^{T} \\
\addlinespace\wBB & -\wMM[\Lp] \\
\end{bmatrix}  ,\quad \wfstar = \begin{bmatrix}
\addlinespace \widehat{\bm{f}}_\xylem \\
\addlinespace \mathbf{0} \\
\end{bmatrix},
\end{gather*} so that
according to \eqref{eq:picard_constraints_1}-\eqref{eq:picard_constraints_2} we have \begin{equation}
\rph[]{h}^{n,\ell}(\phix[]{h}^{n,\ell}) = \left(\AAstar\right)^{-1}  \left( 
 \sxDD[\Lp] \phix[]{h}^{n,\ell} + \fstar\right),\quad \begin{bmatrix} 
\addlinespace\uux[]{h}^{n,\ell} (\phis[]{h}^{n,\ell})\\ 
\addlinespace\xph[]{h}^{n,\ell} (\phis[]{h}^{n,\ell})  
\end{bmatrix} = \wAAstar^{-1} \left(\begin{bmatrix} 
 \mathbf{0} \\
 -\xsDD[\Lp]
\end{bmatrix}  \phis[]{h}^{n,\ell} + \wfstar\right).\label{eq:inv_constraints}
\end{equation}
We then rewrite the discrete cost functional \eqref{eq:discrete_cost_functional} as
\begin{align}
\nonumber \widetilde{\mathcal{J}}(\phis[]{h}^{n,\ell},\phix[h]{h}^{n,\ell}) 
&= \frac{1}{2} \Big[ (\rph[]{h}^{n,\ell})^T \MM \rph[]{h}^{n,\ell} - 2(\rph[]{h}^{n,\ell})^T \ssDD \phis[]{h}^{n,\ell} + (\phis[]{h}^{n,\ell})^T\ \GG\ \phis[]{h}^{n,\ell}\\
\nonumber &\qquad + \begin{bmatrix} 
\addlinespace\uux[]{h}^{n,\ell}\\ 
\addlinespace\xph[]{h}^{n,\ell} 
\end{bmatrix} ^T  \begin{bmatrix}
\mathbf{0} & \mathbf{0} \\
\mathbf{0} & \wMM
\end{bmatrix}  \begin{bmatrix} 
\addlinespace\uux[]{h}^{n,\ell}\\ 
\addlinespace\xph[]{h}^{n,\ell} 
\end{bmatrix} - 2\begin{bmatrix} 
\addlinespace\uux[]{h}^{n,\ell}\\ 
\addlinespace\xph[]{h}^{n,\ell} 
\end{bmatrix}^T  \begin{bmatrix} 
\addlinespace \bm{0}\\ \addlinespace \xxDD \end{bmatrix} \phix[]{h}^{n,\ell} +(\phix[h]{h}^{n,\ell})^T\ \wGG\ \phix[]{h}^{n,\ell} ,\Big]\end{align}
so that, setting $\Jx = \begin{bmatrix}
\addlinespace \phis[]{h}^{n,\ell}\\
\addlinespace \phix[]{h}^{n,\ell}
\end{bmatrix}$ and using \eqref{eq:inv_constraints}, we obtain 
\begin{align}
\widetilde{\mathcal{J}}(\phis[]{h}^{n,\ell},\phix[h]{h}^{n,\ell}) 
%&= \frac{1}{2}\begin{bmatrix}
%\addlinespace \phis[]{h}^{n,\ell}\\
%\addlinespace \phix[]{h}^{n,\ell}
%\end{bmatrix}^T 
%\begin{bmatrix}
%\GG & \mathbf{0}\\
%\mathbf{0} & \wGG
%\end{bmatrix} \begin{bmatrix}
%\addlinespace \phis[]{h}^{n,\ell}\\
%\addlinespace \phix[]{h}^{n,\ell}
%\end{bmatrix}
%\nonumber&- \begin{bmatrix}
%\addlinespace \phis[]{h}^{n,\ell}\\
%\addlinespace \phix[]{h}^{n,\ell}
%\end{bmatrix}^T 
%\begin{bmatrix} (\ssDD)^T (\AAstar)^{-1} \left( 
% \sxDD[\Lp] \phix[]{h}^{n,\ell} + \fstar\right)\\
%\begin{bmatrix}
%\addlinespace \mathbf{0}\\
%\addlinespace \xxDD
%\end{bmatrix}^T
%(\wAAstar)^{-1}
%\left(\begin{bmatrix} 
% \mathbf{0} \\
% -\xsDD[\Lp]
%\end{bmatrix}  \phis[]{h}^{n,\ell} + \wfstar\right)
%\end{bmatrix} \\
%\nonumber&+ \frac{1}{2}\left( 
% \sxDD[\Lp]\ \phix[]{h}^{n,\ell} + \fstar\right)^T (\AAstar)^{-1} \MM  (\AAstar)^{-1} \left( 
% \sxDD[\Lp]\ \phix[]{h}^{n,\ell} + \fstar\right)\\
%\nonumber&+\frac{1}{2}\left(\begin{bmatrix} 
% \mathbf{0} \\
% -\xsDD[\Lp]
%\end{bmatrix} \phis[]{h}^{n,\ell} + \wfstar\right)^T (\wAAstar)^{-1}
%\begin{bmatrix}
%\mathbf{0} & \mathbf{0}\\
%\mathbf{0} & \wMM
%\end{bmatrix} (\wAAstar)^{-1}
%\left(\begin{bmatrix} 
% \mathbf{0} \\
% -\xsDD[\Lp]
%\end{bmatrix}  \phis[]{h}^{n,\ell} + %\wfstar\right)\\
&= \frac{1}{2} \left(\Jx^T \JM\Jx + 2\Jx^T \Jd + \Jb\right),
\label{eq:unconstrained_problem}
\end{align}
where
\begin{equation}
\resizebox{0.91\hsize}{!}{$
\JM = 
\begin{bmatrix}
\GG + \begin{bmatrix} 
 \mathbf{0} \\
 -\xsDD[\Lp]
\end{bmatrix}^T(\wAAstar)^{-1}\begin{bmatrix}\mathbf{0} & \mathbf{0}\\ \mathbf{0} & \wMM \end{bmatrix} (\wAAstar)^{-1}\begin{bmatrix} 
 \mathbf{0} \\
 -\xsDD[\Lp]
\end{bmatrix}& 
-(\ssDD)^T (\AAstar)^{-1}  \sxDD[\Lp] 
- \begin{bmatrix} 
 \mathbf{0} \\
 -\xsDD[\Lp]
\end{bmatrix}^T 
(\wAAstar)^{-1} \begin{bmatrix}
\addlinespace \mathbf{0}\\
\addlinespace \xxDD
\end{bmatrix}\\
-(\sxDD[\Lp] )^T(\AAstar)^{-1}\ssDD
  - \begin{bmatrix}
\addlinespace \mathbf{0}\\
\addlinespace \xxDD
\end{bmatrix}^T
(\wAAstar)^{-1}
\begin{bmatrix} 
 \mathbf{0} \\
 -\xsDD[\Lp]
\end{bmatrix}  & 
\wGG + (\sxDD[\Lp] )^T (\AAstar)^{-1} \MM (\AAstar)^{-1}  \sxDD[\Lp] 
\end{bmatrix},$}
\label{eq:matrix_cg}
\end{equation}
\begin{equation*}
\Jd = \begin{bmatrix}
 -(\ssDD)^T (\AAstar)^{-1}  \fstar
 + \begin{bmatrix} 
 \mathbf{0} \\
 -\xsDD[\Lp]
\end{bmatrix}^T (\wAAstar)^{-1}\begin{bmatrix}\mathbf{0} & \mathbf{0}\\ \mathbf{0} & \wMM \end{bmatrix} (\wAAstar)^{-1}\wfstar\\
-\begin{bmatrix}
\addlinespace \mathbf{0}\\
\addlinespace \xxDD
\end{bmatrix}^T
(\wAAstar)^{-1}
\wfstar + (\sxDD[\Lp] )^T  (\AAstar)^{-1} \MM (\AAstar)^{-1} \fstar 
\end{bmatrix},
\end{equation*}
\begin{equation*}
\Jb = (\fstar)^T  (\AAstar)^{-1} \MM (\AAstar)^{-1} \fstar+ (\wfstar)^T (\wAAstar)^{-1}\begin{bmatrix}\mathbf{0} & \mathbf{0}\\ \mathbf{0} & \wMM \end{bmatrix} (\wAAstar)^{-1}\wfstar.
\end{equation*}

We observe that the matrix $\JM$ is symmetric and positive semi-definite by construction. In \cite{Grappein2023} it has been formally shown that the matrix $\JM$ is symmetric positive definite for the linear and stationary case. Thus, assuming the non-singularity of matrix $\JM$, we solve the unconstrained problem \eqref{eq:unconstrained_problem} via the Conjugate Gradient method (CG in short), looking for the solution of the linear system
\begin{equation}
\nabla \widetilde{\mathcal{J}}( \Jx ) = \JM \Jx + \Jd = \bm{0}.
\label{eq:conjugate_gradient_equation}
\end{equation}
\begin{comment}
\begin{algorithm}[t]
\caption{Preconditioned Conjugate Gradient method to solve $ \JM \Jx + \Jd = \bm{0}$.}\label{alg:PCG}
\KwData{The initial guess $\Jx_0$.}
\KwResult{$\Jx$.}
$\bm{r}_0 = \bm{\mathcal{M}} \mathcal{X}_0 + \bm{\mathit{d}}$\;
Solve $\bm{\mathcal{P}} \bm{z}_0 = \bm{r}_0$\;
Set $\delta \Jx_0 = - \bm{z}_0$ and $k = 0$\;
\While{PCG stopping criteria is not satisfied}{
$\alpha_k = \frac{\bm{r}_k^T \bm{z}_k}{\delta \mathcal{X}_k^T (\bm{\mathcal{M}}\delta \mathcal{X}_k)}$\;
$\mathcal{X}_{k+1} = \mathcal{X}_k + \alpha_k \delta \mathcal{X}_k$\;
$\bm{r}_{k+1} = \bm{r}_k + \alpha_k \bm{\mathcal{M}}\delta \mathcal{X}_k$\;
Solve $\bm{\mathcal{P}} \bm{z}_{k+1} = \bm{r}_{k+1}$\;
$\beta_{k+1} = \frac{\bm{r}_{k+1}^T \bm{z}_{k+1}}{\bm{r}_k^T \bm{z}_k}$\;
$\delta \mathcal{X}_{k+1} = - \bm{z}_{k+1} + \beta_{k+1}\delta \mathcal{X}_k$\;
$k = k+1$\;
}
\end{algorithm}
\end{comment}
Let us remark that the conjugate gradient scheme can be applied without building explicitly matrix $\JM$. Indeed, given a descent direction $\delta\Jx = \begin{bmatrix}
\addlinespace \delta\phis[]{h}\\
\addlinespace \delta\phix[]{h}
\end{bmatrix}$  it is possible to prove that
\begin{align*}
    \JM \delta \Jx &=   \begin{bmatrix}
    \addlinespace\GG\ \delta \phis[]{h} -(\ssDD)^T \delta \rph[]{h}
    -(\xsDD[\Lp])^T\delta \xph[]{d} \\
    \addlinespace \wGG\ \delta \phix[]{h}  - (\xxDD)^T  \delta \xph[]{h}
     +(\sxDD[\Lp])^T \delta \rph[]{d}
    \end{bmatrix},
\end{align*}
where $\delta \rph[]{h}\in\R^{\Nps}$ and $\delta \xph[]{h}\in\R^{\Npx}$ are the solutions of 
\begin{equation*}
    \AAstar \delta \rph[]{h} = \sxDD[\Lp] \delta \phix[]{h}, \qquad \wAAstar \begin{bmatrix} 
    \addlinespace \delta \uux[]{h} \\ 
    \addlinespace \delta \xph[]{h}
    \end{bmatrix} = \begin{bmatrix} 
    \addlinespace\mathbf{0} \\
    \addlinespace-\xsDD[\Lp] \ \delta \phis[]{h}
    \end{bmatrix},    
\end{equation*}
with $\delta \uux[]{h}\in \mathbb{R}^{\Nux}$. Similarly, 
$\delta \rph[]{d}\in \mathbb{R}^{\Nps},$ and $\delta \xph[]{d} \in \mathbb{R}^{\Npx}$ are the solutions of 
\begin{equation*}
\begin{gathered}
\AAstar \delta \rph[]{d}^{n,\ell} = \MM  \delta \rph[]{h}^{n,\ell} -  \ssDD  \delta \phis[]{h}^{n,\ell},\quad 
\wAAstar
\begin{bmatrix} 
\addlinespace \delta \uux[]{d}^{n,\ell}\\ 
\addlinespace \delta \xph[]{d}^{n,\ell} 
\end{bmatrix} =
\begin{bmatrix}
\mathbf{0}\\
 \wMM\ \delta \xph[]{h}^{n,\ell} - \xxDD\  \delta \phix[]{h}^{n,\ell}
\end{bmatrix},
\end{gathered}
\label{eq:subproblems}
\end{equation*}
with $\delta \uux[]{d}\in \mathbb{R}^{\Nux}$.
To reduce the number of iterations required by the conjugate gradient scheme a suitable preconditioning strategy could be devised. We refer to Section \ref{sec:TP3} in the numerical examples, in which a very simple, but effective preconditioner is tested on a large scale simulation.

\begin{algorithm}[t]
\caption{The Solving Strategy.}\label{alg:optAlgo}
\KwData{$\rph[]{h}^{0},\ \Lambdaj[0]$ \Comment*[r]{Initial conditions}}
\For{$j=1,\dots,J$ }{ 
$\Lambdaj[j] \gets \mathrm{Evolve}(\Lambdaj[j-1], \rph[]{h}^{j-1},\Theta^{j-1}, \nabla \Theta^{j-1})$ \;
\For{$n=1,\dots,N$ }{ 
$\ell \gets 0$\;
$ \rph[]{h}^{n} \gets  \rph[]{h}^{n-1}$ \;
\While{non-linear stopping criteria is not satisfied}{
Set an initial guess for control variables $\Jx_0$ \;
$\begin{bmatrix} \phis[]{h}^{n,\ell} \\ \phix[]{h}^{n,\ell} \end{bmatrix} = \Jx = \mathrm{CG}(\Jx_0)$ \Comment*[r]{Solve \eqref{eq:conjugate_gradient_equation} via the CG}  
Compute $ \rph[]{h}^{n,\ell}(\phix[]{h}^{n,\ell})$ and $\begin{bmatrix}
\uux[]{h}^{n,\ell}(\phis[]{h}^{n,\ell})\\
\xph[]{h}^{n,\ell}(\phix[]{h}^{n,\ell})
\end{bmatrix}$ from the systems in Eq. \eqref{eq:picard_constraints_1}-\eqref{eq:picard_constraints_2}\;
$\ell = \ell + 1$ \;
}
}
$\rph[]{h}^{0} \gets \rph[]{h}^{N, \ell}$ \;
}
\end{algorithm}

The whole solving process, accounting both for root growth and for the evolution of the quantities of interest is summarized in Algorithm \ref{alg:optAlgo}. The function $\mathrm{Evolve}(\Lambdaj[j-1], \rph[]{h}^{j-1},\Theta^{j-1}, \nabla \Theta^{j-1})$ updates the current root system architecture $\Lambdaj[j-1]$ with the information given by the discrete soil pressure head $\rph[]{h}^{j-1}$ and by the discrete water content $\Theta^{j-1}$ and its gradient at time $t_{j-1}$, according to the RSA growth strategy detailed in Section \ref{sec:RSA}. Since we are using a Virtual Element method to discretize the soil pressure head, we need to resort to the $L^2$-projection operator defined on the elements $E \in \Th$ to access the point-wise evaluation of $\rph[]{h}^{j-1}$. Thus, we actually approximate $ \rph[]{h}^{j-1}(\xx)$, $\Theta^{j-1}(\bm{x})$ and $\nabla \Theta^{j-1}(\xx)$ in a point $\xx \in \Omega$ with the weighted sums
\begin{gather*}
\rph[]{h}^{j-1}(\xx) \approx \displaystyle \frac{1}{\overline{\omega}(\xx)} \sum_{E \in \Th : \xx \in E} \omega_E \proj{0, E}{k} (\rph[]{h}^{j-1}(\xx)), \\
\Theta^{j-1}(\xx) \approx \displaystyle \frac{1}{\overline{\omega}(\xx)} \sum_{E \in \Th : \xx \in E} \omega_E\ \Theta\!\left(\proj{0, E}{k} (\rph[]{h}^{j-1}(\xx))\right),\\
 \nabla \Theta^{j-1}(\xx) \approx \displaystyle \frac{1}{\overline{\omega}(\xx)} \sum_{E \in \Th : \xx \in E} \omega_E\  \rc\!\left(\proj{0, E}{k} \rph[]{h}^{j-1}(\xx)\right) \proj{0, E}{k-1} (\nabla \rph[]{h}^{j-1}(\xx)),
\end{gather*}
where we choose as weights $\omega_E = \vert E \vert$, for each $E \in \Th$, and $\overline{\omega}(\xx) = \sum_{E \in \Th : \xx \in E} \omega_E$.

\begin{remark}[Numerical evaluation of the repulsion function]\label{rem:repulsionfunction}
The repulsion function \eqref{eq:repulsionfunction} is here treated as a VEM function of order $k=1$ over a tetrahedral mesh of $\Omega$. We observe that this kind of function is completely determined by assigning its values at the vertices of the tetrahedra. In particular, we assign the value $1$ at vertices belonging to the domain boundaries considered impenetrable and zero otherwise. As a consequence, the value of $d_{\max}$ in \eqref{eq:repulsionfunction} is determined by the size of the underlying mesh. The choice of building the repulsion function on a tetrahedral mesh is related to the fact that its VEM linear projection may assume negative values on differently shaped elements \cite{PintoreTeora2024}.
\end{remark}

\section{Numerical results}\label{sec:numericalresults}
In this section, we propose some numerical examples to validate the proposed approach. Given $E\in \Th$, we denote by $h_E:=\mathrm{diam}(E)$ the diameter of $E$, and we set $h=\max_{E \in \Th} h_E$. The refinement level of the 1D partitions is instead given in relation to the refinement level of the mesh induced on each $\Lambdaj_i$ by the intersection with the polyhedrons in $\Th$. We hence introduce three parameters $\hat{\delta}_i,\ \hat{\delta}^\phi_{\xylem,i}$ and $\hat{\delta}^\phi_{\soil,i}$, which express the ratio between the number of elements in $\xTh[i], \xxTh[i]$ and $\ssTh[i]$ and the number of elements of the mesh induced on $\Lambdaj_i$. For simplicity, in the proposed numerical examples we always consider equispaced partitions and we choose unique values of the 1D mesh parameters for the different root segments. In particular $\hat{\delta}_i=1,\ \mathrm{ and } \ \hat{\delta}^\phi_{\xylem,i}=\hat{\delta}^\phi_{\soil,i}=0.5$, $\forall i=1,...,J$. For what concerns the order of accuracy we choose $k=1$ for all the variables.

\subsection{Test Problem 1 (TP1): Convergence results}
Let us consider a cubic domain $\Omega=(-1,1)^3$ and let $\Sigmaj[ ]$ be a straight root segment crossing $\Omega$ from side to side. In particular, $\Sigmaj[ ]$ coincides with a cylinder of radius $R = 10^{-2}$ and height $2$ whose centerline lies on the $z$-axis . In this numerical test, we are not accounting for root growth: we consider instead a simplified configuration with a known analytical solution, to analyze accuracy.

The considered test problem is adapted from  Test Problem 1 in \cite{Grappein2023}, which is modified by introducing time dependency and non-linearity. In particular, we aim at solving the following 3D-1D coupled problem \begin{equation}\label{probTest1}
\begin{cases}
\scal[{\Dj[ ]}]{\rc({\rph{}})\frac{\partial \rph{}}{\partial t}}{\qs{}}  + \scal[{\Dj[ ]}]{\rK(\rph{}) \nabla \rph{}}{\nabla \qs{}} + \scal[{\Lambdaj[ ]}]{ 2 \pi R \Lp(\rph{a} - \xph{a}) }{\qs{a}} = \scal[{\Dj[ ]}]{\mathcal{S}}{\qs{}}&\forall \qs{} \in \Qsj,\\
\\
\addlinespace \scal[{\Lambdaj[ ]}]{ \kx \ux }{\vx} - \scal[{\Lambdaj[ ]}]{ \xph{a} }{ \pi R^2 \frac{\partial \vx}{\partial s}} = \pi R^2 \Big[(\widehat{\mathcal{S}}_{\psi}  \vx)_{|[0,0,1]^T}  - (\widehat{\mathcal{S}}_{\psi}  \vx)_{|[0,0,-1]^T}\Big] & \forall \vx \in \Vrxj,\\
\addlinespace \scal[{\Lambdaj[ ]}]{\pi R^2 \frac{\partial \ux}{\partial s}}{ \qx{a}} + \scal[{\Lambdaj[ ]}]{ 2\pi R \Lp (\xph{a} - \rph{a})}{\qx{a}} = \scal[{\Lambdaj[ ]}]{ \widehat{\mathcal{S}}}{\qx{a}}& \forall \qx{a} \in \Qrj,\\
\end{cases}
\end{equation} with
\begin{equation*}
\rc(\rph{}) = \frac{-\rph{}}{(1 + \rph{}^2)^{3/2}} + 4,\quad \rK(\rph{}) = 1,\quad \kx = \pi R^2 \left(\frac{z^2}{3} + \frac{1}{2}\right)^{-1}, \quad \Lp=\frac{2R}{R^2+2}.
\end{equation*} 

Problem \eqref{probTest1} is actually obtained from \eqref{eq:3D1Dproblem1}-\eqref{eq:constr1D_2} by neglecting the gravity terms, introducing proper forcing terms $\mathcal{S}$, $\widehat{\mathcal{S}}$, imposing boundary Dirichlet conditions $\widehat{\mathcal{S}}_{\psi}$ for the xylem pressure head, and imposing initial and Dirichlet boundary conditions for the soil pressure head such that the exact solution of \eqref{probTest1} is 
\begin{gather*}
\rph{} = \frac{1}{2}(x^2 + y^2)(1 - z^2) - 1 -t,\\
\xph{a} = z^2 - 2 -t,\qquad \ux = - 2 z \left(\frac{z^2}{3} + \frac{1}{2}\right).
\end{gather*}  
The accuracy of the method is measured by the following error indicators:
\begin{equation}
    \begin{gathered}
        e_{L^2, \rph{}}^2 =  \frac{ \sum_{E \in \Th} \norm[E]{\rph{} - \proj{0}{k} \rph[]{h}}^2 }{\norm[\Omega]{\rph{}}^2}, \quad   e_{L^2, \xph{a}}  = \frac{ \norm[{\Lambdaj[ ]}]{\xph{a} -  \xph[]{h}} }{\norm[{\Lambdaj[ ]}]{\xph{a}}}, \\
        e_{H^1, \rph{}}^2 =  \frac{ \sum_{E \in \Th} \norm[E]{\nabla\rph{} - \proj{0}{k-1}\nabla \rph[]{h}}^2 }{\norm[\Omega]{\nabla \rph{}}^2},  \quad  e_{L^2, \ux}  = \frac{ \norm[{\Lambdaj[ ]}]{\ux -  \uux[]{h}} }{\norm[{\Lambdaj[ ]}]{\ux}},
    \end{gathered}
    \label{eq:errors_variable}
\end{equation}
corresponding, respectively, to the relative errors in $L^2$-norm for the soil and xylem pressure head, in $H^1$-seminorm for the soil pressure head and in $L^2$-norm for the xylem velocity.  Concerning the interface variables, we introduce two additional error indicators
\begin{equation}
     e_{L^2, \rph{a}}  = \frac{ \norm[{\Lambdaj[ ]}]{\tg \rph{} -  \phis[]{h}} }{\norm[{\Lambdaj[ ]}]{\tg \rph{} }}, \quad 
    e_{L^2, \xph{a}}  = \frac{ \norm[{\Lambdaj[ ]}]{\xph{a} -  \phix[]{h}} }{\norm[{\Lambdaj[ ]}]{\xph{a}}},
    \label{eq:errors_control}
\end{equation}
which measure the $L^2$-distance of the control variables from the trace of the corresponding pressure head variable. We consider 4 different tetrahedral meshes, with mesh size $h \approx 4.67e\text{-}01, 3.56e\text{-}01, 2.82e\text{-}01, 2.33e\text{-}01$, and we solve the problem for $t \in (0,1]$.

\begin{figure}[!ht]
    \begin{subfigure}[t]{0.45\textwidth}
            \centering
            \includegraphics[width=\textwidth]{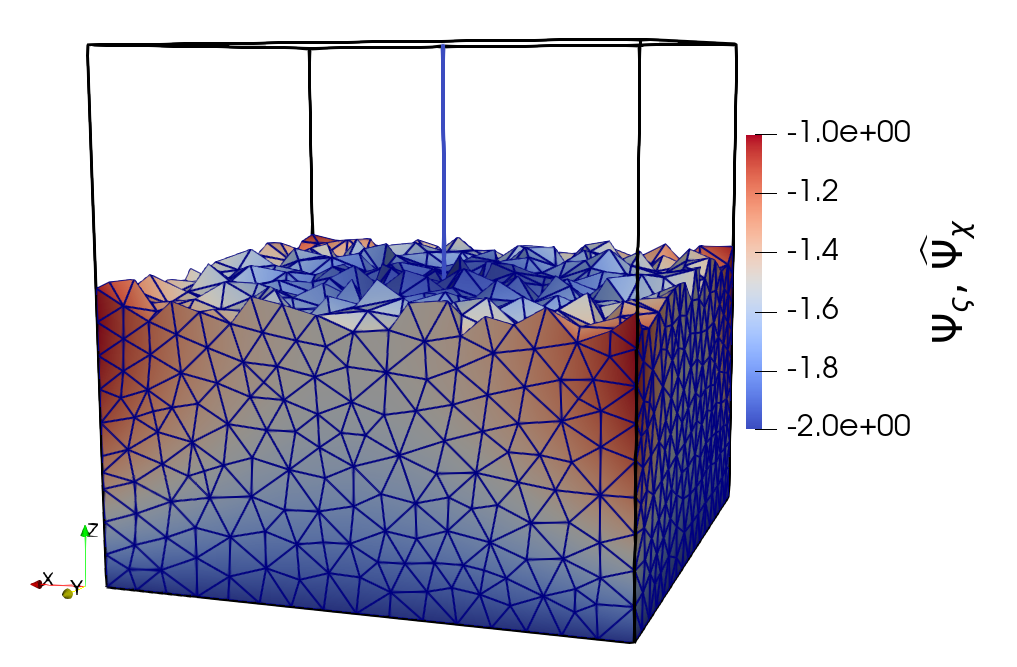}
            \caption{Numerical pressure head. The view is obtained by a crinkle cut of the 3D domain with a plane orthogonal to the $z$-axis).}
            \label{fig:test1_press}
    \end{subfigure}\hfill
    \begin{subfigure}[t]{0.45\textwidth}
            \centering
            \includegraphics[width=\textwidth]{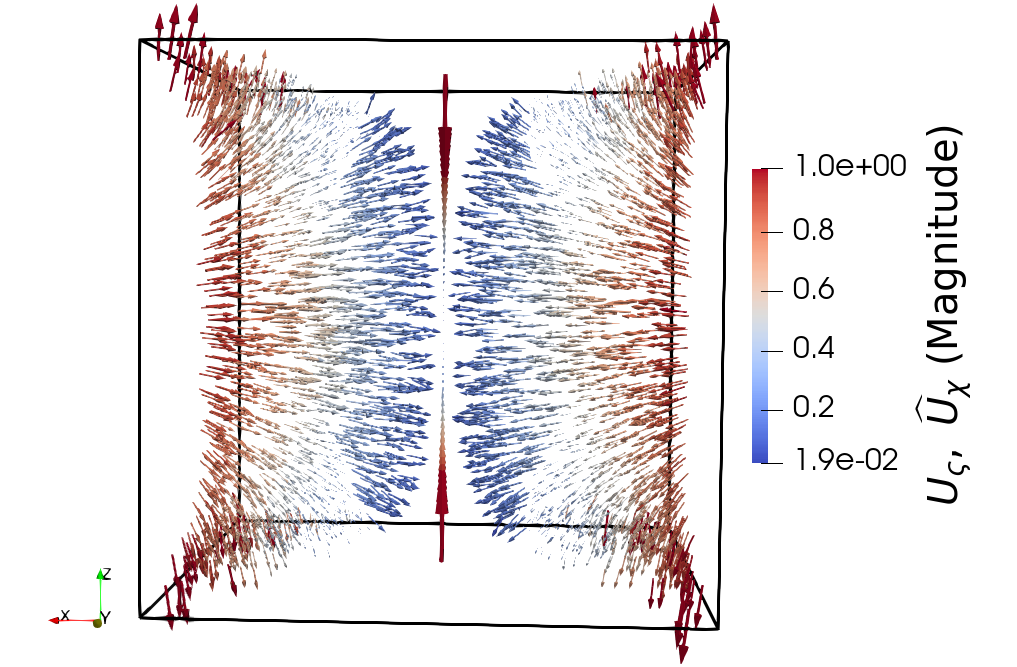}
            \caption{Numerical velocity field.}
            \label{fig:test1_vel}
    \end{subfigure}
    \caption{TP1: Numerical solution in the 3D soil sample and in the root segment at the final time $t=1$ on the finest mesh.}
    \label{fig:test1_sol}
\end{figure}

Figure \ref{fig:test1_sol} shows the numerical solution of the problem obtained at $t=1$ on the finest mesh. Figure \ref{fig:test1_press} reports the pressure head distribution in the soil sample and in the root segment, and it is obtained by a crinkle cut of the soil sample with a plane orthogonal to the $z$-axis to show the pressure distribution in the interior of the domain; \ref{fig:test1_vel} shows, instead, the velocity field, where the 3D velocity variable is reconstructed as
\begin{equation*}
\uux{h} \approx - \proj{0,E}{k-1} \nabla \rph[]{h} \quad \forall E \in \Th,
\end{equation*}
provided that $\rK(\rph{}) \equiv 1$ and that the gravity term is neglected.

\begin{figure}[!ht]
    \begin{subfigure}[b]{0.327\textwidth}
            \centering
            \includegraphics[width=\textwidth]{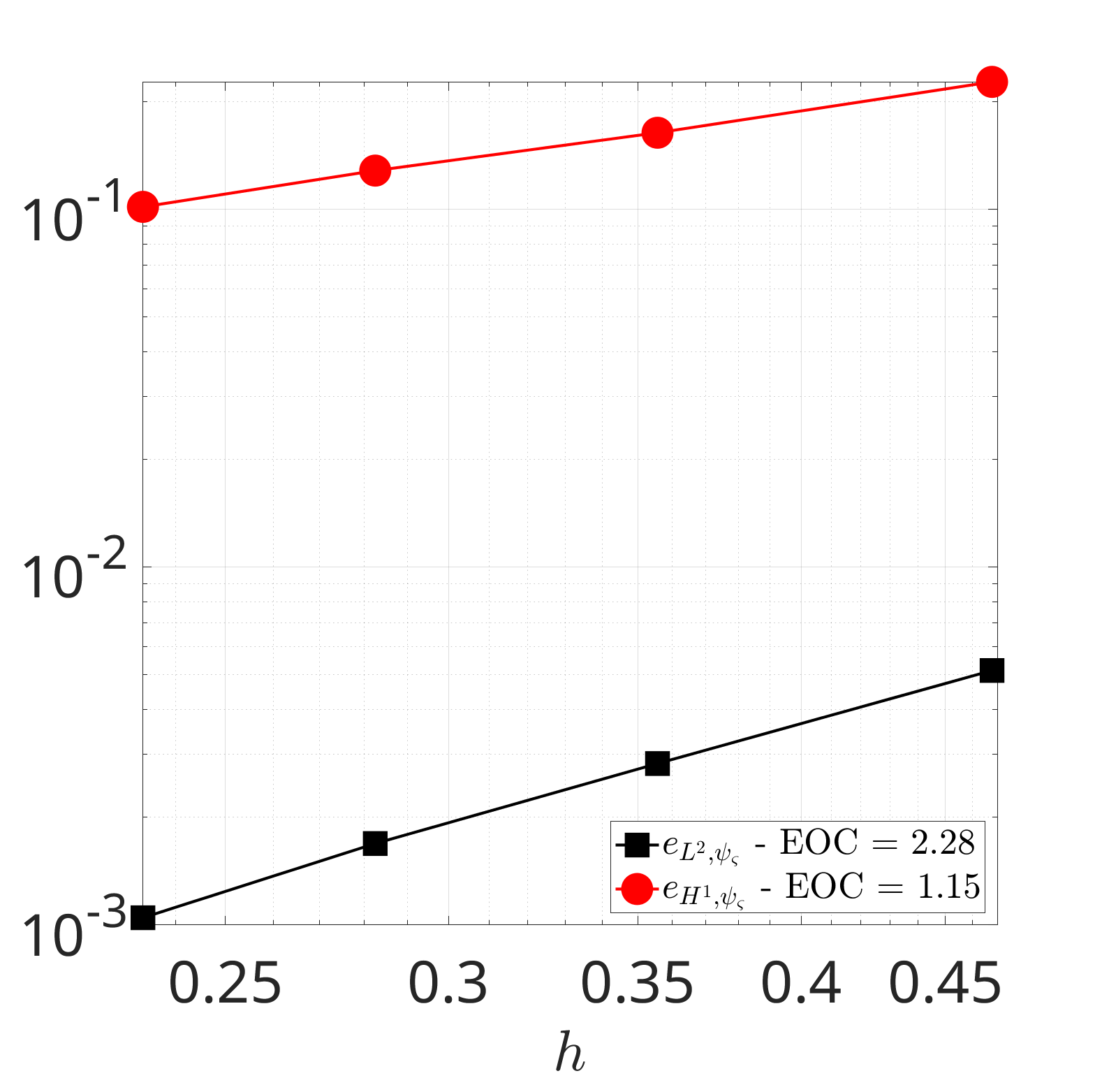}
            \caption{}
            \label{fig:test1_nl_soil_errors}
    \end{subfigure}\hfill
    \begin{subfigure}[b]{0.32\textwidth}
            \centering
            \includegraphics[width=\textwidth]{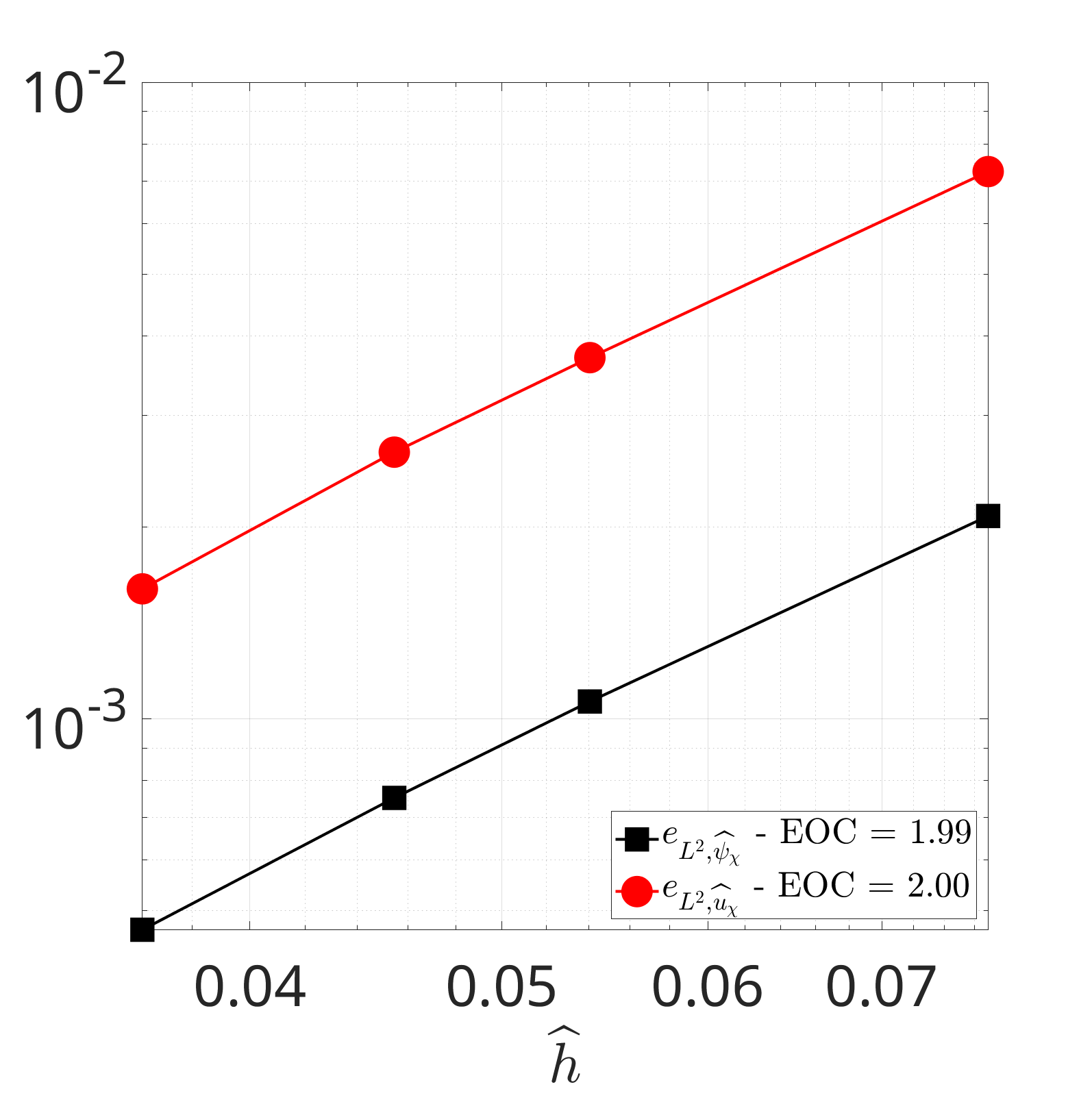}
            \caption{}
            \label{fig:test1_nl_xylem_errors}
    \end{subfigure}\hfill
        \begin{subfigure}[b]{0.327\textwidth}
            \centering
            \includegraphics[width=\textwidth]{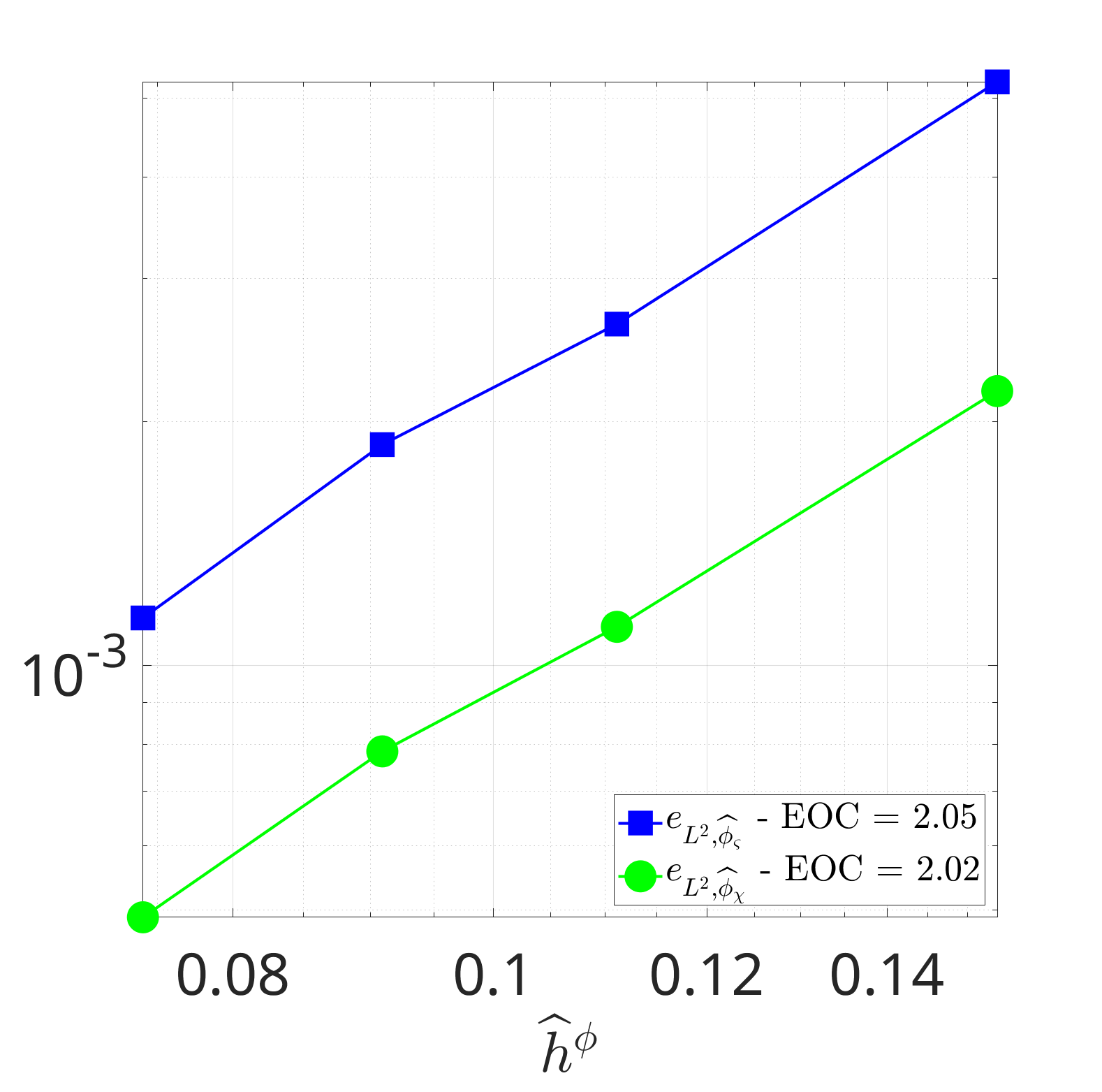}
            \caption{}
            \label{fig:test1_nl_control_errors}
    \end{subfigure}
    \caption{TP1: Behaviour of the errors \eqref{eq:errors_variable}-\eqref{eq:errors_control} as the related mesh parameter decreases at the final time $t=1$.}
    \label{fig:test1_errors}
\end{figure}

Figure \ref{fig:test1_errors} displays the convergence curves for the error indicators \eqref{eq:errors_variable} and \eqref{eq:errors_control} as the related mesh parameters vary. In addition, Figure \ref{fig:test1_errors} reports the Empirical Order of Convergence (EOC in short) measured for each variable. Figure \ref{fig:test1_nl_soil_errors} refers to the soil error indicators, which appear to converge at the expected rate. Indeed, since we are using a primal formulation for the soil pressure variable, $ e_{L^2, \rph{}}$ is expected to decay as $O(h^2)$, whereas $ e_{H^1, \rph{}}$ as $O(h)$.
For what concerns the xylem error indicators, let us observe that, since the 1D mesh parameters $\hat{\delta},\ \hat{\delta}^\phi_{\xylem},\ \hat{\delta}^\phi_{\soil}$ are kept fixed, a finer 3D mesh induces finer 1D partitions. In particular we denote by $\hat{h},\ \hat{h}^\phi_\xylem$ and $\hat{h}^\phi_\soil$ the width of the equispaced intervals respectively of $\xTh,\ \xxTh$ and $\ssTh$ and we analyze the convergence of $e_{L^2, \xph{a}}$ and $e_{L^2, \ux}$ with respect to $\hat{h}$, and of $e_{L^2, \xph{a}}$ and $e_{L^2, \rph{a}}$ with respect to $\hat{h}^\phi_\xylem$ and $\hat{h}^\phi_\soil$, respectively. Actually, since we are considering $\hat{\delta}^\phi_{\xylem}=\hat{\delta}^\phi_{\soil}$ we here have $\hat{h}^\phi:=\hat{h}^\phi_\soil=\hat{h}^\phi_\xylem$. Also in this case we can observe how xylem pressure head and velocity converge at the expected rate: indeed we expect both $e_{L^2, \xph{a}}$ and $e_{L^2, \ux}$ to converge as $O(\hat{h}^2)$. For what concerns the control variables, we have no theoretical results. However, we observe an empirical convergence rate equal to $2$ for both variables, which is quite predictable since they represent the finite element interpolation of their continuous counterparts.

Concerning the time-advancing scheme, let us remark that, since the exact pressure variables are linear in time,  the proposed results can be obtained by performing just one backward Euler step. Concerning instead the non-linearities and the conjugate gradient algorithm, we want to observe that, for each mesh refinement,  Algorithm \ref{alg:optAlgo} stops after $8$ non-linear iterations $\ell$ and, as $\ell$ increases, we observe a decrease in the number of CG iterations required to reach convergence, as shown in Table \ref{tab:test1_CGstep} for the finest mesh.

\begin{table}[!ht]
\centering
\begin{tabular}{@{}ccccccccc@{}}
\toprule
$\bm{\ell}$ & 0  & 1  & 2 & 3 & 4 & 5 & 6 & 7 \\ \midrule
\textbf{\# CG}    & 22 & 20 & 16 & 12 & 9 & 5 & 2 & 0 \\ \bottomrule
\end{tabular}
\vspace{4pt}
\caption{TP1: Number of CG iterations for each non-linear iteration $\ell$ for the finest mesh.}
\label{tab:test1_CGstep}
\end{table}

\subsection{Test Problem 2 (TP2): Root Architecture development in bounded soil sample}
For this second numerical example we consider a small soil sample $\Omega=(0,B)^2\times (-H,0)$, with $B=$3~cm and $H=$6~cm. We assume the seed to be in position $\bm{x}_s=(1.5,1.5,0.1)$ cm, i.e. below the top surface of the soil sample, and to be connected to the upper face by a vertical mesocotyl. We do not aim at reproducing a specific plant genotype, but at testing the proposed 3D-1D coupling approach along with the RSA growth in a controlled setting, in which all the conditions are optimal for root growth. 

We consider a uniform hexaedral mesh with size $h=0.15$~cm, whereas we simulate root growth in a time span of 9 days, divided into uniform time intervals of width $\Delta \mathcal{I}_j=0.2$ days. For the time discretization of the constraint equations we are not further subpartitioning the time steps $\Ij$ used for root growth, i.e. we choose $\Delta t=\Delta \mathcal{I}_j$, $\forall j \geq 1$ (see Section \ref{sec:time_discretization}).

We assume the soil sample to be bounded into a box whose walls are impenetrable for the roots. From a practical standpoint, we treat the walls as an obstacle, defining a linear repulsion function $\repulsionfunction$ as the one in \eqref{eq:repulsionfunction}. In particular we set $\repulsionfunction(\bm{x})=1$ for $\bm{x} \in \partial \Omega$ and we assume $\repulsionfunction$ to linearly decay from 1 to 0 within the first layer elements, i.e. at a distance $h$ from the faces of $\Omega$.

\begin{figure}[!ht]
    \centering
    \includegraphics[width=0.6\textwidth]{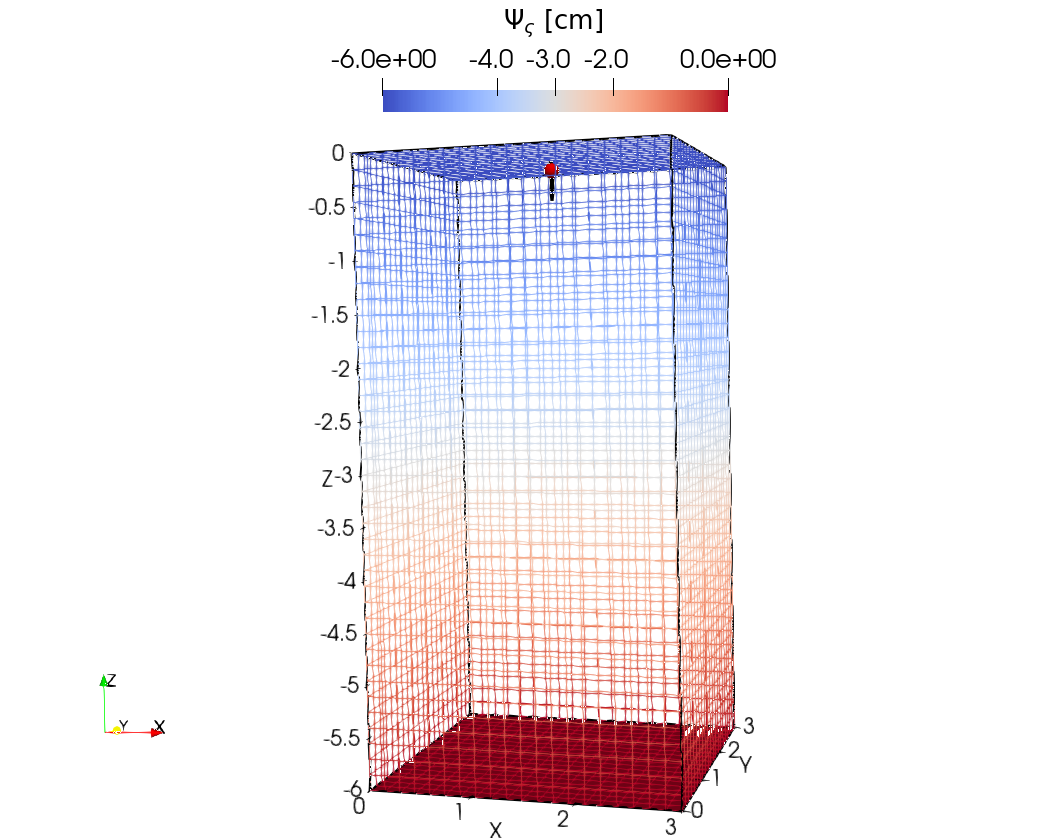}
    \caption{TP2: Hexaedral mesh and (initial) soil pressure head distribution. The red face corresponds to the water table, whereas the red dot corresponds to the seed position $\bm{x}_S$.}
    \label{fig:initial_configtp2}
\end{figure}

Similarly to \cite{Jin2020}, at the initial time $t_0 = 0$, we consider a linear soil pressure head increasing from $\rph{h}=-6\ \mathrm{cm}$ at the top to $\rph{h}=0\ \mathrm{cm}$ at the soil bottom, as reported in Figure \ref{fig:initial_configtp2}. For each time $t > 0$ we set no flux boundary conditions on all the sample faces except for the bottom one, on which we keep a constant water table $\rph{h}=0\ \mathrm{cm}$. At the root collar we set a constant transpiration rate of $0.2\ \mathrm{cm}^3/ \mathrm{day}$ .

At each time-step $\Delta \Ij$ we solve a 3D-1D coupled problem in the form of \eqref{eq:3D1Dproblem1}-\eqref{eq:constr1D_2} applying the solving strategy described in Section \ref{sec:solving_strategy}. In particular, we choose the Van Genuchten-Mualem model \cite{VanGenuchten1980}
\begin{equation}
\wc(\rph{}) = \begin{cases}
\wc_r + \frac{\wc_s - \wc_r}{[1 + (a \vert \rph{} \vert)^n]^m} & \text{if }  \rph{} < 0,\\
\wc_s & \text{if }  \rph{} \geq 0,
\end{cases}\quad \text{with } m = 1 - \frac{1}{n},\label{eq:theta_exp}
\end{equation}
\begin{equation}
\rK(\rph{}) = \begin{cases}
\rK_s  \frac{\left(1 - (a \vert \rph{}\vert)^{n-1}[1 + (a \vert \rph{} \vert)^n]^{-m}\right)^2}{[1 + (a \vert \rph{} \vert)^n]^{\frac{m}{2}}} & \text{if } \rph{} < 0,\\
\rK_s & \text{if }  \rph{} \geq 0,
\end{cases}\label{eq:K_exp}
\end{equation}
\begin{equation}
\rc(\rph{}) = \begin{cases}
anm(a\vert \rph{}\vert)^{n-1} \frac{\wc_s - \wc_r}{[1 + (a \vert \rph{} \vert)^n]^{m+1}}  & \text{if } \rph{} < 0,\\
0 & \text{if }  \rph{} \geq 0.
\end{cases}\label{eq:C_exp}
\end{equation}
The value assigned to the parameters involved in \eqref{eq:theta_exp}-\eqref{eq:C_exp} and in \eqref{eq:3D1Dproblem1}-\eqref{eq:constr1D_2} are summarized in Table \ref{tab:pde_parameter_tp2}.

\begin{table}[!ht]
    \renewcommand*{\arraystretch}{1.2}
    \centering
    \begin{tabular}{ccclc}
    \hline
    \textbf{Parameter} & \textbf{Value} &\textbf{Unit} &\textbf{Description} &\textbf{Reference}\\
    \hline
    $a$&0.03&1/cm& Water retention curve shape parameter&\cite{Jin2020}\\
    $n$&2.5& -& Water retention curve shape parameter&\cite{Jin2020}\\
    $\wc_r$&0.06&-&Residual volumetric water content&\cite{Jin2020}\\
    $\wc_s$&0.41&-&Saturated volumetric water content&\cite{Jin2020}\\
    $K_s$&10.24&cm/day&Saturated hydraulic conductivity&\cite{Jin2020}\\
    $R$&0.05&cm&Root radius&\cite{javaux2008}\\
    $\kx$&0.18&day/cm&Reciprocal of specific root axial conductance&\cite{Jin2020,schnepf2023}\\
    $\Lp$&$1.36\cdot 10^{-6}$&cm$^{2}$/day&Permeability of root wall &\cite{schnepf2023}\\
    $\sigma_{\mathrm{max}}$&$1$&MPa &Soil strength parameter & \\
    \hline
    \end{tabular}
    \vspace{4pt}
    \caption{TP2: PDE parameters used for simulation.}
    \label{tab:pde_parameter_tp2}
\end{table}

For what concerns the development of the RSA, we allow for a maximum of 3 root orders and, as mentioned before, we assume to be in optimal growth conditions: we set $\mathrm{Imp}(\rph{})=1$ for all the values of $\rph{}$ while, given the small pressure head absolute values, the soil strength impedance tends automatically to 1, according to \eqref{eq:theta_exp}, \eqref{eq:soil_strength} and \eqref{eq:soil_strength_impedance}. In this simplified setting, we assume that the maximum root elongation rate $V_a$, defined in \eqref{eq:growth_rate}, is constant in time and that it depends only on the root order, i.e. we set 
\begin{equation*}
    V_a(t) =\overline{V}_a(\omega) \quad \forall t.
\end{equation*}
with $\omega\in \lbrace0,1,2\rbrace$. The values of all the parameters related to root growth are reported in Table \ref{tab:parameter_tp2}. 

\begin{table}[!ht]
\renewcommand*{\arraystretch}{1.2}
\centering
\begin{tabular}{cccccl}
\hline
\multirow{2}{*}{\textbf{Parameter}} & \multicolumn{3}{c}{\textbf{Value}}   & \multirow{2}{*}{\textbf{Unit}}                           & \multirow{2}{*}{\textbf{Description}}                            \\
\cline{2-4} & \textbf{0} & \textbf{1} & \textbf{2} & &   \\ 
\hline
$\overline{V}_a$&1 & 0.8 &0.6 &cm& Growth rate\\
$L_B$&0.8& 0.8 & - & cm & Length of basal non-branching zone\\
$L_A$&2&0.5 & - &cm&Length of apical non-branching zone\\
$I$&1&0.4&-&cm&Inter-branch distance\\
$\alpha_I$&1.4&1.2&-&rad&Branching insertion angle\\
$X$&5&3&-&-&Number of xylem poles\\
$b_c$&1&1&-&-&Branching probability parameter\\
$k_g$&0.1&0.1&0.1&-&Geotropism weight\\
$k_s$&0.25&0.25&0.25&-&Hydrotropism weight\\
$k_w$&1&1&1&-&Exotropism weight\\
$m_a$&&[0,1]&&-&Space exploration parameter\\
\hline
\end{tabular}
\vspace{4pt}
\caption{TP2: RSA growth parameters used for simulation.}
\label{tab:parameter_tp2}
\end{table}

\begin{figure}[!ht]
    \centering
    \begin{subfigure}{0.33\textwidth}
        \centering
        \includegraphics[width=1\textwidth]{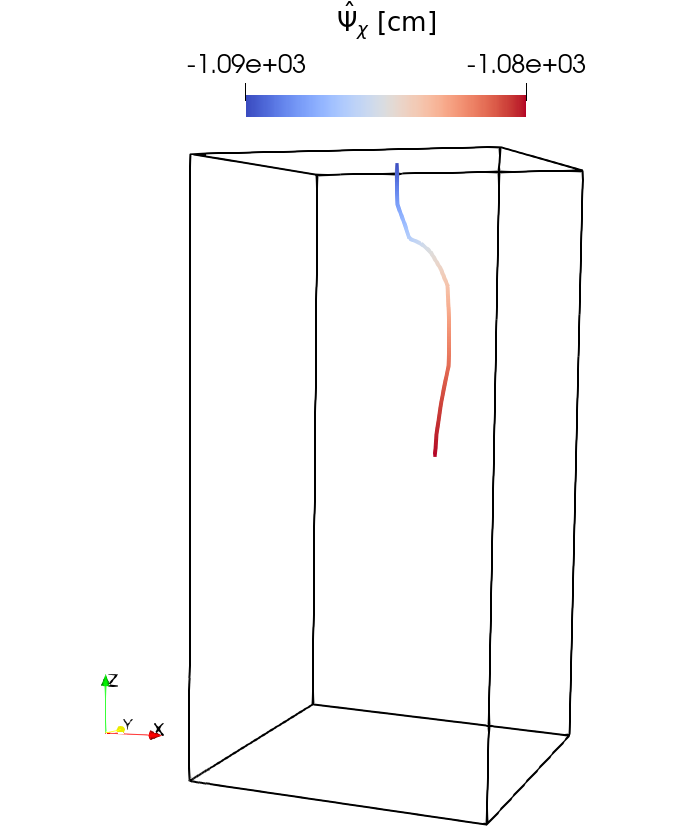}
        \caption{$t=3$ days.}
        \label{}
    \end{subfigure}
    \begin{subfigure}{0.33\textwidth}
        \centering
        \includegraphics[width=1\textwidth]{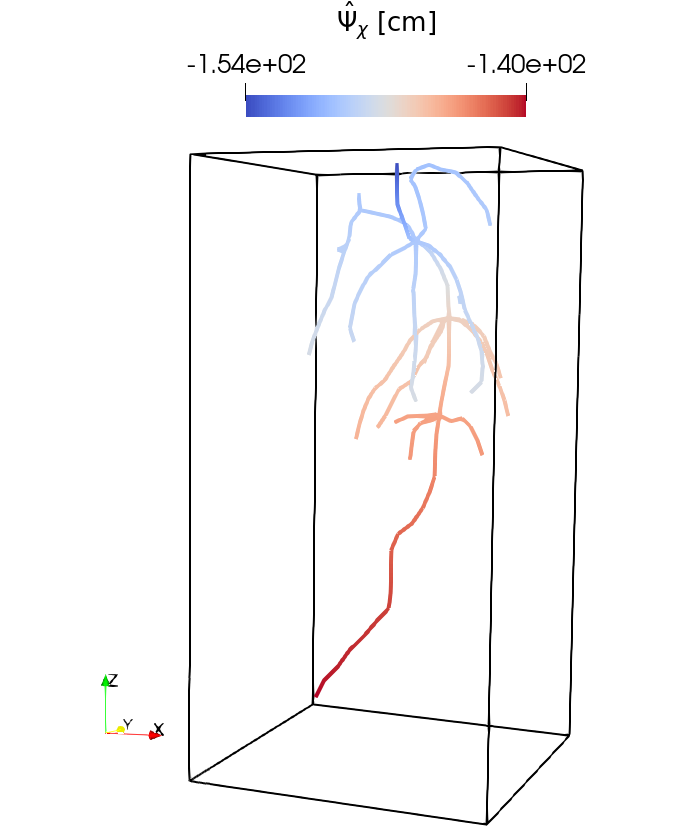}
        \caption{$t=6$ days.}
        \label{}
    \end{subfigure}
    \begin{subfigure}{0.33\textwidth}
        \centering
        \includegraphics[width=1\textwidth]{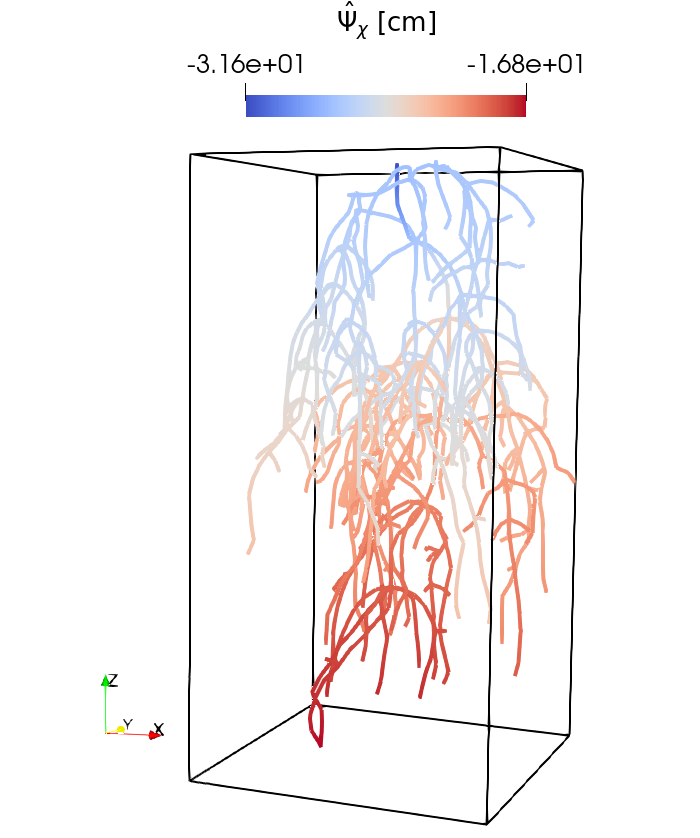}
        \caption{$t=9$ days.}
        \label{fig:Tp2_growth-c}
    \end{subfigure}
    \caption{TP2: RSA at three different time instants.}
    \label{fig:Tp2_growth}
\end{figure}

Figure \ref{fig:Tp2_growth} reports the root network at three different time instants, namely for $t=3,6$ and $9$ days, together with the xylem pressure head distribution. In Figure~\ref{fig:Tp2_growth-c} the effect of the repulsion function on the growth direction is particularly evident close to the bottom face of the domain, where the thigmotropic behavior of the primary root can clearly be observed. 

Concerning the soil variables, we notice from Figure~\ref{fig:Tp2_vel} that the velocity field is oriented in the radial direction around the roots, pointing towards the roots themselves. The figure shows the magnitude and orientation of the velocity field close to the primal root at $t=3$~days. Since the soil flux is entering in roots, we can assert that the root water uptake $2R\pi\Lp(\tg \rph{}-\tg \xph{})$ is positive at that stage. However, one of the strengths of the optimization-based approach is that the approximation of the interface quantities $\tg \rph{},\ \tg \xph{}$ is directly computed and stored in the auxiliary variables $\phis{h},\ \phix{h}$, so that the water uptake along the RSA can easily be computed without the need of post-processing, at any time of the simulation. 

Figure \ref{fig:Tp2_uptake}-left shows the magnitude of the water uptake at $t=9$~days.
Due to the constant transpiration rate imposed, we observe that the roots near the root collar are primarily responsible for water uptake. This is a consequence of the fact that the permeability $\Lp$ is assumed to be the uniform across all root segments. Recent works have shown that older root segments can play an important role in water uptake \cite{Cuneo2018}. However, in many scenarios, this process mostly occurs at the youngest root segments (the ones closest to the root tips), while the older ones are mainly responsible for the transport of water towards the collar. To reflect this behavior in the model, $\Lp$ can be chosen as a decreasing function of time, and in particular of the age of the root segments. For simplicity, we here propose a case in which we set  $\Lp=0$ for all the root segments $\Lambda_i^j$ whose age $a_i^j$ is greater than 3 days, i.e. for all the segments that have been existing for more than 3 days from the beginning of the simulation. Conversely, we impose $\Lp$ as specified in Table \ref{tab:pde_parameter_tp2} for root segments with $a_i^j< 3\ \mathrm{days}$. Figure \ref{fig:Tp2_uptake}-right shows the water uptake resulting from this last setting. As expected, this adjustment increases the water uptake in the youngest root segments, which now exclusively sustain the imposed transpiration rate. On the other hand, older root segments, marked in purple in Figure \ref{fig:Tp2_uptake}-right, no longer contribute to water uptake. It is worth noting that smoother dependencies of $\Lp$ on segment age could be employed to better capture specific experimental observations.

\begin{figure}[!ht]
    \centering
    \savebox{\imagebox}{\includegraphics[width=0.35\textwidth]{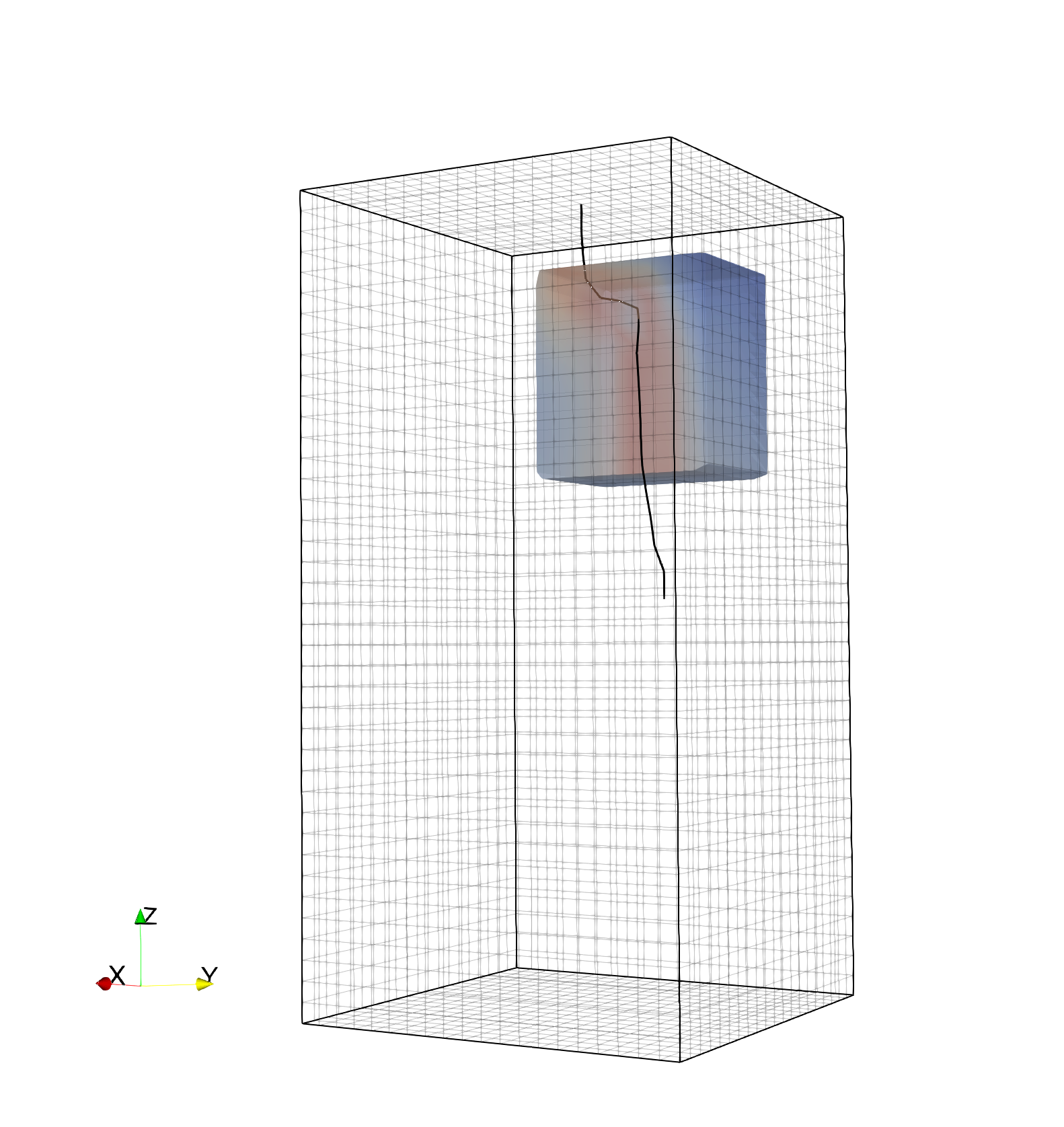}}%
    \begin{subfigure}[t]{.4\linewidth}
    \centering\usebox{\imagebox}% Place largest image
    \end{subfigure}
    \begin{subfigure}[t]{.4\linewidth}   \centering\raisebox{\dimexpr.5\ht\imagebox-.35\height}{% Raise smaller image into place
    \includegraphics[width=0.7\textwidth]{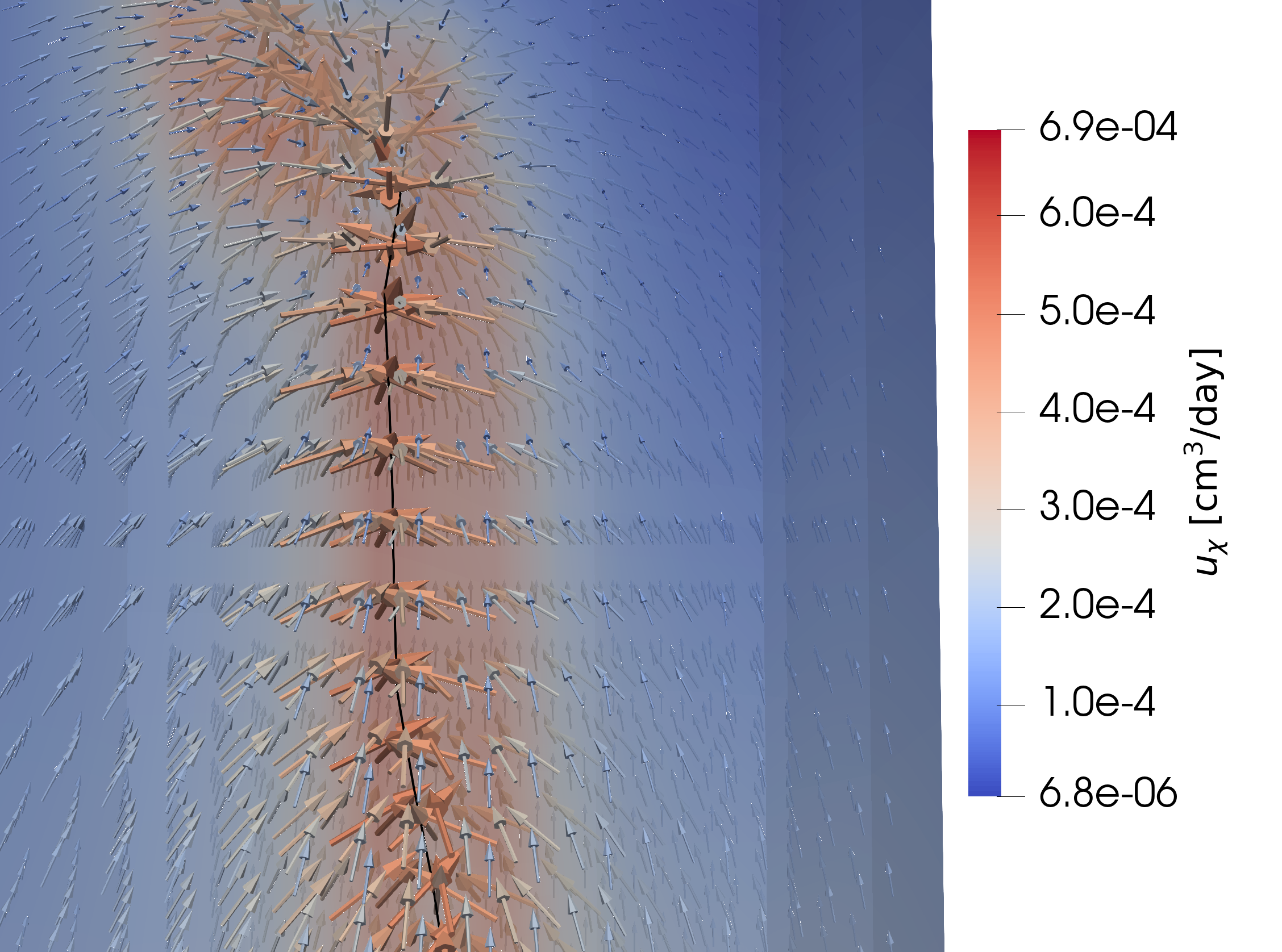}}%
    \end{subfigure}
    % Draw an arrow between the two subfigures
    \begin{tikzpicture}[overlay, remember picture]
    \draw[->] (-9,4) to [bend left] (-5.8,4);
    \end{tikzpicture}
    \caption{TP2: Detail on 3D velocity field at $t=3$ days in proximity of the RSA.}
    \label{fig:Tp2_vel}
    \end{figure}
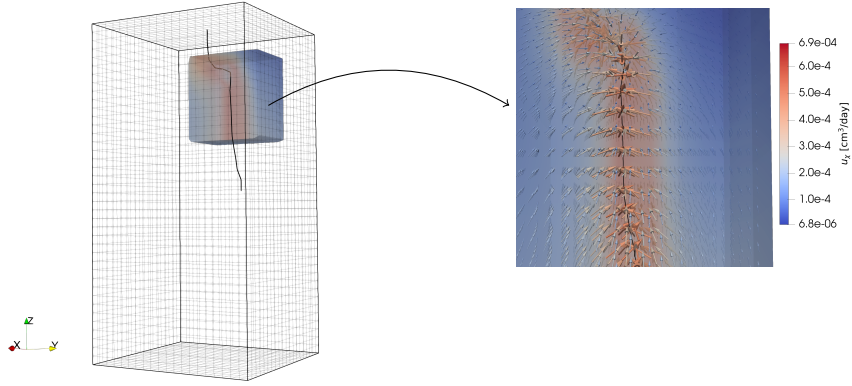
    \begin{figure}[t]
    \centering
    \includegraphics[width=0.8\textwidth]{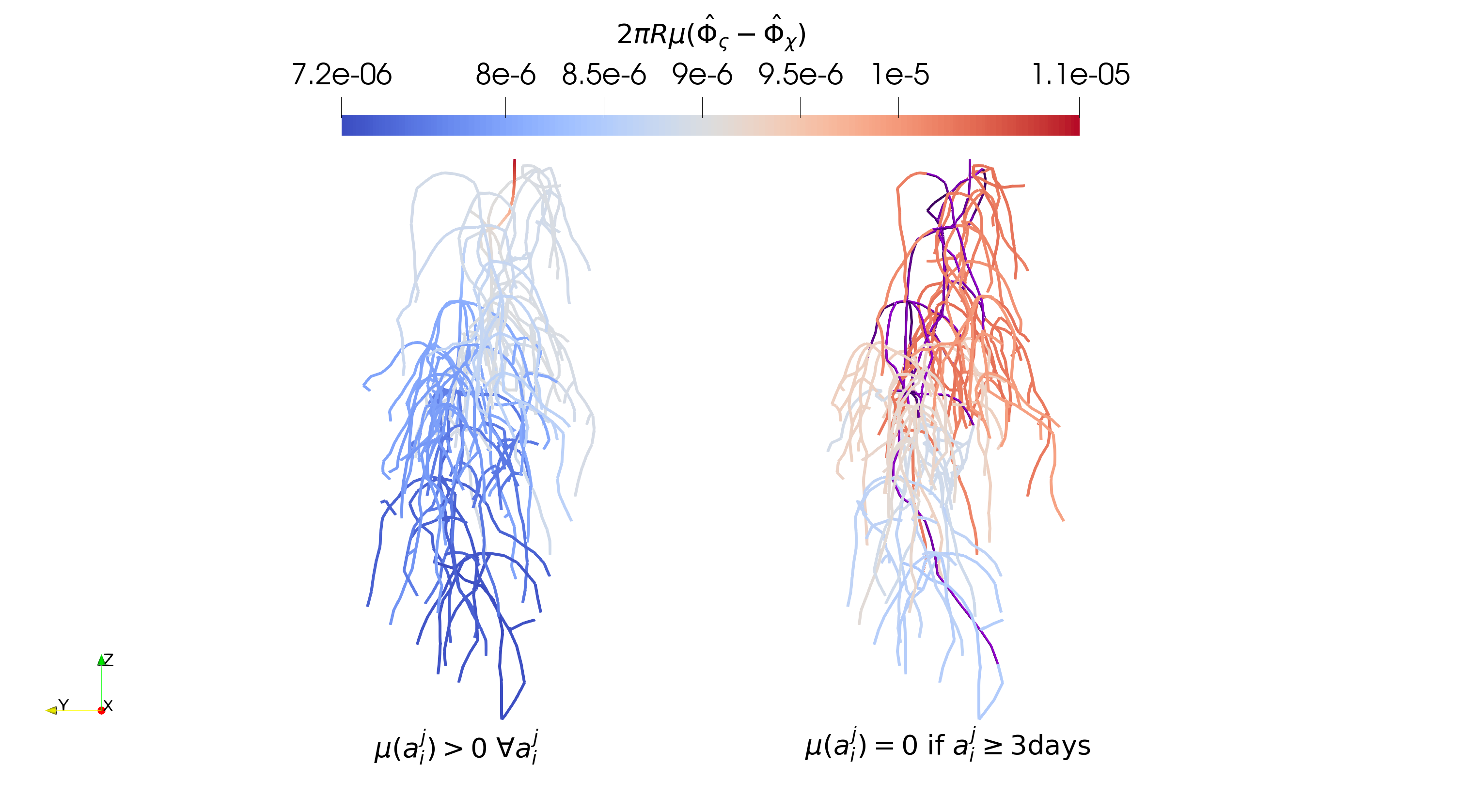}
    \caption{TP2: Water uptake along the RSA. On the left: uptake occurs along the whole root network. On the right: uptake occurs only at root segments which are less than 3 days old; older root segments are marked in purple.}
    \label{fig:Tp2_uptake}
\end{figure}

\subsection{Test Problem 3 (TP3): Root Architecture development in stony soil} \label{sec:TP3}
\begin{figure}[ht]
	\begin{subfigure}[t]{0.5\textwidth}
		\centering
		\includegraphics[width=\textwidth]{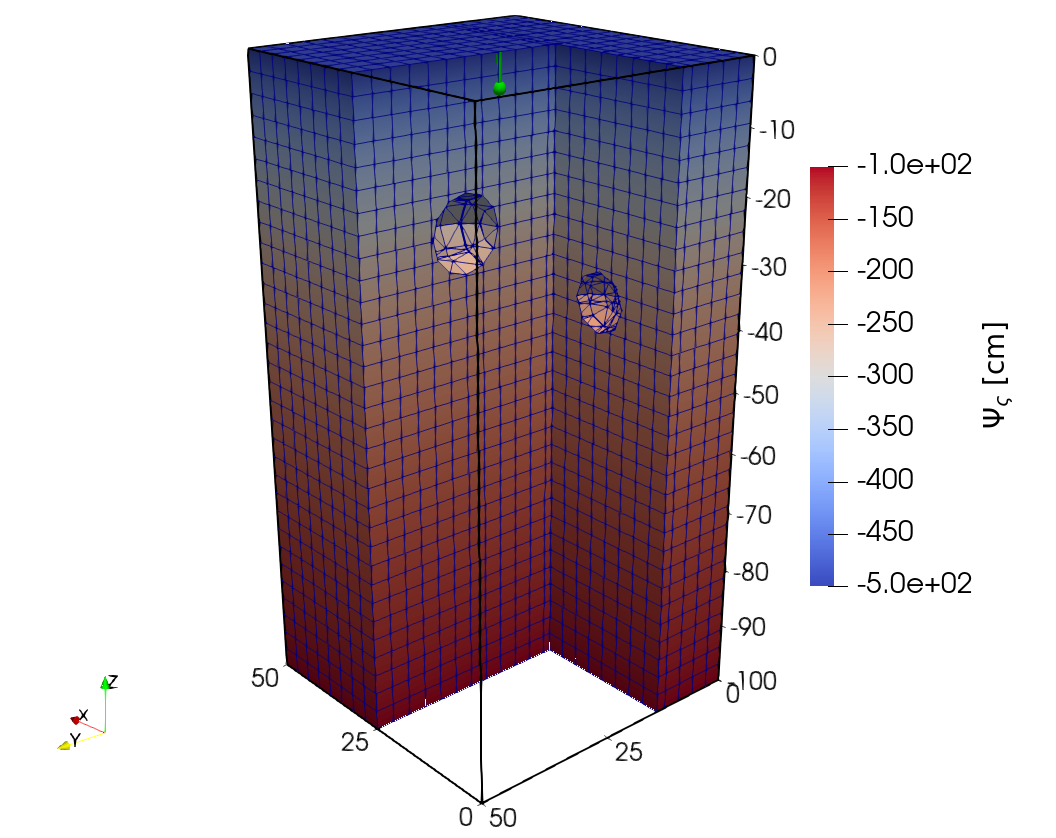}
		\caption{Clipped view of the 3D mesh. Color scale refers to the initial soil pressure head distribution. Seed and mesocotyl are colored in green.}
		\label{fig:stony_test:mesh:detail}
	\end{subfigure}
	\begin{subfigure}[t]{0.4\textwidth}
		\centering
		\includegraphics[width=\textwidth]{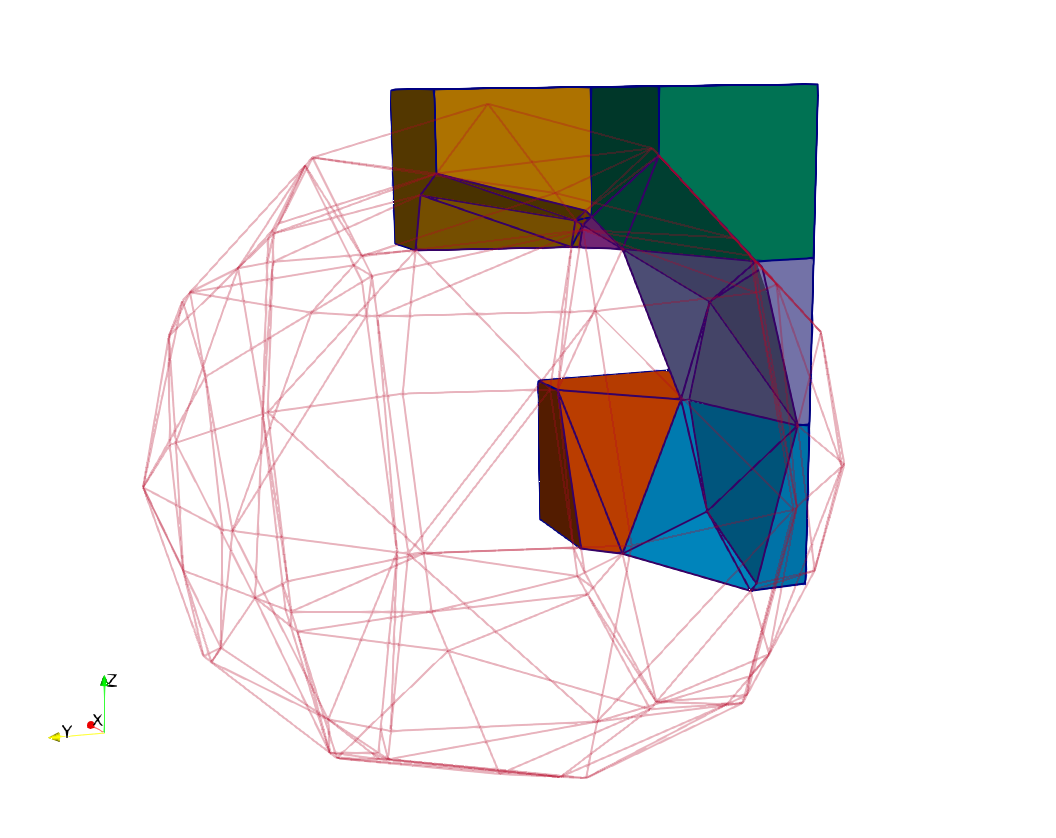}
		\caption{Sample of cells around $\Omega_{\mathrm{stone}}^1$.}
		\label{fig:stony_test:mesh:cell}
	\end{subfigure}
	\caption{TP3: The three-dimensional mesh of the stony soil.}
	\label{fig:stony_test:mesh}
\end{figure}
This last numerical experiment examines root growth in a stony soil on a wider space and time scales with respect to the previous test case. 
%Stones in the soil are impermeable and act as a barrier, preventing root tips from penetrating the rocks. 
%Additionally, a high stone content also affects soil properties such as moisture levels, which in turn impacts root growth \cite{Morandage2021}.

We consider a hexahedral soil sample $\Omega_{\text{soil}} := (0,B)^2 \times (-H, 0)$, with $B=50$~cm and $H=100$~cm, which contains two spherical stones 
\begin{equation*}
\Omega^i_{\text{stone}} := \{ \xx \in \R^3, \norm{\xx - \bm{x}_{\text{stone}}^i} \leq r_i \}
\end{equation*}
with
\begin{eqnarray*}
	&\bm{x}_{\text{stone}}^1=[15, 15, -36.5]~\mathrm{cm}, &\quad r_1=5 ~\mathrm{cm}, \\
	&\bm{x}_{\text{stone}}^2=[25, 31.25, -25]~\mathrm{cm}, &\quad r_2=6~\mathrm{cm}.
\end{eqnarray*}
The spheres are discretized as \textit{UV spheres}, which is a standard shape available in many 3D modeling tools. In particular, we tessellate the surface of the spheres along 8 meridians and 6 parallels each.
Since the stones are assumed to be impervious, we exclude them from the computational domain, i.e. $\Omega = \Omega_{\text{soil}} \setminus \{\Omega^i_{\text{stone}}\}_{i \in \{1, 2\}}$.

Mesh generation on $\Omega$ exploits the possibility of VEM to handle generic polyhedral cells.
We begin by creating a regular 3D structured grid mesh for the $\Omega_{\text{soil}}$ domain, using cubes of edge $3.125$ cm.
Then, since the stones are not involved in the simulation, we subtract them from the mesh by removing the material from the 3D cells.
The result of this process is shown in Figure~\ref{fig:stony_test:mesh:detail}, where part of the global domain mesh is clipped to reveal the internal removed cells, while the color scale refers to the initial soil water pressure head distribution.  As we can observe from the image, the resulting mesh contains complex polyhedral cells, which can be easily handled by the numerical method proposed in this work. Figure~\ref{fig:stony_test:mesh:cell} illustrates a sample of the generated polyhedral cells, highlighting the presence of concave shapes. Using the VEM, such badly-shaped elements can easily be handled, without the need of generating good quality sub-tetrahedrons and without consequently refine the neighboring elements to keep conformity.  %Figure \ref{fig:stony_test:repfunction} shows a clipped view of the soil domain colored by the repulsion function \eqref{eq:repulsionfunction}, which is used to avoid the roots penetrating the stones.

On the top and bottom faces of $\Omega$ we impose Dirichlet boundary conditions for the soil water pressure head $\rph{}$. In particular we set $\rph{}=-500$~cm at the top and $\rph{}=-100$~cm at the bottom. No flux boundary conditions are imposed on all the other faces, on the stone surface and at root tips, while we set a constant transpiration rate of 0.2 cm$^3$/day at the root collar \cite{Jin2020}. 

The initial condition for the soil water pressure head is obtained by solving \eqref{eq:3D1Dproblem1} with $\mathcal{C}(\rph{})=0~ \forall \rph{}$ and $\Lp=0$, i.e. by solving a stationary problem is absence of the RSA. We then simulate root growth in $\Omega$ for 160 days, starting from a seed positioned in 
$\xx_S=[25,25,-5]$ cm. We divide the considered time-span into uniform time intervals of width $\Delta \mathcal{I}_j=1$~day. For the time discretization of the constraint equations we choose $\Delta t=\Delta \mathcal{I}_j$. The PDE parameters used in this simulation are reported in Table \ref{tab_tp3}, whereas the water content $\theta(\rph{})$ and the non-linear coefficients $K(\rph{})$ and $\mathcal{C}(\rph{})$ are defined as in \eqref{eq:theta_exp}, \eqref{eq:K_exp} and \eqref{eq:C_exp}, respectively. With respect to the previous test case, we are choosing lower values of $a$ and $n$. According to \eqref{eq:theta_exp}, this promotes root growth as it enhances water retention at high values of $|\rph{}|$.
\begin{table}
    \renewcommand*{\arraystretch}{1.2}
    \centering
    \begin{tabular}{ccclc}
        \hline
        \textbf{Parameter } & \textbf{Value} &\textbf{Unit} &\textbf{Description} &\textbf{Reference}\\
        \hline
        $a$&0.02&$1/\mathrm{cm}$& Water retention curve shape parameter&\cite{CarselParrish_1988}\\
        $n$&1.2& -& Water retention curve shape parameter&\cite{CarselParrish_1988}\\
        $\wc_r$&0.06&-&Residual volumetric water content&\cite{Jin2020,CarselParrish_1988}\\
        $\wc_s$&0.41&-&Saturated volumetric water content&\cite{Jin2020,CarselParrish_1988}\\
        $K_s$&10.24&cm/day&Saturated hydraulic conductivity&\cite{Jin2020,CarselParrish_1988}\\
        $R$&0.05&cm&Root radius&\cite{javaux2008}\\
        $\kx$&0.18&day/cm&Reciprocal of specific root axial conductance&\cite{Jin2020,schnepf2023}\\
        $\Lp$&$1.36\cdot 10^{-6}$&cm$^{2}$/day&Permeability of root wall &\cite{schnepf2023}\\
        $\sigma_{\mathrm{max}}$&$1$&MPa &Soil strength parameter & \\
        \hline
    \end{tabular}
    \vspace{4pt}
    \caption{TP3: PDE parameters used for simulation.}
    \label{tab_tp3}
\end{table}

The parameters concerning the development of the RSA are instead reported in Table \ref{tab2_tp3}. They are mostly taken from \cite{Morandage2019} and slightly adapted to the proposed growth model. Although most parameters refer to winter wheat, we underline that the primary aim of this numerical example is to show how the proposed 3D-1D VEM approach fits a larger-scale simulation, regardless of the choice of the actual plant genotype. Some of the parameters reported in Table \ref{tab2_tp3} are stochastic and follow either a lognormal or a uniform distribution. In the first case we report the mean value and the standard deviation inside round brackets, whereas in the second we provide the bounds of the uniform distribution inside square brackets. Let us remark that most of the parameters vary with the root order. Those for which we report only a single value or the data of a single distribution are instead assumed to be invariant with respect to root order. 
\begin{table}

\renewcommand*{\arraystretch}{1.2}
    \centering
    \begin{tabular}{cccccl}
        \hline
        \multirow{2}{*}{\textbf{Parameter}} & \multicolumn{3}{c}{\textbf{Value}}   & \multirow{2}{*}{\textbf{Unit}}                           & \multirow{2}{*}{\textbf{Description}}                     \\
        \cline{2-4} & \textbf{0} & \textbf{1} & \textbf{2} &    \\ 
        \hline
        $\vert (\rph{})_1\vert$  &&1 &&cm& Hypoxia threshold\\
        $\vert (\rph{})_2\vert$  &&510& &cm&Lower bound for $\mathrm{Imp}(\rph{})=1$\\
        $\vert (\rph{})_3\vert$  &&920& &cm&Upper bound for $\mathrm{Imp}(\rph{})=1$\\
        $\vert (\rph{})_4\vert$  &&$1.6\cdot 10^{4}$ &&cm&Drought threshold\\
        $\overline{V}_a$&1.2 (0.6) & 1 (0.2) &0.4 (0.12)&cm& Growth rate\\
        $L_B$&0.8 (1.2)& 0.8 (1) & - & cm & Length of basal non-branching zone\\
        $L_A$&4.2 (6.4) &1.8 (2.4) & - &cm&Length of apical non-branching zone\\
        $I$&0.8 (0.4)&1 (0.5) &-&cm&Inter-branch distance\\
        $\alpha_I$&1.4 (0.2)&1.2(0.4)&-&rad&Branching insertion angle\\
        $X$&5&3&-&-&Number of xylem poles\\
        $b_c$&1&1&-&-&Branching probability parameter\\
        $k_g$&&[0.1, 0.2]&&-&Geotropism weight\\
        $k_s$&&[0.35, 0.45]&&-&Hydrotropism weight\\
        $k_s$&2&1&1&-&Exotropism weight\\
        $m_a$&&[0,1]&&-&Space exploration parameter\\
        \hline
    \end{tabular}
    \vspace{4pt}
    \caption{TP3: RSA growth parameters used for simulation. Stochastic parameters follow either a lognormal distribution characterized by $mean~(standard ~deviation)$ or a uniform distribution in $[lower ~bound,~upper ~bound]$ \cite{Morandage2019}.}
    \label{tab2_tp3}
\end{table} 

Unlike the previous test case, we allow for more than a zero-order root sprouting from the seed. At the initial time we consider a single straight root segment $\Lambda^0$, representing the vertical mesocotyl. It is parallel to the $z$-axis and connects the seed with the top surface of the computational soil domain, as shown in Figure \ref{fig:stony_test:mesh:detail}. From a computational standpoint, at the first time step the seed is assumed to be both an element of $\mathcal{P}_{\mathrm{tip}}^0$ and of $\mathcal{P}_{\mathrm{branch}}^0$, with $X=19$ potential zero order root that can emerge and which are treated as if they were lateral branches of the vertical mesocotyl, emerging with probability $p_{br}(0)$ with a certain insertion and radial angle (see Section \ref{sec:branching}). In the following time-steps the seed keeps being an element of $\mathcal{P}_{\mathrm{branch}}^j$, until the maximum number of sprouting roots is saturated. For what concerns root tips, we link their growth rate to root order and root age. In particular we set
\begin{equation*}
    V_a(t_j)=\begin{cases}\overline{V}_a(0) &\text {if } \omega=0\\
\overline{V}_{a}(\omega) e^{-\frac{\overline{V}_{a}(\omega)}{L_{\max}(\omega)} a_p^j} &\text {if } \omega=1,2
    \end{cases}
\end{equation*}
where $a_P^j$ represents the age of the root ending in the tip $P$ at time $t_j$ and $L_\mathrm{max}$ is the maximum root elongation, depending on the order.

\begin{figure}[t]
    \centering
    \begin{subfigure}{0.33\textwidth}
        \centering   \hspace{-0.6cm}\includegraphics[width=1.3\textwidth]{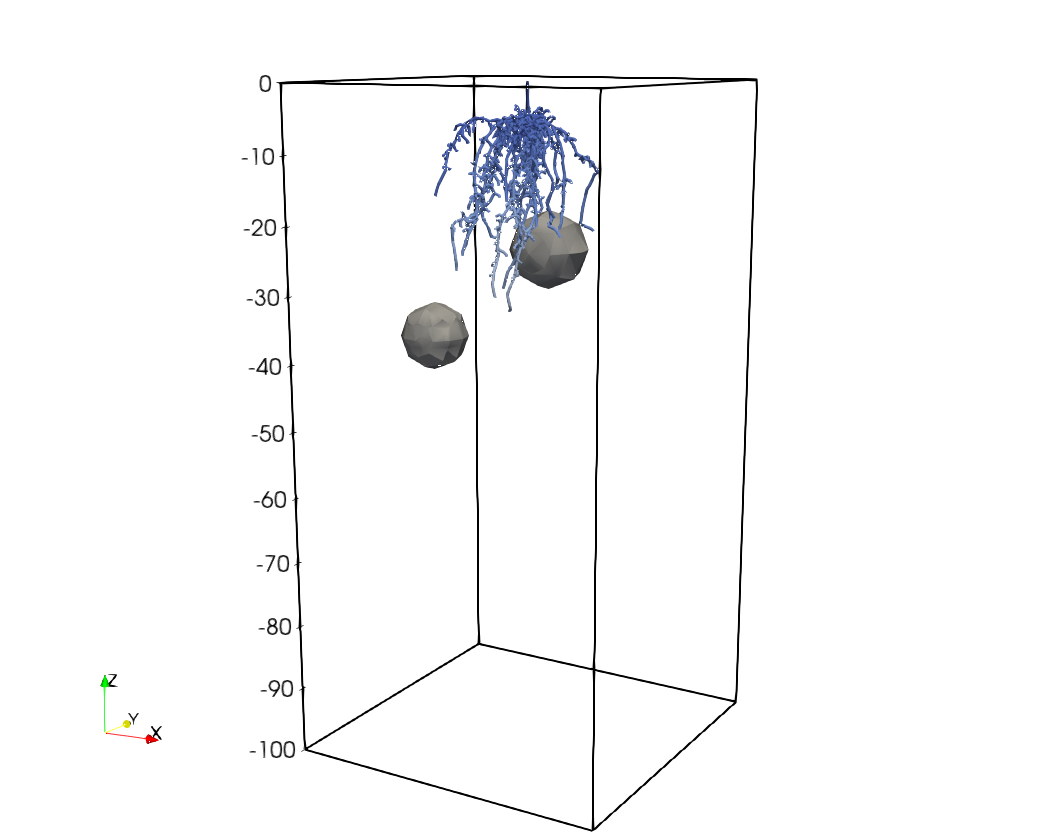}
        \caption{$t=40$ days}
        \label{}
    \end{subfigure}\hspace{-0.6cm}%
    \begin{subfigure}{0.33\textwidth}
        \centering
        \includegraphics[width=1.3\textwidth]{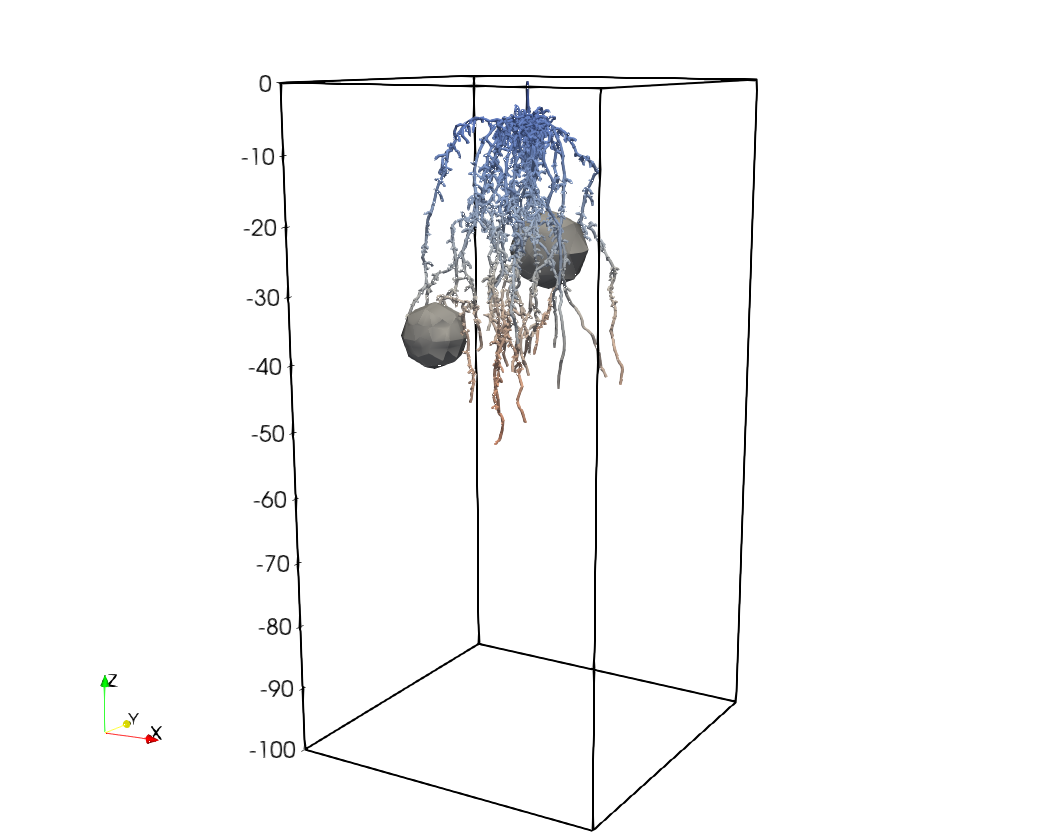}
        \caption{$t=80$ days}
        \label{}
    \end{subfigure}\hspace{-0.6cm}%
    \begin{subfigure}{0.33\textwidth}
        \centering
        \includegraphics[width=1.3\textwidth]{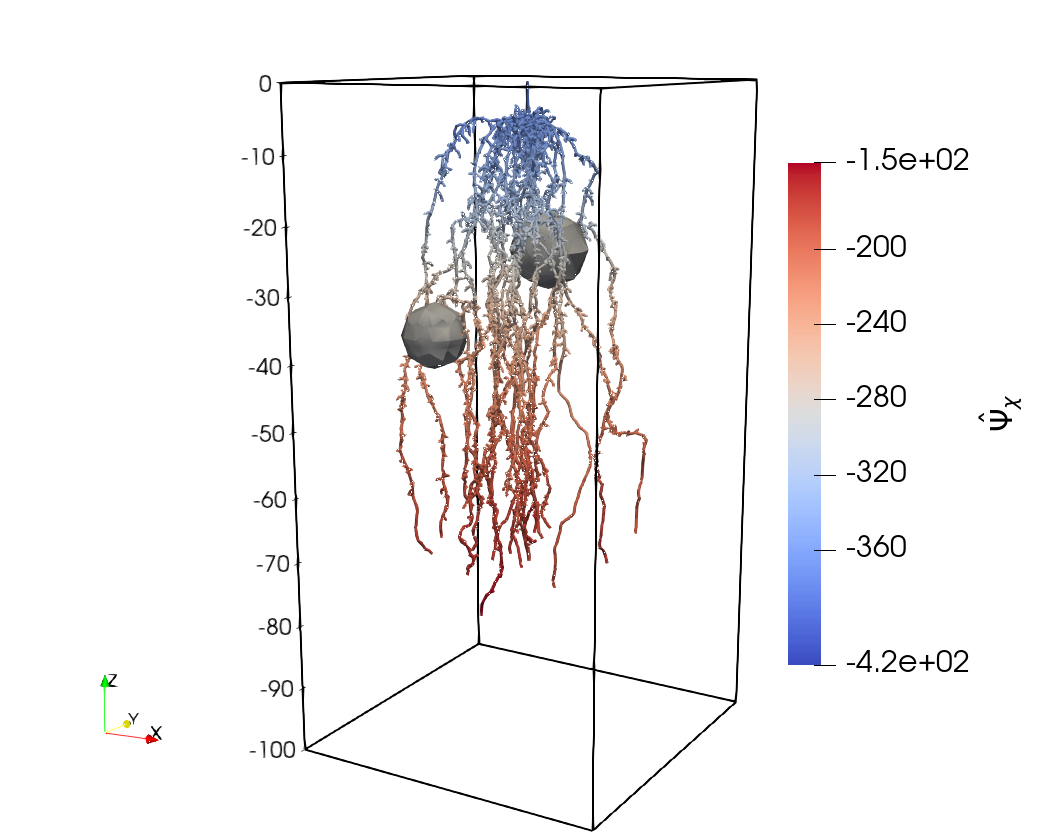}
        \caption{$t=160$ days}
        \label{}
    \end{subfigure}
    \caption{TP3: RSA at three different time instants. Color scale refers to xylem pressure head.}
    \label{fig:Tp3_growth}
\end{figure}

Figure \ref{fig:Tp3_growth} shows the evolution of the RSA at different stages of the simulation, namely after 40, 80 and 160 days, along with the xylem pressure head distribution.
The root wall permeability was kept constant, but as in the previous numerical example it could be linked, for example, to the root age.

The contribution of the repulsion function in the computation of the growth direction allows us to avoid the roots to enter the stones. In particular, for the proposed experiment we defined a function $\repulsionfunction(\bm{x})$ that varies linearly between 1 and 0 in the first layer of elements of a sub-tetrahedralization of $\Th$ around the stones themselves, but the actual repulsion function was then defined as
\begin{equation*}
\tilde{\repulsionfunction}(\bm{x})=\begin{cases}
    \repulsionfunction (\bm{x}) &\text{if } \repulsionfunction(\bm{x})>0.9\\
    0 &\text{otherwise}.
\end{cases}
\end{equation*}
This allows the roots to get very close to the stone surface without the need of refining the mesh or to increase the polynomial order of the repulsion function. The same repulsion function was imposed in the first layer of elements near the top surface of the domain to avoid the roots growing above the soil surface.

\begin{figure}
        \centering
\includegraphics[width=0.4\textwidth]{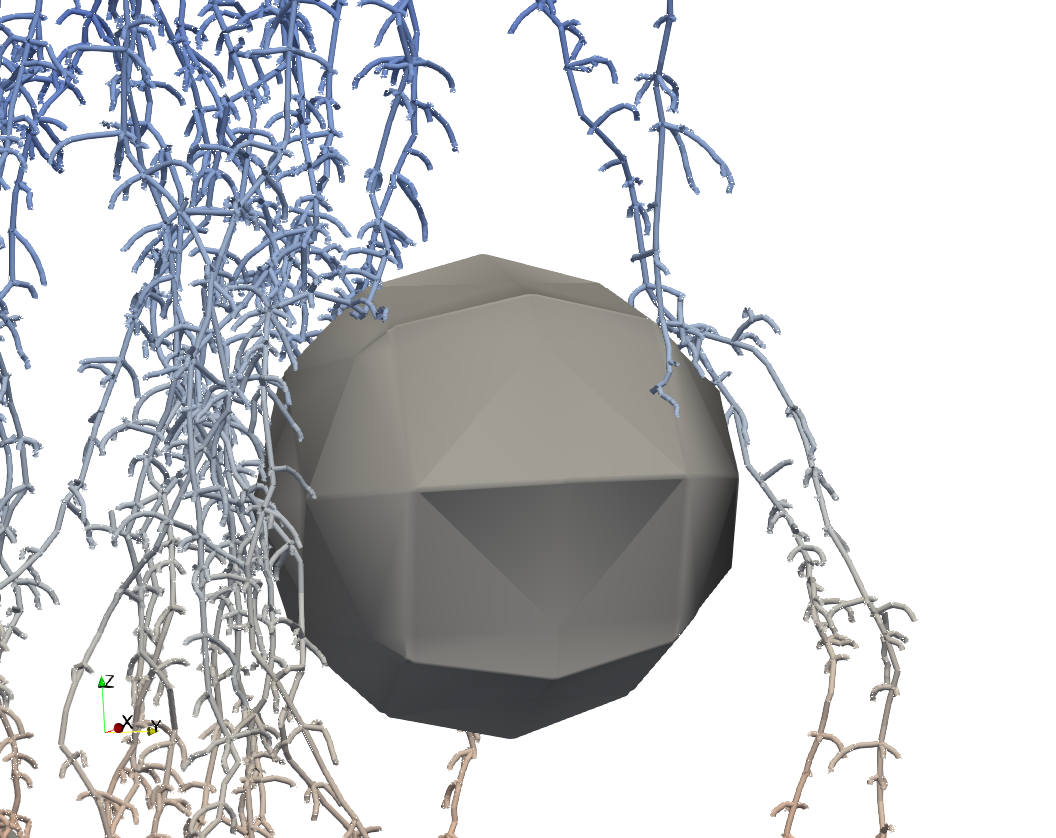}
        \caption{TP3: detail of the RSA around a stone}
        \label{fig:Tp3_stonedetail}
    \end{figure}

\begin{figure}[t]
    \centering
    \begin{subfigure}{1\textwidth}
        \centering
        \begin{subfigure}{0.5\textwidth}
        \includegraphics[width=1.1\textwidth]{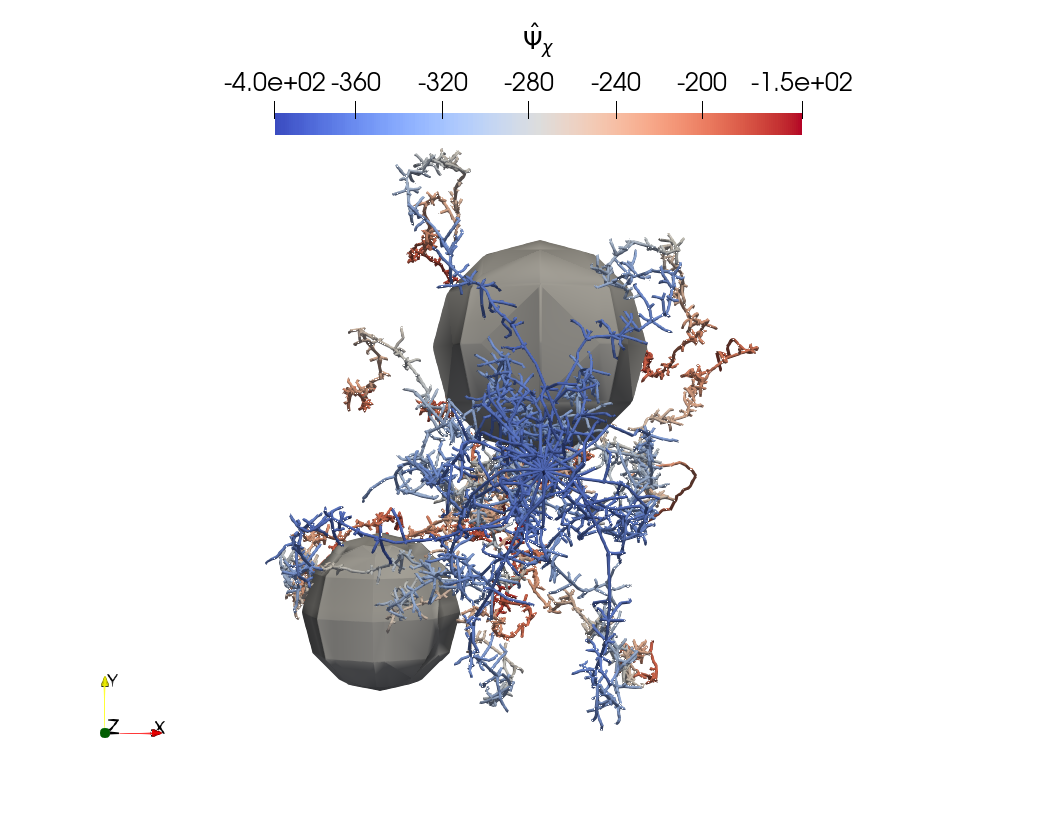}
        \caption{$t=160$ days, view from top face}
        \label{}
    \end{subfigure}\hspace{-0.5cm}%
    \begin{subfigure}{0.5\textwidth}
        \centering
        \includegraphics[width=1.1\textwidth]{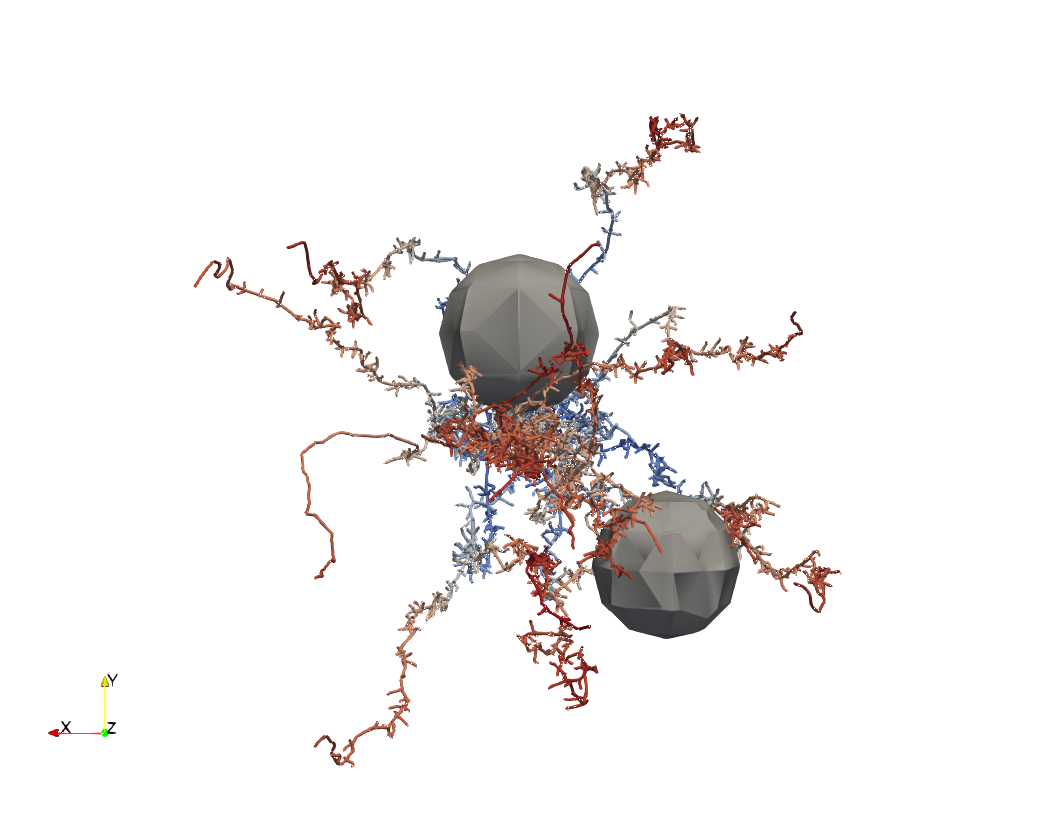}
        \caption{$t=160$ days, view from bottom face}
        \label{}
    \end{subfigure}
    \end{subfigure}
    \caption{TP3: Top and bottom view of the RSA at $t=$ 160 days}
    \label{fig:Tp3_xy}
\end{figure}

A detail on the behavior of the RSA in the surroundings of one of the stones is shown in Figure \ref{fig:Tp3_stonedetail}. Finally, Figure \ref{fig:Tp3_xy} provides an $xy$ plane perspective of the final RSA, from the top and bottom faces of the computational domain.

The number of degrees of freedom involved in this numerical example tends to increase fast, as we are simulating the growth of the RSA for a quite long period of time, allowing for a high number of lateral roots. To enhance the efficiency of the proposed solving strategy we hence adopted a preconditioned conjugate gradient method to solve~\eqref{eq:conjugate_gradient_equation}. Referring to the notation introduced in Section \ref{sec:discr_cost_fun}, we used a very simple preconditioner, defined as
\begin{equation}
\bm{\mathcal{P}}=\begin{bmatrix} \bm{G} &\bm{0}\\\bm{0} &\widehat{\bm{G}}\end{bmatrix}.
\end{equation} Let us remark that the matrices $\bm{G}$ and $\widehat{\bm{G}}$ are also part of the diagonal blocks of matrix $\bm{\mathcal{M}}$ (see \eqref{eq:matrix_cg}). The proposed preconditioner can hence be seen as a very coarse approximation of matrix $\bm{\mathcal{M}}$, retaining only the blocks that do not require the inversion of any matrix. Table \ref{tab3_tp3} shows the number of conjugate gradient iterations required, with and without preconditioning, for each iteration $\ell$ of the non-linear solver during 4 different timesteps $\Delta \mathcal{I}_j$, with $j=1,5,10,15$. We can observe that in both cases, and as observed also in TP1, the number of required CG iterations decreases from the first to the last non-linear iteration. However, in the non-preconditioned case, the number of CG iterations at the first non-linear steps tends to increase very rapidly with the number of degrees of freedom of the control variables ($\Ncs+\Ncx$). On the contrary, using the preconditioner $\bm{\mathcal{P}}$, the number of CG iterations tends to remain much lower and almost unchanged. It is interesting to remark that at the last considered time step, characterized by $\Ncs+\Ncx=422978 $ degrees of freedom for the control variables, the preconditioned algorithm still requires only 15, 11, 2 and 0 iterations of the conjugate gradient scheme respectively for the 4 non-linear steps performed. The non-preconditioned version was not pushed that far. For the sake of completeness, let us specify that the conjugate gradient iterations are stopped, for both the preconditioned and non-preconditioned case, when the norm of the absolute residual $r(\mathcal{X})=\bm{\mathcal{M}}\mathcal{X}+\bm{d}$ goes under a certain threshold. In particular when
\begin{equation*}
\|r(\mathcal{X})\|< 10^{-6}(1+\|r(\mathcal{X}_0)\|).
\end{equation*}
At the beginning of each time step, after the RSA has been updated, the initial guess $\mathcal{X}_0$ on the new root network is computed by averaging the value of the control variables on the old network and by distributing this constant value on the whole new RSA, while at each non-linear step the initial guess is given by the value of the control variables produced at the previous non-linear step.

% \usepackage{booktabs}
% \usepackage{multirow}
% \usepackage{graphicx}
\begin{table}[!ht]
\centering
    \renewcommand*{\arraystretch}{1.3}
    \setlength{\tabcolsep}{5pt}
    \begin{tabular}{c|cccccc|ccccc|ccccc|cccc}
        \hline
        $\bm{j}$                           & \multicolumn{6}{c|}{1} & \multicolumn{5}{c|}{5}  & \multicolumn{5}{c|}{10}  & \multicolumn{4}{c}{15}  \\ 
        $\bm{\Ncs+\Ncx}$              & \multicolumn{6}{c|}{6} & \multicolumn{5}{c|}{36} & \multicolumn{5}{c|}{228} & \multicolumn{4}{c}{954} \\ \cline{1-21}
        $\bm{\ell}$                   & 0  & 1 & 2 & 3 & 4 & 5 & 0    & 1   & 2  & 3 & 4 & 0     & 1    & 2 & 3 & 4 & 0     & 1    & 2   & 3  \\ \cline{1-21}
        \textbf{ \# CG (no prec.)}  & 9  & 8 & 5 & 2 & 2 & 0 & 23   & 15  & 4  & 4 & 0 & 128   & 105  & 8 & 2 & 0 & 897   & 894  & 12  & 2  \\
         \textbf{\# CG} $(\bm{\mathcal{P}})$ & 4  & 4 & 3 & 2 & 0 & 0 & 5    & 4   & 3  & 2 & 1 & 7     & 4    & 3 & 1 & 0 & 8     & 5    & 4   & 0 \\ \hline
    \end{tabular}%
    \vspace{4pt}
    \caption{TP3: Number of iterations of the conjugate gradient scheme required at each non-linear step $\ell$ for 4 different time intervals $\Delta\mathcal{I}_j$. Non-preconditioned (No prec.) and preconditioned ($\bm{\mathcal{P}}$) case. The total number of DOFs for the control variables is also reported for each considered time-step.}
    \label{tab3_tp3}
\end{table}

\section{Conclusions}
In this work we propose the application of an optimization-based 3D-1D coupling approach to the modeling and simulation of the exchanges between a soil sample and a growing root-network. In particular we aim at computing the water pressure head in the soil and in the root-xylem, by coupling Richards equation in primal form in the soil domain with Stokes equation in the root xylem. Thanks to proper assumptions on the regularity of the quantities of interest, the root xylem is identified with the root centerline, leading to a well-posed 3D-1D reduced formulation of the problem. The Virtual Element Method is used to discretize the problem in the bulk soil domain, while Mixed Finite Elements are adopted to discretize the problem in the xylem. The imposition of coupling conditions at the soil-root interface is tackled by means of a PDE-constrained optimization approach, which ends up in a direct computation of interface quantities. The use of the VEM for the 3D discretization extends the applicability of the method to complex shaped computational domains, easing the meshing process in presence of physical obstacles such as pot walls or stones. Root growth is modeled by means of a 
discrete-hybrid tip-tracking strategy, endowed with proper rules for the emergence of lateral roots.
The proposed approach is validated by several numerical examples, showing both the accuracy of the method and its applicability to realistic and large scale simulations.

The discretization of Richards equation by Mixed Virtual Elements will be the subject of a forthcoming work, as well as a deeper investigation on the computational tools which may increase the efficiency of the solving strategy, such as tailored preconditioners and parallelization.

\section*{Acknowledgements}
The author S.B. kindly acknowledges partial financial support provided by PRIN project ``Advanced polyhedral discretisations of heterogeneous PDEs for multiphysics problems'' (No. 20204LN5N5\_003), by PNRR M4C2 project of CN00000013 National Centre for HPC, Big Data and Quantum Computing (HPC) (CUP: E13C22000990001) and the funding by the European Union through project Next Generation EU, M4C2, PRIN 2022 PNRR project P2022BH5CB\_001 ``Polyhedral Galerkin methods for engineering applications to improve disaster risk forecast and management: stabilization-free operator-preserving methods and optimal stabilization methods.''. The author S.F. kindly acknowledges partial financial support provided by Project NODES, which received funding from the MUR-M4C2 1.5 of the National Recovery and Resilience Plan (PNRR) with grant agreement no. ECS00000036. The author D.G. kindly acknowledges to be holder of a Research Grant from INdAM (Istituto Nazionale di Alta Matematica) at the research unit of Politecnico di Torino. The author G.T. kindly acknowledges financial support provided by the MIUR programme ``Programma Operativo Nazionale Ricerca e Innovazione 2014 - 2020''~(CUP: E11B21006490005). The authors G.T. and F.V. kindly acknowledge financial support provided by INdAM - GNCS Project CUP\_E53C23001670001. 


\bibliographystyle{IEEEtran}
\bibliography{biblio.bib}

\end{document}